%% file: main_XVEM_elasticity.tex
\documentclass[final,12pt,times]{elsarticle}

\usepackage[a4paper, total={6.3in, 9.05in}]{geometry}

\makeatletter
\def\ps@pprintTitle{%
 \let\@oddhead\@empty
 \let\@evenhead\@empty
 \def\@oddfoot{}%
 \let\@evenfoot\@oddfoot}
\makeatother

\usepackage{graphicx}
\usepackage{mathrsfs}
\usepackage{multirow}
\usepackage{comment}

\input{defsVEM_new.tex}

\graphicspath{ {./Figures/}}

\begin{document}
\begin{frontmatter}
\title{Extended virtual element method for two-dimensional linear elastic fracture}

  \cortext[cor1]{Corresponding author}

  % Authors' names
  % ----------------------------
  \author[UNIFE]  {E. Benvenuti}
  \ead{elena.benvenuti@unife.it}

  \author[UNIFE]  {A. Chiozzi\corref{cor1}}
  \ead{andrea.chiozzi@unife.it}

  \author[IMATI]   {G. Manzini}
  \ead{marco.manzini@imati.cnr.it}
  
  \author[UCDAVIS]{N. Sukumar}
  \ead{nsukumar@ucdavis.edu}

  % Affiliations
  % ----------------------------
  \address[UNIFE]{
    Department of Engineering,
    University of Ferrara, Via Saragat 1, 44122 Ferrara, Italy
  }
  \address[IMATI]{
  Istituto di Matematica Applicata e Tecnologie Informatiche, 
  Consiglio Nazionale delle Ricerche, Pavia, Italy,
  }
  \address[UCDAVIS]{
    Department of Civil and Environmental Engineering,
    University of California, Davis, CA 95616, USA
  }

  % Abstract
  % ----------------------------
  \begin{abstract}  
 
    In this paper, we propose an eXtended Virtual Element Method
    (X-VEM) for two-dimensional linear elastic fracture. 
    This approach, which is an extension of the standard Virtual Element Method (VEM), facilitates mesh-independent modeling of crack discontinuities and elastic crack-tip singularities on general polygonal meshes.
    For elastic fracture in the X-VEM, the standard virtual element
    space is augmented by additional basis
    functions that are constructed by multiplying standard virtual basis functions by
    suitable enrichment fields, such as asymptotic mixed-mode crack-tip solutions.
    The design of the X-VEM requires an extended projector that maps
    functions lying in the extended virtual element space onto a set
    spanned by linear polynomials and the enrichment fields.
    An efficient scheme to compute the mixed-mode
    stress intensity factors using the domain form of the interaction
    integral is described. The formulation permits integration of 
    weakly singular  functions to be performed over the boundary edges of the element.
    Numerical experiments are conducted on benchmark mixed-mode linear elastic fracture problems that demonstrate the
    sound accuracy and optimal convergence in energy of the proposed
    formulation.
    
  \end{abstract}
  
  % Keywords
  % ----------------------------

\end{frontmatter}

%% paper

\input{sec1_intro.tex}

\input{sec2_model_new.tex}
\input{sec3_xvem_new.tex}

\input{sec4_imple_new.tex}
\input{sec5_numer.tex}
\input{sec6_final.tex}
\section*{Acknowledgements}
Elena Benvenuti and Andrea Chiozzi gratefully acknowledge the support
of PRIN, Italy: Progetti di Ricerca di Rilevante Interesse Nazionale
(Bando 2015) Prot.2015LYYXA8. Andrea Chiozzi acknowledges the support
of the research fund FIR 2020 of the University of Ferrara,
Italy. 
Gianmarco Manzini gratefully acknowledges the financial support of the ERC Project CHANGE, which has received funding from the European Research Council under the European Union’s Horizon 2020 research and innovation program (grant agreement no.~694515).

% Bibliography
% ----------------------------
%\bibliographystyle{elsarticle-num}

\clearpage

\bibliographystyle{unsrt}
\biboptions{sort&compress}
\bibliography{main_XVEM_elasticity}

\end{document}

%% file: defsVEM_new.tex
\usepackage{amsmath,amsthm,amsfonts}  %   AMS equation environments
\allowdisplaybreaks

\usepackage{bm}
\usepackage{bbm}

\usepackage{booktabs}

\usepackage{color}
\usepackage{xcolor}

\usepackage{graphicx}
\usepackage{caption}
\usepackage[font=normalsize]{subcaption}
\usepackage{nicefrac}

\newtheorem{remark}{Remark}[section]

\definecolor{MyDarkGreen}{rgb}{0,0.45,0}

\usepackage[normalem]{ulem}
\newcounter{corr}
\definecolor{violet}{rgb}{0.580,0.,0.827}
\newcommand{\corr}[3]{\typeout{Warning : a correction remains in page
\thepage}				
\stepcounter{corr}
{\color{blue}\ifmmode\text{\,\sout{\ensuremath{#1}}\,}\else\sout{#1}\fi}
{\color{red}#2}
{\color{violet} #3}}

\newcommand{\VEM}{\text{VEM}}

\newcommand{\XVEM}{\text{X-VEM}}

\newcommand{\INTP}{I} %% {\footnotesize{\texttt{I}}}
\newcommand{\SPAN}{\mbox{\textrm{span}}}
\newcommand{\REAL}{\mathbbm{R}}

\newcommand{\gv}{\bm{g}}

\newcommand{\mv}{\bm{m}}
\newcommand{\nv}{\bm{n}}

\newcommand{\qv}{\bm{q}}
\newcommand{\rv}{\bm{r}}

\newcommand{\tv}{\bm{t}}
\newcommand{\uv}{\bm{u}}
\newcommand{\vv}{\bm{v}}
\newcommand{\wv}{\bm{w}}
\newcommand{\xv}{\bm{x}}
\newcommand{\yv}{\bm{y}}

\newcommand{\Bv}{\bm{B}}
\newcommand{\Cv}{\bm{C}}
\newcommand{\Dv}{\bm{D}}
\newcommand{\Ev}{\bm{E}}

\newcommand{\Gv}{\bm{G}}

\newcommand{\Iv}{\bm{I}}
\newcommand{\Jv}{\bm{J}}
\newcommand{\Kv}{\bm{K}}

\newcommand{\Mv}{\bm{M}}
\newcommand{\Nv}{\bm{N}}

\newcommand{\Sv}{\bm{S}}

\newcommand{\Vv}{\bm{V}}

\newcommand{\qvX}{\qv_{X}}

\newcommand{\Vvh} {\Vv^{\hh}}
\newcommand{\Vvha}{\Vv^{\hh,\ast}}
\newcommand{\Vvhg}{\Vv^{\hh}_{\gv}}
\newcommand{\Vvhz}{\Vv^{\hh}_{\bm{0}}}

\newcommand{\VvhX}{\Vv_{X}^{\hh}}
\newcommand{\VvhXg}{\Vv_{X,\gv}^{\hh}}
\newcommand{\VvhXz}{\Vv_{X,\bm{0}}^{\hh}}

\newcommand{\Vvhat}{\widetilde{\Vv}^{\hh,\ast}}
\newcommand{\VvhXt}{\widetilde{\Vv}_{X}^{\hh}}

\newcommand{\Vvhp} {\Vv^{\hh,+}}
\newcommand{\Vvhm} {\Vv^{\hh,-}}

\newcommand{\NvX}{\Nv_{X}}
\newcommand{\MvX}{\Mv_{X}}

\newcommand{\GvX}{\Gv_{X}}
\newcommand{\BvX}{\Bv_{X}}
\newcommand{\DvX}{\Dv_{X}}
\newcommand{\SvX}{\Sv_{X}}
\newcommand{\EvX}{\Ev_{X}}

\newcommand{\TILDE}[1]{\widehat{#1}}
\newcommand{\GvXt}{\TILDE{\Gv}_{X}}
\newcommand{\BvXt}{\TILDE{\Bv}_{X}}
\newcommand{\DvXt}{\widehat{\Dv}_{X}}

\newcommand{\matPiaX} {\bm{\Pi}^a_{X}}

\newcommand{\KvPX}{\Kv^{\P}_{X}}
\newcommand{\KvPXt}{\widetilde{\Kv}^{\P}_{X}}
\newcommand{\KvPXc}{\Kv^{\P}_{X,c}}
\newcommand{\KvPXs}{\Kv^{\P}_{X,s}}

\newcommand{\as}{a}
\newcommand{\bs}{b}

\newcommand{\is}{i}
\newcommand{\js}{j}

\newcommand{\rs}{r}

\newcommand{\us}{u}
\newcommand{\vs}{v}

\newcommand{\xs}{x}
\newcommand{\ys}{y}

\newcommand{\Bs}{B}

\newcommand{\Gs}{G}

\newcommand{\Vs}{V}

\newcommand{\Vsh}{\Vs^{\hh}}

\newcommand{\GsX}{\Gs_{X}}
\newcommand{\BsX}{\Bs_{X}}

\newcommand{\GsXt}{\TILDE{\Gs}_{X}}
\newcommand{\BsXt}{\TILDE{\Bs}_{X}}

\newcommand{\vshXi}{\vs^{\hh}_{X,i}}

\newcommand{\calE}{\mathcal{E}}

\newcommand{\calI}{\mathcal{I}}

\newcommand{\kP}{k_{\P}}
\newcommand{\extkP}{\calE_{\kP}}

\newcommand{\PS} [1]{\mathbbm{P}^{#1}}
\newcommand{\PSX}[1]{\mathbbm{P}^{#1}_X}

\newcommand{\HONE}  {H^1}

\newcommand{\HS}[1] {H^{#1}}

\newcommand{\CS}[1] {C^{#1}}

\newcommand{\VS}[1] {V^{#1}}

\renewcommand{\P} {{E}}            % polyhedral element
\newcommand  {\E} {{e}}            % edge
\newcommand{\hh}{h}
\newcommand{\Th}{\Omega_{\hh}}

\newcommand{\xvP}{\xv_{\P}}        % element
\newcommand{\hP}{\hh_{\P}}

\newcommand{\hE}{\hh_{\E}}

\newcommand{\mP}{\ABS{\P}}

\newcommand{\NMB}{N}
\newcommand{\NP}{\NMB_{\P}}         % number of elemental edge

\newcommand{\Ndofs}{N_{\textrm{dofs}}}

\newcommand{\dx}{\,d\xv}

\newcommand{\norP} {\bm{n}_{\P}}

\newcommand{\usxI} {\us_{\xs}^{\INTP}}
\newcommand{\usyI} {\us_{\ys}^{\INTP}}
\newcommand{\usxII}{\us_{\xs}^{\INTP\INTP}}
\newcommand{\usyII}{\us_{\ys}^{\INTP\INTP}}

\newcommand{\HACECK}[1]{\mbox{\it\v{#1}}}
\newcommand{\usxIt} {\HACECK{\us}_{\xs}^{\INTP}}
\newcommand{\usyIt} {\HACECK{\us}_{\ys}^{\INTP}}
\newcommand{\usxIIt}{\HACECK{\us}_{\xs}^{\INTP\INTP}}
\newcommand{\usyIIt}{\HACECK{\us}_{\ys}^{\INTP\INTP}}

\newcommand{\vsh} {\vs^{\hh}}
\newcommand{\vshx}{\vs^{\hh}_{x}}
\newcommand{\vshy}{\vs^{\hh}_{y}}

\newcommand{\uvh}{\uv^{\hh}}
\newcommand{\uvhX}{\uv^{\hh}_{X}}

\newcommand{\uvI} {\uv^{\INTP}}
\newcommand{\uvII}{\uv^{\INTP\INTP}}

\newcommand{\uvIt} {\mbox{\it\textbf{\v{u}}}^{\INTP}}
\newcommand{\uvIIt}{\mbox{\it\textbf{\v{u}}}^{\INTP\INTP}}

\newcommand{\vvh}{\vv^{\hh}}
\newcommand{\vvhm}{\vv^{\hh,-}}
\newcommand{\vvhp}{\vv^{\hh,+}}
\newcommand{\vvhX}{\vv^{\hh}_{X}}
\newcommand{\vvhXt}{\widetilde{\vv}^{\hh}_{X}}

\newcommand{\wvh}{\wv^{\hh}}
\newcommand{\wvhX}{\wv^{\hh}_{X}}
\newcommand{\wvhXt}{\widetilde{\wv}^{\hh}_{X}}

\newcommand{\psivI} {\bm{\psi}^{\INTP}}
\newcommand{\psivII}{\bm{\psi}^{\INTP\INTP}}

\newcommand{\asP}{\as^{\P}}
\newcommand{\bsP}{\bs^{\P}}

\newcommand{\ash}{\as^{\hh}}
\newcommand{\bsh}{\bs^{\hh}}

\newcommand{\ashX}{\as^{\hh}_{X}}
\newcommand{\bshX}{\bs^{\hh}_{X}}

\newcommand{\ashPX}{\as^{\hh,\P}_{X}}

\newcommand{\ashPXt}{\widetilde{\as}^{\hh,\P}_{X}}

\newcommand{\ashP}{\as^{\hh,\P}}
\newcommand{\bshP}{\bs^{\hh,\P}}

\newcommand{\SP} {S^{\P}}
\newcommand{\SPX}{S^{\P}_{X}}

\newcommand{\abs}    [1]{|#1|}

\newcommand{\ABS}    [1]{\left|#1\right|}

\newcommand{\Vh}  {\VS{\hh}}

\newcommand{\Pia}{\Pi^{a}}
\newcommand{\PiaX}{\Pi^{a}_{X}}

\newcommand{\restrict}[2]{{#1}_{|{#2}}} % not using this !

\newcommand{\EOD}{\end{document}}

\newcommand{\xsP}{\xs_{\P}}
\newcommand{\ysP}{\ys_{\P}}

\newcommand{\x}{\bm{x}}

\newcommand{\sigmav}{\bm{\sigma}}
\newcommand{\varepsilonv}{\bm{\varepsilon}}
\newcommand{\phiv} {\bm{\varphi}}
\newcommand{\phiva}{\bm{\varphi}^{\ast}}
\newcommand{\EQUIV}{\equiv}
\newcommand{\piv} {\bm{\pi}}

\newcommand{\DOFS}[1]{\mbox{\textbf{dofs}}\big(#1\big)}

%% file: sec1_intro.tex
\section{Introduction}
\label{sec:intro}

Over the past two decades, significant attention has been devoted to the
development of numerical techniques to solve problems that admit
singular or discontinuous solutions such as fracture propagation in
solids.
Among these techniques,  enriched  finite  element  approximations based on the partition-of-unity framework~\cite{Melenk:1996:PUF,Babuska:1997:PUM} have 
received considerable attention. The eXtended Finite Element Method (X-FEM)~\cite{Moes:1999:FEM} is one
of the most successful methods to analyse fracture
problems on unstructured triangular and quadrilateral
meshes without requiring remeshing.
For fracture simulations on polygonal meshes, extended finite element formulations have been proposed~\cite{Tabarraei:2008:EFE,Zamani:2011:EIP} as well as the
scaled boundary element method~\cite{Song:1997:CMAME, Song:2002:CS,Song:2018:RSB}.
However, construction of shape functions that are defined on general
polygons renders extended finite element formulations to be more involved and
numerical integration of regular and weakly singular functions over polygons is also an issue that requires special
attention~\cite{Mousavi:2010:DT,Chin:2015:NIH,chin:2017}.

The Virtual Element Method (VEM)~\cite{BeiraodaVeiga-Brezzi-Cangiani-Manzini-Marini-Russo:2013} 
is a stabilized Galerkin formulation to solve  partial differential equations on very general polygonal
meshes that overcomes the many difficulties and challenges that are associated with polygonal
finite element formulations.
The VEM derives from the mimetic finite difference
method~\cite{Lipnikov:2014,BeiraodaVeiga-Lipnikov-Manzini:2014} and is
a generalization of the Finite Element Method (FEM) in which the explicit knowledge of the
basis functions is not needed.
Such basis functions are defined as the solution of a local
elliptic partial differential equation, and are never explicitly
computed in the implementation of the method.
Indeed, the VEM uses the elliptic projections of the basis functions
onto suitable polynomial spaces to discretize the bilinear form and
the continuous linear functional deriving from the variational
formulation.
Such projections are computable because of a careful choice of the
degrees of freedom.
The discretized bilinear form is conveniently decomposed as the sum of
a consistent term, which ensures polynomial consistency, and a
correction term, which guarantees stability.
Moreover, the VEM requires the same element-wise assembly procedure of
the FEM for the construction of the global stiffness matrix, thus
resulting in a linear system of equations from which the solution is
obtained. 

In recent years, the VEM has also been used 
to solve problems in solid mechanics, such as
%%$
two- and three-dimensional linear
elasticity~\cite{Beirao-Brezzi-Marini:2013,Paulino:2014:CMAME},
nearly incompressible elasticity~\cite{Paulino:2017:CMAME,
  Wriggers:2017:CM, Paulino:2021:MECC}, inelastic problems~\cite{Beirao-Lovadina-Mora:2015, Hudobivnik:2019:CM}, mixed variational formulations for linear elasticity~\cite{Artioli:2017:CMAME, Dassi:2020}, linear elasticity on curvilinear elements~\cite{Artioli:2020:CMAME}, and elastodynamics~\cite{Paulino:2019:CMAME, Paulino:2020:IJNME, Antonietti:2021:IJNME}.
However,
very few studies have exploited the flexibility of the method to
deal with meshes that are cut by discontinuities and/or contain interior
singularities.
Among these we mention the virtual element modeling of flow in
fracture networks~\cite{Benedetto-Berrone-Pieraccini-Scialo:2014}
and the application of the VEM to 2D elastic fracture
problems~\cite{Nguyen:2018:VEM, Hussein:2019, Artioli:2020}.
In these studies, hanging nodes are inserted at locations where each
discontinuity intersects an element, so that each cut element is
partitioned into a collection of polygonal elements. 

Approximating spaces that consist of the product of low-order virtual element basis functions and a nonpolynomial function were first
proposed in~\cite{Perugia-Pietra-Russo:2016} for the Helmholtz problem,
where the nonpolynomial function is chosen as a planewave in the
two directions. More recently, drawing inspiration from the X-FEM,  an eXtended
Virtual Element Method (X-VEM) is presented
in~\cite{Benvenuti:2019:CMAME} to treat singularities and crack
discontinuities in the scalar Laplace problem, which also
governs the deformation of a stretched membrane or torsion
in a prismatic beam~\cite{Chiozzi:2020:MECC}. 
An enriched nonconforming virtual element method is proposed in~\cite{Artioli-Mascotto:2021:CMAME}, where 
the approximation spaces is enriched with special singular functions (without using the partition-of-unity framework) to solve the Poisson problem with reentrant corners.

In this paper, we develop an extended virtual element formulation for
linear elastic fracture problems, in which the displacement field
features both discontinuities and crack-tip singularities.
For the X-VEM, we construct an enriched virtual element
space by introducing an additional set of virtual basis functions,
which are built on vectorial enrichment fields that are suitably chosen so that they
reproduce the nature of the  weak singularity in the neighborhood of the crack tip.
Hence, additional information about the exact solution is incorporated
in the computational method, mitigating the effects of the singularity
on the numerical accuracy.
In principle, any number of auxiliary fields can be considered to
enrich the virtual element space.
In the X-FEM, near-tip crack functions are used as enrichment functions in
the discrete space~\cite{Moes:1999:FEM}, whereas in the X-VEM we require the enriched stress fields to be divergence-free and hence choose the asymptotic mode I
and mode II crack-tip displacement solutions as vectorial enrichments.
The use of vectorial enrichments was first proposed in the generalized finite element method~\cite{Duarte:2000:CS}. 
Furthermore, as introduced in Benvenuti et al.~\cite{Benvenuti:2019:CMAME}, discontinuities in the displacement field are
incorporated in the virtual element space using the approach proposed for finite
elements by Hansbo and Hansbo~\cite{Hansbo:2004:FEM}.
In contrast to the X-FEM, the X-VEM for elastic fracture provides 
greater
flexibility since it is applicable to arbitrary (simple and nonsimple) polygonal meshes. Furthermore, unlike the X-FEM where special integration schemes~\cite{Chin:2015:NIH} are needed to accurately evaluate the weak form (domain) integrals, in the X-VEM a one-dimensional quadrature rule on the boundary of the polygonal element suffices to compute such integrals.
As in the VEM, the explicit knowledge of virtual shape functions on general polygons is not required, and as we will detail, in this particular instance of the X-VEM, weak form integrals are computed only on the boundary of the element, where the virtual shape functions are known.

%% OUTLINE
The remainder of this article is organized as follows.
In Section~\ref{sec:model:problem}, we introduce the strong and weak
forms for two-dimensional linear elastic fracture problems.
In Section~\ref{sec:XVEM:formulation}, we describe the extended
virtual element formulation.
For crack tip singularities, we devise an extended projector that maps
functions that lie in the extended virtual element space onto the
space spanned by the basis of linear polynomials and the enrichment
fields.
The approach of Hansbo and Hansbo~\cite{Hansbo:2004:FEM} is used to model crack discontinuities in the X-VEM. The implementation of the X-VEM is discussed in Section~\ref{sec:implementation}. In Section~\ref{sec5:numerical_examples}, we presents results for the discontinuous and extended patch tests, and show that the method delivers optimal rate of convergence in energy for benchmark mixed-mode crack problems.

Final remarks and suggestions for future work are discussed in
Section~\ref{conclusions}.

%% file: sec2_model_new.tex
% Hey Emacs, this is -*-latex-*-

\section{Governing equations for 2D linear elasticity}
\label{sec:model:problem}
\begin{figure}[!htb]
  \centering
  \includegraphics[width=0.4\textwidth]{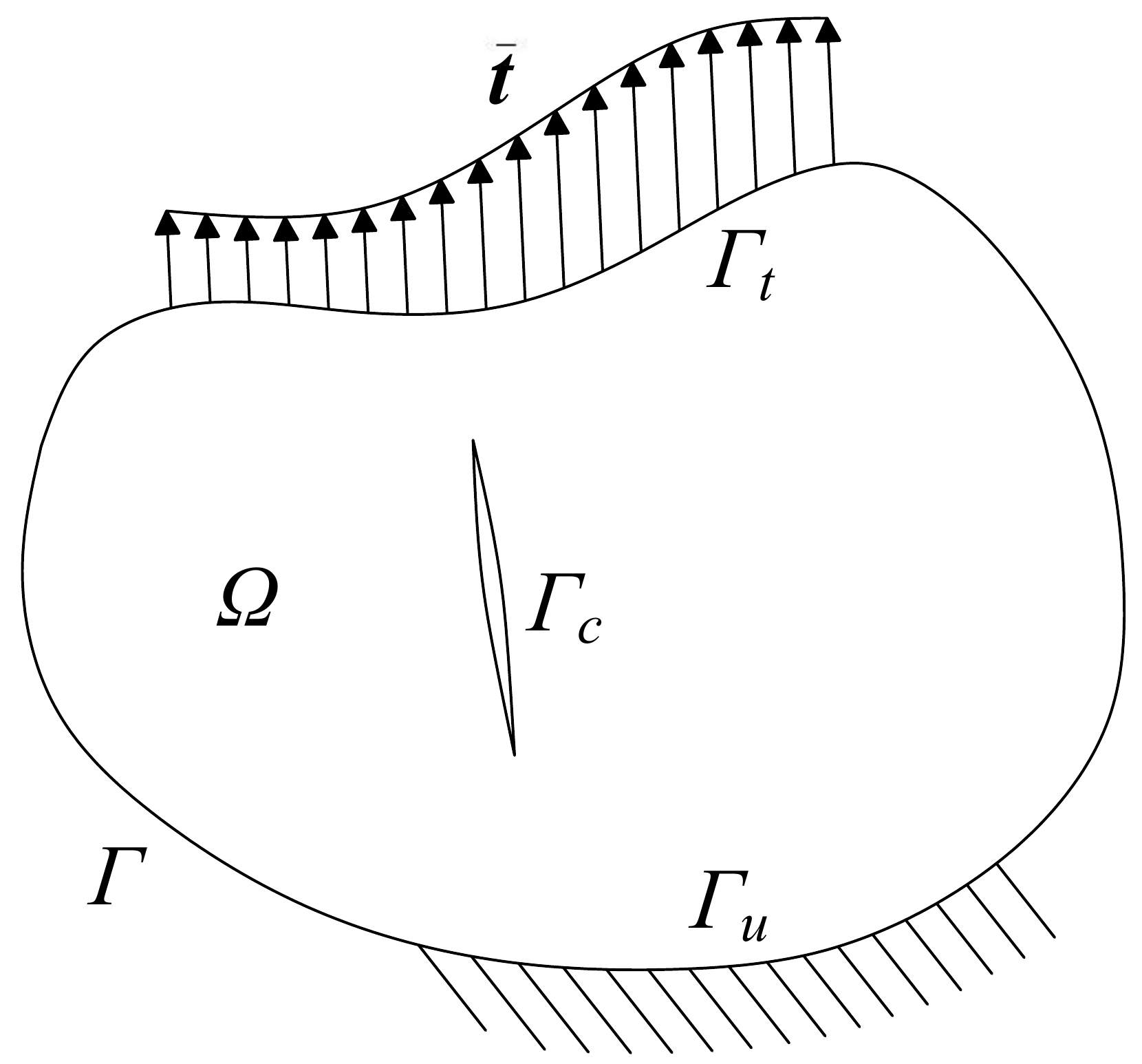}
  \caption{Elastostatic boundary-value problem for an embedded crack.}\label{fig:elastic_body}
\end{figure}

We consider a linear elastic body occupying the two-dimensional domain
$\Omega\subset\mathbb{R}^2$, bounded by $\Gamma$ (see Fig.~\ref{fig:elastic_body}).
We denote the displacement field on $\Omega$ by $\bm{u}(\bm{x})$ and
assume small strains and displacements.
The boundary 
$\Gamma = \Gamma_u \cup \Gamma_t \cup \Gamma_c$, where
$\Gamma_u$, $\Gamma_t$ and $\Gamma_c$ are nonoverlapping, i.e.,
$\Gamma_u\cap\Gamma_t\cap\Gamma_c=\varnothing$. 

Prescribed displacements $\bm{g} \in C^0(\Gamma_u)$ are imposed on
$\Gamma_u$, whereas tractions $\bar{\bm{t}}\in C^0(\Gamma_t)$ are imposed on $\Gamma_t$.
Here, $\Gamma_c$ represents a traction-free internal crack.

We now summarize the governing equations of the elastic problem under
the assumptions of small strains and displacements. Let $\bm{\sigma}$
be the Cauchy stress tensor.
In the absence of body forces, the equilibrium equations are
\begin{subequations}
  \begin{align}
    \nabla\cdot\sigmav  &= \bm{0}\phantom{\bar{\tv}\gv} \ \ \text{in}\;\Omega,   \label{eq:equilibrium:a}\\
    \intertext{with the natural boundary conditions}
    \sigmav\cdot\nv     &= \bar{\tv}\phantom{\bm{0}\gv} \ \ \text{on}\;\Gamma_t, \label{eq:equilibrium:b}\\
    \sigmav\cdot\nv     &= \bm{0}\phantom{\bar{\tv\gv}} \ \ \text{on}\;\Gamma_c, \label{eq:equilibrium:c}\\
    \intertext{where $\nv$ is the unit outward normal, and the essential boundary condition}
    \uv                 &= \bm{g}\phantom{\bm{0}\bar{\tv}} \ \ \text{on}\;\Gamma_u. \label{eq:essentialBC}\\
    \intertext{The small strain tensor $\bm{\varepsilon}$ is related to the displacement field $\bm{u}$ by the compatibility equation}			
    \varepsilonv(\uv)   &= \nabla_s\uv,\label{eq:compatibility} \\
    \intertext{where $\nabla_s$ is the symmetric part of the gradient operator, which is defined
      as
    }
    \nabla_s(\cdot)     &= \frac{1}{2}\left( \nabla(\cdot) + \nabla^T(\cdot) \right).\nonumber \label{eq:grad_sym}
    \intertext{Lastly, the isotropic linear elastic constitutive law is}
    \sigmav(\uv)        &= \Cv : \varepsilonv({\uv}),
  \end{align}
\end{subequations}
where $\Cv$ is the fourth-order elasticity tensor for a homogeneous isotropic material. 

The weak form of the problem is constructed by defining the space of
admissible displacement fields as
\begin{equation}\label{eq:adm_displ}
  \mathscr{U} = \Big\{\vv \in [H^1(\Omega)]^2 : \vv = \bm{g} \:\text{on}\:\Gamma_u,\:\vv\:\text{discontinuous on}\:\Gamma_c \Big\},
\end{equation}
where the space $\mathscr{V}$ is related to the regularity of the
solution, and admits discontinuous functions across the crack.
Similarly, the test function space is defined as:
\begin{equation}\label{eq:test_func}
  \mathscr{U}_0 = \Big\{ \vv \in [H^1(\Omega)]^2 : \vv=0 \:\text{on}\:\Gamma_u,\:\vv\:\text{discontinuous on}\:\Gamma_c \Big\}.
\end{equation}
The weak form of the equilibrium equation reads as: Find
$\uv\in\mathscr{U}$ such that
\begin{equation}
  \label{eq:weak_form}
  \as( \uv, \vv )
  := \int_{\Omega} \sigmav(\vv) : \varepsilonv(\uv) \dx
  =  \int_{\Gamma_t} \bar{\tv}\cdot\vv d\Gamma
  =: \bs(\vv) \quad \forall \vv\in\mathscr{U}_0.
\end{equation}
The above statement is equivalent to the strong
form~\eqref{eq:equilibrium:a} and in a finite element framework it is
solved approximately on a sequence of appropriately nested
finite-dimensional subspaces of $\mathscr{U}$.
 

%% file: sec3_xvem_new.tex
% Hey Emacs, this is -*-latex-*-

\section{Extended virtual element formulation}
\label{sec:XVEM:formulation}

We now discuss the formulation of the extended virtual element method
for two-dimensional elasticity problems. We start, in
Section~\ref{subsec:mesh:regularity:assumptions}, from the definition
and regularity properties of the mesh families for the \XVEM{}, and
after reviewing the `nonenriched' VEM in
Section~\ref{subsec:conforming:VEM:regular}, we provide the design of
the \XVEM{} for full and partial local enrichments in
Sections~\ref{subsec:conforming:VEM:extended}
and~\ref{subsec:partial-enrichment}. 

\subsection{Mesh definition and regularity assumptions}
\label{subsec:mesh:regularity:assumptions}
Let $\mathcal{T}=\{\Th\}_{\hh}$ be a family of decompositions of
$\Omega$ into nonoverlapping polygonal elements $\P$ with
nonintersecting boundary $\partial\P$, barycenter
$\xvP\EQUIV(\xs_{\P},\ys_{\P})^T$, area $\mP$, and diameter
$\hP=\sup_{\xv,\yv\in\P}\abs{\xv-\yv}$.
The subindex $\hh$ that labels each mesh $\Th$ is the maximum of the
diameters $\hP$ of the elements of that mesh.
The boundary of $\P$ is formed by $\NP$ straight edges connecting
$\NP$ vertices.
%%g
The sequence of the vertices on $\partial\P$ is oriented in the
counter-clockwise order and the vertex coordinates are denoted by
$\xv_{i}\EQUIV(\xs_{i},\ys_{i})^T$, $i = 1,2,\ldots,\NP$.
We denote the unit normal vector to $\partial\P$ pointing out of $\P$
by $\norP$.

Usually, in the convergence analysis of the conforming \VEM{}, it is
assumed that there exists a positive constant $\varrho$ independent of
$\hh$ (hence, also of $\Th$) such that for every polygonal element
$\P\in\Th$ it holds that:
\begin{enumerate}[(i)]
\item $\P$ is star-shaped with respect to a disk with
  radius greater than $\varrho\hP$;
\item  for every edge $\E\in\partial\P$ it holds that
  $\hE\geq\varrho\hP$.
\end{enumerate}

Although the convergence analysis of the \XVEM{} is beyond the scope
of this paper, we present such mesh regularity assumptions to
characterize the geometry of the elements in the polygonal meshes,
which is pertinent to our formulation.
We also note that condition (i) implies that all the mesh elements
have a finite number of vertices and edges for $\hh\to0$ and are
simply connected subset of $\REAL^2$.
In turn, condition (ii) excludes the possibility of collapsing
vertices in the refinement process, i.e., vertices whose distance
becomes zero faster than $\hh$.

\subsection{Conforming virtual element space, elliptic projection and bilinear form} %%%  $\ash(\cdot,\cdot)$}
\label{subsec:conforming:VEM:regular}
Let $\Gamma_c=\emptyset$.
On every polygonal element $\P$ with boundary $\partial\P$, we first
define the following scalar virtual element space
\begin{align}
  %%\label{eq:basis2}
  \Vsh(\P)
  \EQUIV \Big\{ \vsh\in\HS{1}(\P) : \Delta\vsh=0,\,\,
  \restrict{\vsh}{\partial\P}\in\CS{0}(\partial\P),\,\,
  \restrict{\vsh}{\E}\in\PS{1}(\E)\,\,\forall\E\in\partial\P)
  \Big\},
\end{align}
where $\PS{1}(\E)$ is the set of linear polynomials on the element
edge $\E\in\partial\P$ and $\Delta$ is the Laplace operator.
We denote the canonical basis of $\Vh(\P)$ by
$\{\varphi_i\}_{i=1}^{N_E}$, so that each $\varphi_i$ is the harmonic
function on $\P$ with continuous piecewise linear trace on the
boundary $\partial\P$ that takes value $1$ on the $i$-th node and $0$
on the remaining nodes.
The linear polynomials $\PS{1}(\P)$ are a subspace of $\Vsh(\P)$, and
the basis functions $\varphi_i$ satisfies the partition-of-unity
property
\begin{align}
  \sum_{i=1}^{\NP}\varphi_i(\xv) = 1
  \quad\forall\xv\in\P.
  \label{eq:PUM:scalar:basis:functions}
\end{align}

For the linear elasticity (vectorial) problem, on
every polygonal element $\P\in\Th$ we define the local virtual element
space of vector-valued functions as $\Vvh(\P) = \big[\Vsh(\P)\big]^2$.
Every vector-valued virtual element function $\vvh\in
\Vvh(\P)$ is uniquely characterized by its vertex values, also
known as the \emph{degrees of freedom} (DOFs) of the method.
In the framework of two-dimensional elasticity, such degrees of
freedom represent the two components of the displacement field at the
mesh vertices.
Therefore, we have $2\NP$ degrees of freedom per mesh element
$\P$. Such degrees of freedom are unisolvent in $\Vvh(\P)$~\cite{Beirao-Brezzi-Marini:2013}.

We define the set of `canonical' basis functions of $\Vvh(\P)$ by
$\big\{\phiv_i\big\}_{i=1}^{2N_{\P}}$ so that
$\phiv_{2i-1}=(\varphi_i,0)^T$ and $\phi_{2i}=(0,\varphi_i)^T$ for
$i=1,\ldots,\NP$.
These functions are made explicit by the following expression
\begin{equation}\label{eq:stdbasis}
  \Vvh(\P) = \textrm{span}
  \left \{
  \left( \begin{array}{c} \varphi_1      \\ 0              \end{array} \right), 
  \left( \begin{array}{c} 0              \\ \varphi_1      \end{array} \right),
  \ldots,
  \left( \begin{array}{c} \varphi_{i}    \\ 0              \end{array} \right), 
  \left( \begin{array}{c} 0              \\ \varphi_{i}    \end{array} \right),
  \ldots,
  \left( \begin{array}{c} \varphi_{N_{\P}} \\ 0              \end{array} \right),
  \left( \begin{array}{c} 0              \\ \varphi_{N_{\P}} \end{array} \right) 
  \right\},
\end{equation}
and the partition-of-unity
property~\eqref{eq:PUM:scalar:basis:functions} implies that
\begin{align*}
  \sum_{i=1}^{\NP}\phiv_{2i-1}(\xv)
  = \left(\begin{array}{c} \sum_{i=1}^{\NP}\varphi_i(\xv) \\[0.2em] 0 \end{array}\right)
  = \left(\begin{array}{c} 1 \\[0.2em] 0 \end{array}\right)
  \quad\text{and}\quad
  \sum_{i=1}^{\NP}\phiv_{2i}(\xv)
  = \left(\begin{array}{c} 0 \\[0.2em] \sum_{i=1}^{\NP}\varphi_i(\xv) \end{array}\right)
  = \left(\begin{array}{c} 0 \\[0.2em] 1 \end{array}\right)
  \quad\forall\xv\in\P.
\end{align*}
We collect all the element spaces $\Vsh(\P)$ in a conforming way and
define the global virtual element space $\Vvh\subset\mathscr{U}_0$ as
follows
\begin{equation*}
  %%\label{eq:globalspace}
  \Vvh =
  \Big\{
  \vvh\in\big[\HS{1}(\Omega)\big]^2:\restrict{\vvh}{\P}\in\Vvh(\P)\,\,\forall\P\in\Th
  \Big\}.
\end{equation*}
Let $\ash(\cdot,\cdot)$ and $\bsh(\cdot)$ denote computable
counterparts of the exact bilinear form $\as(\cdot,\cdot)$ and the
linear functional $\bs(\cdot)$ acting on $\Vvh$, and consider the
virtual element affine subspace of $\Vvh$ given by
\begin{align*}
  \Vvhg =
  \Big\{
  \vvh\in\Vvh:\,\vvh=\gv^h\textrm{~on~}\Gamma_u
  \Big\}, 
\end{align*}
which incorporates the essential boundary
condition~\eqref{eq:essentialBC} in the space definition by taking the linear interpolant $\gv^h$ of $\gv$, and the
linear subspace $\Vvhz\subset\Vvh$ that is obtained by setting
$\gv^h=\bm{0}$ in $\Vvhg$.
With this caveat, the virtual element approximation of the variational
problem~\eqref{eq:weak_form} reads as: Find $\uvh\in\Vvhg$ such
that
\begin{align}
  \ash( \uvh, \vvh ) = \bsh( \vvh )
  \quad 
  \forall \vvh\in\Vvhz.
  %%\subset\mathscr{U}_0.
  \label{eq:VEM:pblm:regular}
\end{align}

To construct the bilinear form $\ash(\cdot,\cdot)$ and the linear
functional $\bsh(\cdot)$, we first split them as the sum of element
terms $\ashP(\cdot,\cdot)$ and $\bshP(\cdot)$ so that
\begin{align*}
  \ash( \uvh, \vvh ) &= \sum_{\P\in\Omega} \ashP( \uvh, \vvh ) \quad \forall \uvh, \vvh\in\Vvh,\\[0.5em]
  \bsh( \uvh )       &= \sum_{\P\in\Omega} \bshP( \vvh       ) \quad \forall       \vvh\in\Vvh.
\end{align*}

It is well established in the VEM literature that a crucial
requirement for every $\ashP(\cdot,\cdot)$ to deliver an accurate and
stable formulation is to satisfy the properties of \emph{linear
consistency} and \emph{stability}
\cite{BeiraodaVeiga-Brezzi-Cangiani-Manzini-Marini-Russo:2013}.
To construct such $\ashP(\cdot,\cdot)$, we resort to the elliptic
projection operator $\Pia:\Vvh(\P)\to\big[\PS{1}(\P)\big]^2$, which
maps vector-valued functions from $\Vvh(\P)$ onto linear vector
polynomials.
To fix the nontrivial kernel in the definition of such elliptic
projector, we introduce the \emph{average translation operator} over
the $\NP$ element vertices $\big\{\xv_j\big\}_{j=1}^{\NP}$ defined
as
\begin{equation}
  \overline{\wv} = \frac{1}{\NP}\sum_{j=1}^{\NP}\wv(\xv_j),
  \label{eq:average:translation:operator} %% \label{eq:condition1b}
\end{equation}
and the \textit{average rotation operator} defined as
\begin{align}
  \overline{(\wv)}_R =
  \frac{1}{\NP}\sum_{j=1}^{\NP}\rv(\xv_j)\cdot\wv(\xv_j),
  \quad \rv(\xv) = \big(y,-x\big)^T.
  \label{eq:average:rotation:operator} %% \label{eq:condition1b}
  %% \left[
  %%   \begin{array} {c}
  %%     y \\
  %%     -x
  %%   \end{array}
  %% \right ]^T.
\end{align}
For each $\vvh\in\Vvh(\P)$, the elliptic projection $\Pia(\vvh)$
is the solution of the variational problem
\begin{subequations} 
  \begin{align}
    \int_{\P}\sigmav(\qv):\varepsilonv(\Pia\vvh) \dx &= \int_{\P}\sigmav(\qv):\varepsilonv(\vvh) \dx
    \quad\forall\qv\in\big[\PS{1}(\P)\big]^2,\label{eq:Pia:def}
    \intertext{with the additional conditions}
    \overline{ \Pia\vvh }   &= \overline{ \vvh },  \label{eq:additional:projection:A}\\[0.25em]
    \overline{(\Pia\vvh)}_R &= \overline{(\vvh)}_R.\label{eq:additional:projection:B}
  \end{align}
\end{subequations}
Conditions~\eqref{eq:additional:projection:A}
and~\eqref{eq:additional:projection:B} fix the 
rigid-body modes (two translations and one rotation)
that form the kernel of
$\varepsilonv(\cdot)$.

A requirement for such a projection operator is that it is computable
from the degrees of freedom of $\Vvh$, as we explain below.
In order to compute $\Pia(\vvh)$ it is convenient to choose, as a
basis of $\PS{1}(\P)$, the set of \emph{scaled monomials}
\begin{align}
  \mv(\xv) = \Big \{ 1,\,\xi(\xv),\,\eta(\xv) \Big \},
  \quad\textrm{with}\quad
  \xi (\xv) = \frac{\xs-\xsP}{\hP},\quad
  \eta(\xv) = \frac{\ys-\ysP}{\hP},
  \label{eq:scaled:monomials}
\end{align}
where $\xv=\big(\xs,\ys\big)^T$, so that the basis functions of $\PS{1}(\P)$
scale as $\mathcal{O}(1)$ with respect to $\hh$.
It immediately follows that $\PS{1}(\P)=\SPAN\{1,\xi,\eta\}$, and
%% we can provide the following explicit expression for the basis of
a possible basis of $\big[\PS{1}(\P)\big]^2$ is
\begin{equation}
  \big[  \PS{1}(\P)\big]^2=
  \SPAN
  \left \{
  \left( \begin{array}{c} 1    \\ 0    \end{array} \right), 
  \left( \begin{array}{c} 0    \\ 1    \end{array} \right),
  \left( \begin{array}{c} \eta \\ -\xi \end{array} \right), 
  \left( \begin{array}{c} \xi  \\ 0    \end{array} \right), 
  \left( \begin{array}{c} 0    \\ \eta \end{array} \right), 
  \left( \begin{array}{c} \eta \\ \xi  \end{array} \right)
  \right \}.
  \label{eq:P2}
\end{equation}
The six vector fields in \eqref{eq:P2} represent the three
planar rigid-body modes and the three independent nonzero deformation modes.
%% , and
%% three independent deformation modes, i.e. two axial strains and one
%% pure shear strain.

To prove the computability of $\Pia$, we rewrite \eqref{eq:Pia:def}
with~\eqref{eq:additional:projection:A}-\eqref{eq:additional:projection:B}
as a linear system.
For every $\phiv_i$ from the canonical basis of $\Vvh(\P)$ shown in
\eqref{eq:stdbasis}, we consider the expansion of $\Pia\phiv_i$ on the
basis of $[\PS{1}(\P)]^2$ shown in \eqref{eq:P2}.
A suitable application of the divergence theorem shows that $\Pia\phiv_i$ is 
computable by using only the degrees of freedom of $\phiv_i$ and noting that
$\nabla\cdot\sigmav(\phiv_i)=\bm{0}$.
The polynomial projection $\Pia\vvh$ can readily be computed for all
virtual element fields $\vvh$ from the projections of the basis
functions $\phiv_i$ because the projection operator is a linear
operator.

We will expand on this observation in the next section.

Once computed, operator $\Pia$ allows us to evaluate the local
approximated bilinear form as follows
\begin{align*}
  %\label{approx-bil}
  \ashP( \vvh, \wvh )
  &= \asP\big(\Pia(\vvh),\,\Pia(\wvh)\big)
  + \SP\big(\big(\vvh-\Pia(\vvh)\big),\,\big(\wvh-\Pia(\wvh)\big)\big)
  \\[1em]
  &= \int_{\P} \sigmav(\Pia(\vvh)) : \varepsilonv(\Pia(\wvh) \dx
  + \SP\big(\big(\vvh-\Pia(\vvh)\big),\,\big(\wvh-\Pia(\wvh)\big)\big),
\end{align*}
where $\SP(\cdot,\cdot)$ is a suitable stabilizing term that preserves the coercivity of the system.
According to the virtual element methodology, $\SP(\cdot,\cdot)$ can
be any symmetric, positive definite, continuous
bilinear form defined on the kernel of the projection operator
$\Pia$~\cite{BeiraodaVeiga-Brezzi-Cangiani-Manzini-Marini-Russo:2013}.

We refer the reader to Section~\ref{sec:implementation} for possible
choices of the stabilization term. 

\medskip
Finally, the expression for the virtual
element approximation of the linear functional in the right-hand side
of~\eqref{eq:VEM:pblm:regular} is given by
\begin{align*}
  \bshP(\vvh)
  = \int_{\Gamma^t\cap\partial\P}\bar{\bm{t}}\cdot\vvh\,d\Gamma
  = \bsP(\vvh),
\end{align*}
where $\bshP(\vvh)$ is computable because $\bar{\tv}$ is
known and the trace of $\vvh$ is a linear polynomial on each edge
$\E\in\Gamma^t\cap\partial\P$ that is known through the interpolation
of the edge degrees of freedom.

\subsection{Extended virtual element space, elliptic projection and bilinear form} %%% $\ashX(\cdot,\cdot)$}
\begin{figure}[!htb]
  \centering
  \includegraphics[width=0.5\textwidth]{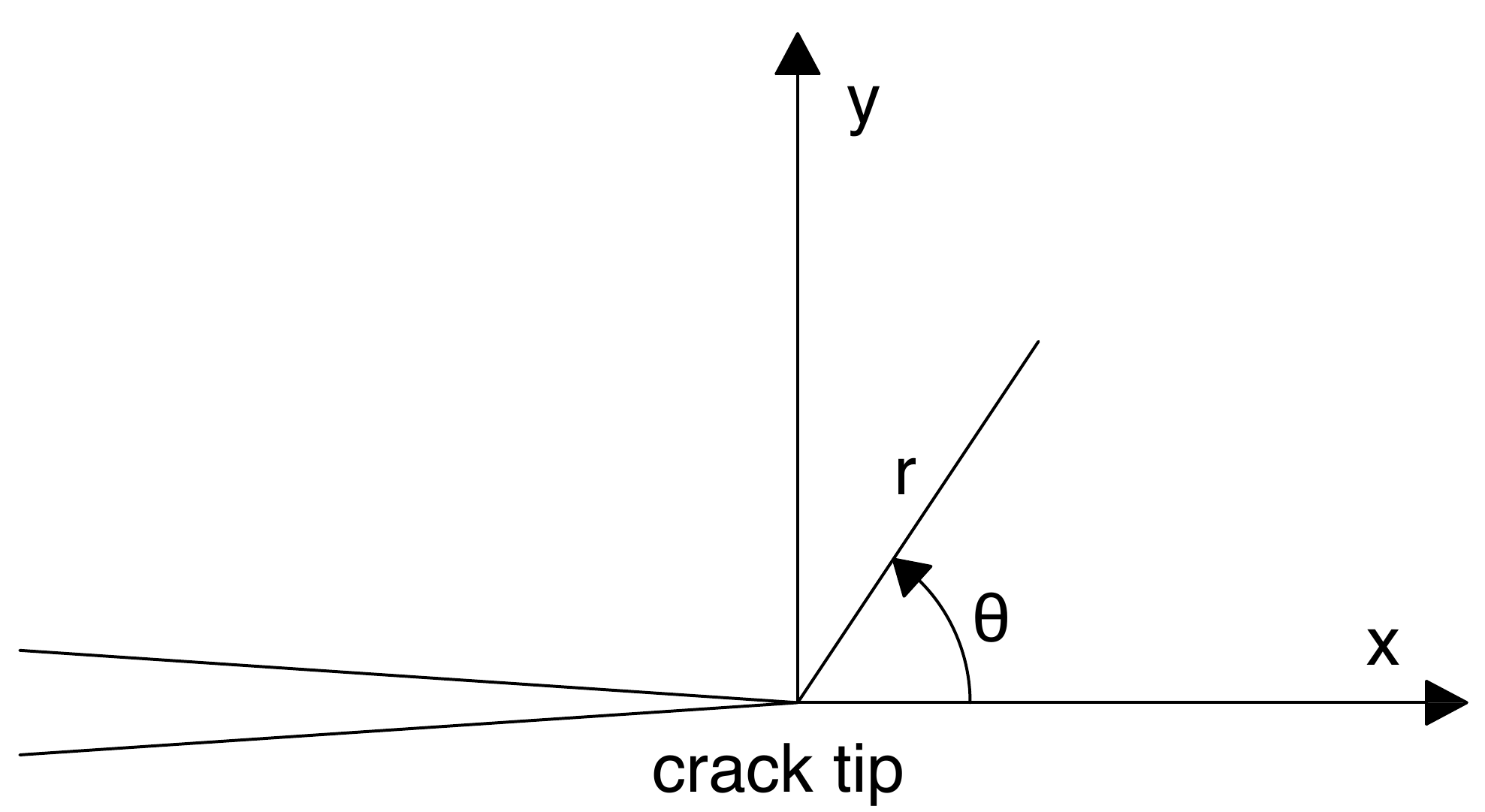}
  \caption{Local crack-tip reference system in polar coordinates.}\label{fig:local_ref}
\end{figure}

\label{subsec:conforming:VEM:extended}
If the exact solution to the selected problem contains singularities,
then similar to the finite element method, the accuracy of the virtual element method is compromised.
For this reason, it is beneficial to enrich the virtual element space
by means of independent fields carrying information about the
singularities affecting the exact solution.
As we discuss later on, such fields are required to satisfy the
equilibrium equations~\eqref{eq:equilibrium:a}.
For two-dimensional elastic fracture problems, we choose the
enrichment fields as a scaled form of the exact asymptotic crack-tip displacement fields for mode I and mode II crack opening, $\uvI=\big(\usxI,\usyI\big)^T$ and
$\uvII=\big(\usxII,\usyII)^T$, respectively.
These enrichment fields are given by the expressions:
\begin{subequations} \label{def:crackmodes}
  \begin{align}
    \usxI  &:= \usxI(r,\theta) = \sqrt{\frac{r}{2\pi}} \left[ (2\kappa - 1)\cos \left( \frac{\theta}{2} \right) - \cos \left( \frac{3\theta}{2} \right) \right], \\[0.75em]
    \usyI  &:= \usyI(r,\theta) = \sqrt{\frac{r}{2\pi}} \left[ (2\kappa + 1)\sin \left( \frac{\theta}{2} \right) - \sin \left( \frac{3\theta}{2} \right) \right], \\[0.75em]
    \usxII &:= \usxII(r,\theta) =\sqrt{\frac{r}{2\pi}} \left[ (2\kappa + 3)\sin \left( \frac{\theta}{2} \right) + \sin \left( \frac{3\theta}{2} \right) \right], \\[0.75em]
    \usyII &:= \usyII(r,\theta) = -\sqrt{\frac{r}{2\pi}} \left[ (2\kappa - 3)\cos \left( \frac{\theta}{2} \right) + \cos \left( \frac{3\theta}{2} \right) \right], 
  \end{align}
\end{subequations}
where $(\rs,\theta)$ are polar coordinates in the local crack tip reference system (see Fig.~\ref{fig:local_ref}) and $\kappa$ is the Kolosov constant. % E_0$ is the Young modulus and $\nu$ is the Poisson ratio of the material.

An explicit computation implies that these fields satisfy equilibrium, i.e., the conditions
$\nabla\cdot\sigmav(\uvI)=0$ and $\nabla\cdot\sigmav(\uvII)=0$ hold.
Note that $\uvI$ and $\uvII$ belong to $\HS{\frac{3}{2}-\eta}(\Omega)$
for any $\eta>0$~\cite{Grisvard:1992}, and this fact reduces the
convergence rate of a standard finite element or virtual element
method to $\mathcal{O}(\hh^{\frac12})$.

Let
$\uvIt = \big(\usxIt, \usyIt \big)^T$ and $\uvIIt = \big(\usxIIt,
\usyIIt\big)^T$ denote
the dimensionless version of fields $\uvI$ and $\uvII$, respectively,
\begin{align}
  %%\uvIt  = \uvI /\hP^{1/2}
  \uvIt  = \uvI /\hh^{1/2}
  \quad\textrm{and}\quad
  %%\uvIIt = \uvII/\hP^{1/2}.
  \uvIIt = \uvII/\hh^{1/2},
  \label{def:scaledcrackmodes}
\end{align}
where $\hh=\max_{\P\in\Th}(\hP)$.
%%
%%
%% \begin{subequations} \label{def:scaledcrackmodes}
%%   \begin{align}
%%     \uvIt  &= \uvI /\hP^{1/2}, \label{def:scaledcrackmodes1} \\[0.25em]
%%     \uvIIt &= \uvII/\hP^{1/2}. \label{def:scaledcrackmodes2}
%%   \end{align}
%% \end{subequations}
%%
In order to define the extended virtual element space, we first
introduce the local virtual element space $\Vvha(\P)\subset\Vvh(\P)$,
which reads as
\begin{align}\label{def:vhstar}
  \Vvha(\P)
  \EQUIV
  \left\{
  \vvh=(\vshx,\vshy)^T\in\Vvh(\P):\vshx=\vshy
  \right\}.
\end{align}
This space is generated by the linear combination of the basis
functions $\phiva_i=(\varphi_i,\varphi_i)^T$, $i=1,\ldots,\NP$, where
the functions $\varphi_i$ are the basis functions of the scalar
virtual element space $\Vsh(\P)$, so that
$\Vvha(\P)=\SPAN\big\{\phiva_1,\ldots,\phiva_{\NP}\big\}$.
The dimension of this space is clearly $\NP$ and the
partition-of-unity property of functions $\varphi_i$ implies that
\begin{align*}
  \sum_{i=1}^{\NP}\phiva_i(\xv)
  = \left(\begin{array}{c} \sum_{i=1}^{\NP}\phi_i(\xv) \\[0.2em] \sum_{i=1}^{\NP}\phi_i(\xv) \end{array}\right)
  = \left(\begin{array}{c} 1 \\[0.2em] 1 \end{array}\right)
  \qquad\forall\xv\in\P.
\end{align*}
Then, we define the matrices $\psivI$ and $\psivII$ as
\begin{align}\label{def:matricespsi}
  \psivI
  \EQUIV
  \left [ \begin{array}{cc} \usxIt  & 0 \\ 0 & \usyIt  \end{array} \right], \quad
  \psivII
  \EQUIV
  \left [ \begin{array}{cc} \usxIIt & 0 \\ 0 & \usyIIt \end{array} \right].
\end{align}

\medskip
We now have all the ingredients to define the \textit{local extended
  virtual element space} $\VvhX(\P)$, which reads as
\begin{align}\label{def:xvemspace}
  \VvhX(\P) \EQUIV \Vvh(\P) \oplus \psivI\Vvha(\P) \oplus \psivII \Vvha(\P).
\end{align}
We obtain a basis of this space as the union of the basis functions of
$\VvhX(\P)$, $\psivI\Vvha(\P)$ and $\psivII \Vvha(\P)$, so that
\begin{align}
  \VvhX(\P)
  &=
  \SPAN\Big\{\, \phiv_1, \phiv_{2}, \ldots, \phiv_{2i-1}, \phiv_{2i}, \ldots, \phiv_{2\NP-1}, \phiv_{2\NP} \,\Big\}
  \cup \psivI \SPAN\Big\{\, \phiva_1, \phiva_{2}, \ldots, \phiva_{\NP} \,\Big\}
  \nonumber\\[0.5em]
  &\qquad
  \cup
  \psivII\SPAN\Big\{\, \phiva_1, \phiva_{2}, \ldots, \phiva_{\NP} \,\Big\},
  \label{eq:basis:xvem}
\end{align}
where we recall that $\phiv_{2i-1}=(\varphi_i,0)^T$,
$\phiv_{2i}=(0,\varphi_i)^T$ and $\phiva_i=(\varphi_i,\varphi_i)^T$,
$i=1,\ldots,\NP$.
Therefore, at every enriched node the vector-valued field $\vvhX(\xv)$
that belongs to the extended virtual element space $\VvhX(\P)$ is
characterized by four values and for an element whose nodes are all
enriched, we have $4\NP$ degrees of freedom.
For example, at the $j$-th node with coordinates $\xv_j$, we find that
\begin{align*}
  \vvh(\xv_j) &= \sum_{i=1}^{\NP}\Bigg[
    \vsh_{i,x}   \left(\begin{array}{c} \varphi_i(\xv_j) \\[0.25em] 0 \end{array}\right)
    + \vsh_{i,y} \left(\begin{array}{c} 0 \\[0.25em] \varphi_i(\xv_j) \end{array}\right)
    + \vsh_{i,I} \left(\begin{array}{c} \usxIt (\xv_j) \varphi_i(\xv_j) \\[0.25em] \usyIt (\xv_j)\varphi_i(\xv_j) \end{array}\right)
    + \vsh_{i,II}\left(\begin{array}{c} \usxIIt(\xv_j) \varphi_i(\xv_j) \\[0.25em] \usyIIt(\xv_j)\varphi_i(\xv_j) \end{array}\right)
    \Bigg]
  \\[0.5em]
  &= \left(\begin{array}{c}
    \vsh_{j,x} + \vsh_{j,I} \usxIt(\xv_j) + \vsh_{j,I} \usxIIt(\xv_j) \\[0.25em]
    \vsh_{j,y} + \vsh_{j,II}\usyIt(\xv_j) + \vsh_{j,II}\usyIIt(\xv_j)
  \end{array}\right),
\end{align*}
since $\varphi_i(\xv_j)= \delta_{ij}$.
\begin{remark}
  Here, $\vsh_{i,x}$, $\vsh_{i,y}$, $\vsh_{i,I}$ and $\vsh_{i,II}$ are
  the coefficients of the basis functions in~\eqref{eq:basis:xvem} and
  can thus be identified with the degrees of freedom of the method.
  Note, however, that the degrees of freedom of an enriched function
  $\vvhX\in\VvhX(\P)$ are no longer the values of $\vvhX$ at the
  vertices of element $E$.
\end{remark}
To ease the exposition, we denote the basis functions of $\VvhX(\P)$
by the symbol $\phiv_i$, $i=1,2,\ldots,4\NP$, so that
$\VvhX(\P)=\SPAN\big\{\,\phiv_{1},\,\phiv_{2},\ldots,\,\phiv_{4\NP}\,\big\}$
where
\begin{align*}
  \phiv_i =
  \begin{cases}
    \Big( \varphi_i,        0             \Big)^T & \mbox{for~$1     \leq\is\leq2\NP,\,i$ odd},  \\[0.5em]
    \Big( 0,             \varphi_i        \Big)^T & \mbox{for~$1     \leq\is\leq2\NP,\,i$ even}, \\[0.5em]
    \Big( \usxIt\varphi_i,  \usyIt\varphi_i  \Big)^T & \mbox{for~$1+2\NP\leq\is\leq3\NP$},          \\[0.5em]
    \Big( \usxIIt\varphi_i, \usyIIt\varphi_i \Big)^T & \mbox{for~$1+3\NP\leq\is\leq4\NP$}.
  \end{cases}
\end{align*}
Finally, the extended global virtual element space $\VvhX$ is defined
as follows:
\begin{equation*}
  %%\label{eq:globalspace}
  \VvhX =
  \Big\{
  \vvhX\in\big[\HS{1}(\Omega)\big]^2:\restrict{\vvhX}{\P}\in\VvhX(\P)\,\quad\forall\P\in\Th
  \Big\}.
\end{equation*}
Again, to consider the essential boundary
condition~\eqref{eq:essentialBC} we consider the affine subspace
$\VvhXg$ of $\VvhX$ defined by
\begin{align*}
  \VvhXg = \Big\{
  \vvhX\in\VvhX:\,\vvhX=\gv^h_X~\textrm{on}~\Gamma_u
  \Big\},
\end{align*}
where $\gv^h_X$ is the extended linear interpolant of $\gv$, and the linear subspace $\VvhXz$, which is defined by setting
$\gv^h_X=\bm{0}$ in the above definition.

Since $\{\phiv_i\}_{i=1}^{4\NP}$ are not known in the interior of the
element, we construct a convenient projection operator that will allow
us to obtain computable approximations
$\ashX(\cdot,\cdot):\VvhX(\P)\times\VvhX(\P)\to\REAL$ and
$\bshX(\cdot):\VvhX(\P)\to\REAL$ of the exact bilinear form
$\as(\cdot,\cdot)$ and the linear functional $\bs(\cdot)$ appearing in
\eqref{eq:weak_form}.
The extended virtual element formulation then reads: Find
$\uvhX\in\VvhXg$ such that
\begin{align}
  \ashX(\uvhX,\vvhX) = \bshX(\vvhX)
  \quad\forall\vvhX\in\VvhXz,
  %%\label{eq:pblm:xvem}
\end{align}
where the bilinear form $\ashX(\cdot,\cdot)$ is built element-wise as
\begin{align}
  \ashX( \uvhX, \vvhX )
  = \sum_{\P\in\Omega} \ashPX( \uvhX, \vvhX )
  \quad \forall\uvhX,\vvhX\in\VvhX,
\end{align}
and again we set $\bshX(\vvhX)=\bs(\vvhX)$.

In order to construct a consistent and stable bilinear form
$\ashPX(\cdot,\cdot)$, we extend the polynomial space $\PS{1}(\P)$ to
a subspace of $\Vsh(\P)$ including the linear polynomials and the
additional enrichment functions $\uvIt$ and $\uvIIt$, so that
\begin{align*}
  \PSX{}(\P) \EQUIV \PS{1}(\P) \oplus \SPAN(\uvIt,\uvIIt).
\end{align*}
Space $\PSX{}(\P)$ is spanned by the eight linearly independent vector fields
\begin{equation}
  \label{eq:PX}
  \PSX{}(\P) = \SPAN
  \left \{\,
  \left( \begin{array}{c} 1        \\ 0       \end{array} \right), 
  \left( \begin{array}{c} 0        \\ 1       \end{array} \right),
  \left( \begin{array}{c} \eta     \\ -\xi    \end{array} \right), 
  \left( \begin{array}{c} \xi      \\ 0       \end{array} \right), 
  \left( \begin{array}{c} 0        \\ \eta    \end{array} \right), 
  \left( \begin{array}{c} \eta     \\ \xi     \end{array} \right), 
  \left( \begin{array}{c} \usxIt   \\ \usyIt  \end{array} \right),
  \left( \begin{array}{c} \usxIIt  \\ \usyIIt \end{array} \right)
  \,\right \}.
\end{equation}
The first six vector fields in \eqref{eq:PX} represent the three
fundamental rigid body motions and the three independent deformation
modes that form $\PS{1}(\P)$, cf.~\eqref{eq:P2}.
The last two vector fields are the scaled enrichment fields chosen to
construct the extended virtual element space $\VvhX(\P)$.
\begin{remark}
  \label{remark:equilibrium:extended:PX}
  All $\qvX\in\PSX{}(\P)$ satisfy the equilibrium equation
  $\nabla\cdot\sigmav(\qvX)=\bm{0}$.
  This property is crucial to determine the computability of the
  extended projection operator $\PiaX$.
\end{remark}

To construct a bilinear form $\ashPX(\cdot,\cdot)$ for which such
properties hold, we define the extended elliptic projection operator
$\PiaX:\VvhX(\P)\to\PSX{}(\P)$ for each element $\P$.
For each $\vvhX\in\VvhX(\P)$, the extended elliptic projection
$\PiaX(\vvhX)$ is the solution of the variational
problem
\begin{subequations} 
  \begin{align}
    \int_{\P} \sigmav(\qvX):\varepsilonv(\Pi_X^a\vvhX) \dx &=
    \int_{\P} \sigmav(\qvX):\varepsilonv(\vvhX) \dx
    \quad\forall\qvX\in\PSX{}(\P),
    \label{eq:cond_extra_xvem:A}
    \intertext{with the additional conditions}
    \overline{\PiaX  \vvhX}    &= \overline{ \vvhX },  \label{eq:cond_extra_xvem:B}\\[0.25em]
    \overline{(\PiaX \vvhX)}_R &= \overline{(\vvhX)}_R,\label{eq:cond_extra_xvem:C}
  \end{align}
\end{subequations}
where $\overline{(\cdot)}$ and $\overline{(\cdot)}_R$ are the average
translation and rotation, respectively, which are defined 
in~\eqref{eq:average:translation:operator} 
and~\eqref{eq:average:rotation:operator}.
Recalling the divergence theorem and
Remark~\ref{remark:equilibrium:extended:PX}, the vector polynomial
$\PiaX\vvhX\in\PSX{}(\P)$ is computable from the degrees of freedom of
$\vvhX$.

The projection operator $\PiaX$ allows us to define the local extended
bilinear form as follows:
\begin{align}
  \ashPX( \vvhX, \wvhX )
  &\EQUIV
  \asP\Big(\PiaX(\vvhX),\,\PiaX(\wvhX)\Big)
  + \SPX\Big(\vvhX-\PiaX(\vvhX),\,\wvhX-\PiaX(\wvhX)\Big)
  \nonumber\\[1em]
  &= \int_{\P} \sigmav\big(\PiaX(\vvhX)\big) : \varepsilonv\big(\PiaX(\wvhX)\big) \dx
  + \SPX\Big(\vvhX-\PiaX(\vvhX),\,\wvhX-\PiaX(\wvhX)\Big),
  \label{eq:ashPX:def}
\end{align}
where $\SPX(\cdot,\cdot)$ is a stabilization term that must be
suitably defined to guarantee linear consistency
(cf.~\eqref{eq:linear-consistency}) and stability
(cf.~\eqref{eq:stability}) of the method.
Again, according to the virtual element methodology,
$\SPX(\cdot,\cdot)$ can be any symmetric, positive definite,
continuous bilinear form defined on the kernel of
the extended projection operator
$\PiaX$~\cite{Beirao-Brezzi-Marini:2013}.
The reader is referred to Section~\ref{sec:implementation} for
possible choices of the stabilization term.

With a suitable choice of the stabilization term, the bilinear form
$\ashPX(\cdot,\cdot)$ has the following properties, which are
fundamental in order to guarantee the convergence of the method:
\begin{description}
\item[$(i)$] {\emph{extended linear consistency}}: for all
  $\vvhX\in\VvhX(\P)$ and $\qvX\in\PSX{}(\P)$ it holds that
  \begin{align}
    \label{eq:linear-consistency}
    \ashPX(\vvhX,\qvX) = \asP(\vvhX,\qvX);
  \end{align}
\item[$(ii)$] {\emph{stability}}: there exist two positive constants
  $\alpha_*,\,\alpha^*$, independent of $\hh$ and $\P$, such that
  \begin{align}
    \label{eq:stability}
    \alpha_*\asP(\vvhX,\vvhX)
    \leq\ashPX(\vvhX,\vvhX)
    \leq\alpha^*\asP(\vvhX,\vvhX)\quad\forall\vvhX\in\VvhX(\P).
  \end{align}
\end{description}
According to the virtual element theory,
cf.~\cite{BeiraodaVeiga-Brezzi-Cangiani-Manzini-Marini-Russo:2013},
the constants $\alpha_*$ and $\alpha^*$ must be independent of the mesh
size parameter $\hh$.
However, they can depend on the other model and discretization
parameters such as the bound on $\Cv$ and the mesh regularity constant
$\rho$.
Here, $\asP(\cdot,\cdot)$ is the local coercive and continuous
bilinear form
\begin{align*}
  \asP(\uv,\vv)
  = \int_{\P}\sigmav(\vv) : \varepsilonv(\uv)\dx
  \quad\forall\uv,\vv\in\mathscr{U}_0.
\end{align*}

\begin{remark}
  In Section~\ref{sec:implementation}, we provide two possible choices
  of the stabilization term by considering the standard
  \emph{dofi-dofi} and \emph{D-recipe} formulations in our extended
  setting.
  Such choices are widely accepted in the VEM literature and in some
  cases they were theoretically proved to be effective to guarantee
  stability relations such as \eqref{eq:stability}.
  However, the choice of the stabilization term in the presence
  of enrichment functions and its impact on the
  behavior of the VEM are still open issues 
  at this time. 
  For example, it would be desirable that the constants of the
  stability relation~\eqref{eq:stability} are independent of the 
  Young's modulus and Poisson's ratio 
  to realize a robust discretization.
  These topics will be the subject of future work.
\end{remark}

\subsection{Partial enrichment}
\label{subsec:partial-enrichment}

Let $\P$ denote an element of mesh $\Th$ and $\kP$ a positive integer
number strictly less than $\NP$ (the case for $\kP=\NP$ is the full
enrichment case).
We select $\kP$ distinct nodes of element $\P$ to be enriched and the
corresponding basis functions $\phiva_{\is_{\ell}}\in\Vvha(\P)$
labeled by the $\kP$ distinct indices $\is_{\ell}\in[1,\NP]$ for
$\ell=1,\ldots,\kP$.
We formally denote the subset of these indices by
$\calI=\{\is_{\ell}\}_{\ell=1}^{\kP}$.
Using these basis functions, we define the \emph{reduced virtual
element space}
\begin{align*}
  \Vvhat(\P) \EQUIV
  \SPAN
  \Big\{\,
  \phiva_{\is_1}, \phiva_{\is_2}, \ldots, \phiva_{\is_{\kP}} 
  \,\Big\}
  \subset
  \Vvha(\P)
\end{align*}
and the \emph{reduced extended virtual element space}
\begin{align*}
  \VvhXt(\P) =
  \Vvh(\P) \oplus \psivI\Vvhat(\P) \oplus \psivII\Vvhat(\P)
  \subset \VvhX,
\end{align*}
%%
%% single quotes for a word !
where a tilde accent as a superscript is used to denote all `reduced' mathematical objects.
Equivalently, we can define the reduced virtual element space
$\VvhXt(\P)$ as the span of the basis functions of $\Vvh(\P)$,
$\psivI\Vvhat(\P)$ and $\psivII\Vvhat(\P)$, so that
\begin{align}
  \VvhXt(\P)
  &=
  \SPAN\Big\{\, \phiv_1, \phiv_{2}, \ldots, \phiv_{2i-1}, \phiv_{2i}, \ldots, \phiv_{2\NP-1}, \phiv_{2\NP} \,\Big\}
  \cup
  \psivI \SPAN\Big\{\, \phiva_1, \phiva_{2}, \ldots, \phiva_{\kP} \,\Big\}
  \nonumber\\[0.5em]
  &\qquad
  \cup
  \psivII\SPAN\Big\{\, \phiva_1, \phiva_{2}, \ldots, \phiva_{\kP} \,\Big\},
  \label{eq:reduced:basis:xvem}
\end{align}
which can be compared to~\eqref{eq:basis:xvem}.
Accordingly, a generic virtual element function that belongs to the
reduced space $\VvhXt(\P)$ is described by $2\NP+2\kP$ degrees of
freedom instead of $4\NP$ degrees of freedom.
The first $2\NP$ degrees of freedom are the vertex values of a
vector-valued field $\vvh\in\Vvh(\P)$.
The other $2\kP$ degrees of freedom correspond to the vertex values
of a virtual vector-valued function that belongs to the enriching
space $\psivI\Vvhat(\P)\oplus\psivII\Vvhat(\P)$ and clearly depends on
$\uvIt$ and $\uvIIt$.
We outline a few important facts that will be crucial in the
implementation of the partially enriched virtual element method.
First, the set of basis functions $\phiva_{\is_{\ell}}$ for
$\ell=1,\ldots,\kP$ does not satisfy a partition-of-unity property.
Consequently, the enriching fields $\uvIt$ and $\uvIIt$ are not
elements of $\psivI\Vvhat(\P)\oplus\psivII\Vvhat(\P)$ and the extended
space $\PSX{}(\P)$ cannot be a subspace of $\Vvhat(\P)$.
However, since $\VvhXt(\P)$ is a linear subspace of $\VvhX(\P)$, we
can still apply the projection operator $\PiaX$ to its functions and
obtain a projection in the extended space $\PSX{}(\P)$, and the
construction of the bilinear form $\ashPX(\cdot,\cdot)$ of the
previous section still holds.
For a proper formal definition, we introduce the \emph{extension (or
injection) operator} $\extkP:\Vvhat(\P)\to\Vvha(\P)$ that remaps any
reduced virtual element function $\vvhXt\in\VvhXt(\P)$ into the fully
enriched function $\extkP(\vvhXt)\in\VvhX(\P)$ such that:
\begin{align}
  \mbox{$i$-th DOF of $\extkP(\vvhXt)$} =
  \begin{cases}
    \mbox{$i$-th DOF of $\vvhXt$}
    & \text{if } 1\leq\is\leq2\NP,
    \\[0.5em]
    \mbox{$\is_{\ell}$-th DOF of $\vvhXt$}
    & \text{if } \is=2\NP+\is_{\ell} \text{~or~} \is=3\NP+\is_{\ell}
    \text{~with~} \is_{\ell}\in\calI,
    \\[0.5em]
    0 & \text{otherwise}.
  \end{cases}
  \label{eq:extkP:def}
\end{align}
Practically speaking, the remapped function has the same degrees of
freedom of the reduced functions and zero at all the additional
degrees of freedom that correspond to the nonenriched nodes.
Then, we define a new stiffness bilinear form
$\ashPXt(\cdot,\cdot):\VvhXt(\P)\times\VvhXt(\P)\to\REAL$ as
\begin{align}
  \ashPXt\Big(\vvhXt,\wvhXt\Big) := \ashPX\Big(\extkP\big(\vvhXt\big),\extkP\big(\wvhXt\big)\Big),
  \label{eq:ashPXt:def}
\end{align}
so that we can reuse the definition of $\ashPX(\cdot,\cdot)$.
Furthermore, the whole construction of the previous section, including
the consistency and stability properties, still holds.

\medskip
As we discuss in the implementation section, this formal approach
also suggests a straightforward way (but perhaps not
the most efficient one) to implement the partial enrichment as all we
need in practice is to apply a matrix representation of the injection
operator $\extkP$ to the element stiffness matrix of a fully
enriched element.
We will see that this procedure is equivalent to first constructing the fully
enriched stiffness matrix, and then simply 
suppressing all rows and
columns that correspond to the degrees of freedom of the
nonenriched nodes.

As we note in Section~\ref{sec5:numerical_examples}, partial enrichment induces a loss of optimal convergence, which also occurs in the 
X-FEM. This consequence is not surprising, since
even though we are projecting onto a space consisting of 
polynomials and nonpolynomial near-tip enrichment fields, 
% in the case of partial enrichment,
in this case the local extended virtual element space is not sufficiently rich to approximate the singular behaviour of the function near the crack tip. Special enrichment strategies can be devised to overcome this issue, for instance using the so-called \textit{geometric enrichment}.

\subsection{Embedding discontinuities}
\label{subsec:disc}

In this section, we show how both the regular and the extended virtual
element formulations presented in Sections
\ref{subsec:conforming:VEM:regular} and
\ref{subsec:conforming:VEM:extended} can be endowed with a structure
that allows discontinuous fields to be embedded within the virtual
element space.
Consider a crack $\gamma$ that intersects some of the elements in a
mesh, and define $d(\xv)$ as the signed distance from a point $\xv$ to
$\gamma$.
For modeling strong discontinuities like a crack, it would be
convenient to consider enrichment with the generalized Heaviside
function $H(\xv)$, which is equal to $+1$ for points with $d(\xv)\ge0$
($\xv$ is on or above the crack) and $-1$ for points with $d(\xv)<0$
($\xv$ is below the crack).
As in the X-FEM, we could enrich those nodes whose basis function's
support intersects the interior of the crack (not including the tips)
with $H(\xv)$.
However, the resulting extended projection $\Pi_X^a$ onto
$\PSX{}(\P)$ would not be directly computable from the degrees of
freedom of the method because the corresponding enriched virtual
element basis functions $H\phiv_i$ are not known along the crack.

To deliver a viable solution, we let the element $\P$ to be
partitioned by the discontinuity $\gamma$ into two subdomains $\P^{-}$
and $\P^{+}$.
Following \cite{Benvenuti:2019:CMAME}, in order to represent two
independent linear polynomials on $\P^{-}$ and $\P^{+}$, we adopt the
approach of Hansbo and Hansbo~\cite{Hansbo:2004:FEM}
and tailor it to
the \XVEM{}.
It is known that the approach of Hansbo and Hansbo is equivalent to
the standard X-FEM approximation with Heaviside
enrichment~\cite{Areias:2006:ACO}.
To this end, let $\Ndofs^{\textrm{VE}}$ denote the number of degrees
of freedom for element $\P$, such that $\Ndofs^{\textrm{VE}}=2\NP$ for
the virtual element formulation in \ref{subsec:conforming:VEM:regular}
and $\Ndofs^{\textrm{VE}}=4\NP$ for the extended virtual element
formulation in Section~\ref{subsec:conforming:VEM:extended}.
Each one of the $\Ndofs^{\textrm{VE}}$ virtual shape functions,
$\phiv_i$ on $\P$, is written as the sum of two new virtual shape
functions $\phiv^{-}_i$ and $\phiv^{+}_i$ that are both discontinuous
across the crack, and are defined as follows:
\begin{align}
  \label{eq:phiHH:def}
  \phiv_i^{+} =
  \begin{cases}
    \displaystyle 0       & \textrm{in}~\P^{-}\\[0.25em]
    \displaystyle \phiv_i & \textrm{in}~\P^{+}
  \end{cases},\quad
  \phiv_i^{-} =
  \begin{cases}
    \displaystyle \phiv_i & \textrm{in}~\P^{-}\\[0.25em]
    \displaystyle 0       & \textrm{in}~\P^{+}
  \end{cases}.
\end{align}

Clearly, $\phiv_i^{-}$ and $\phiv_i^{+}$ are harmonic and continuous
functions in $\P^{-}$ and $\P^{+}$, respectively, and
$\phiv_i=\phiv_i^{-} + \phiv_i^{+}$.
Proceeding likewise for all the degrees of freedom in the element, we
can generate $\Ndofs^{\textrm{HH}} = 2\Ndofs^{\textrm{VE}}$
discontinuous functions, starting from the initial
$\Ndofs^{\textrm{VE}}$ virtual basis functions.
This choice implies doubling the nodal DOFs of the element.
Therefore, the number of degrees of freedom for the element with an
internal discontinuity is twice that of the original element, and a
virtual element basis is constructed by considering two copies of the
original virtual element basis functions, restricted to $\P^{-}$ and
$\P^{+}$ respectively, as defined in~\eqref{eq:phiHH:def}.

We now define the local virtual element space to which the
discontinuous approximate solution belongs.
For the sake of simplicity, we present the derivation with respect to
the formulation presented in
Section~\ref{subsec:conforming:VEM:regular}.
Consider the following spaces:
\begin{align*}
  \Vvhm(\P)
  &\EQUIV \Big\{\, \vvh\in\big[\HONE(\P{-})\big]^2\,:\,
  \Delta\restrict{\vvh}{\P^{-}} = \bm{0}, \ 
  \restrict{\vvh}{\partial\P^{-}}\in[\CS{0}(\partial\P^{-})]^2, \\
  &\quad \ \ \ \restrict{\vvh}{\E}\in\big[\PS{1}(\E)\big]^2\,\,\forall\E\in(\partial\P\cap\partial\P^{-}),\
  \restrict{\vvh}{\P^{+}} = \bm{0}
  \Big\},\\[1em]
  \Vvhp(\P)
  &\EQUIV \Big\{\, 
  \vvh\in\big[\HONE(\P^{+})\big]^2\,:\,
  \Delta\restrict{\vvh}{\P^{+}} = \bm{0}, \ 
  \restrict{\vvh}{\partial\P^{+}}\in[\CS{0}(\partial\P^{+})]^2, \\
  &\quad \ \ \ \restrict{\vvh}{\E}\in\big[\PS{1}(\E)\big]^2\,\,
  \forall\E\in(\partial\P\cap\partial\P^{+}),\
  \restrict{\vvh}{\P^{-}} = \bm{0}
  \Big\}.
\end{align*}
Then, the local virtual element space reads:
\begin{align}
  \label{hhspace}
  \VvhX(\P)
  &\EQUIV \Big\{\,
  \vvhX = (\vvhm + \vvhp) \,:\, \vvhm\in\Vvhm(\P), \, 
  \vvhp\in\Vvhp(\P)
  \Big\}.
\end{align}

\begin{remark}
  The space $\VvhX(\P)$ in \eqref{hhspace} is not a subspace of
  $\HONE(\P)$ as we do not assume any regularity of the virtual
  element functions across the crack, so that a discontinuity is
  admissible.
  This fact implies that also the global virtual element space $\VvhX$
  cannot be a subspace of $\HONE(\Omega)$, but this is not an issue
  since the exact solution contains a discontinuity and thus cannot be
  in $\HONE(\Omega)$.
\end{remark}

An analogous definition of the local virtual element space for
elements cut by a crack can be easily provided also for the enriched
formulation presented in Section \ref{subsec:conforming:VEM:extended}.

As we detail later on, virtual element functions along interface edges
can be reconstructed by a suitable approximation.
We obtain the following representation for the virtual element
approximation on the element $\P$ cut by $\gamma$:
\begin{align}
  \label{eq:XVEM-HH}
  \vvhX(\xv) =
  \sum_{i = 1}^{N_{\textrm{dof}}^{\textrm{VE}}}
  \left[ \phiv_i^{-}(\xv)\vs^{-}_i + \phiv_i^{+}(\xv)\vs^{+}_i \right]
  \quad\forall\xv\in\P ,
\end{align}
where $\vs^{-}_i$ and $\vs^{+}_i$ are the degrees of freedom
associated with $\phiv_i^{-}$ and $\phiv_i^{+}$, respectively.
To provide a feasible solution using~\eqref{eq:XVEM-HH}, it is
necessary to know the trace of the virtual shape functions $\phiv_i$
along the crack.
We also need two distinct \emph{regular} projectors, respectively
$\Pi^{a,-}$ onto $\big[\PS{1}(\P^{-})\big]^2$ and $\Pi^{a,+}$ onto
$\big[\PS{1}(\P^{+})\big]^2$, and two distinct \emph{extended}
projectors, respectively $\Pi^{a,-}_X$ onto $\PSX{}(\P^{-})$ and
$\Pi^{a,+}_X$ onto $\PSX{}(\P^{+})$.
These projection operators must be computable from the
$N_{\textrm{dof}}^{\textrm{HH}}$ nodal degrees of freedom.
A convenient approximation of the trace of the $i$-th virtual element
shape function $\phiv_i$ along the crack is provided by a
vector-valued function $\Nv_i(\xv)$, that is componentwise harmonic
on the cracked element $\P$.
Such a function is built as a first-order polyharmonic
spline~\cite{Duchon:1977:SMR} and the reader is pointed
to~\cite{Benvenuti:2019:CMAME} for further details.

We also point out that the flexibility of the virtual element method allows an element to be cut into two polygonal virtual
elements, regardless of the element shape, and therefore the modeling of crack opening and growth can % could 
follow this alternative route (see~\cite{Nguyen:2018:VEM}). However, 
mesh quality can be affected. For instance, let us consider the case % of one element with a
when partitioning of the element results in one subelement being a quasi-degenerate triangle: this badly-shaped triangle
will worsen matrix-conditioning and/or the interpolation error.
% , and more so in a 3D setting where sliver thetraedra 
This scenario becomes acute in 3D if sliver tetrahedra
appear and the partitioning % is very 
is now much more difficult to handle, both algorithmically and computationally. Moreover, a technique to embed a discontinuous field in the extended virtual element discrete space is required whenever a mesh-independent modeling approach is preferred, % especially in the  domain 
such as in the simulation
of cohesive fracture % modeling 
or when a finite element transitions from a continuous regime to a 
region with discontinuous kinematics~\cite{Benvenuti:2021:CM}.

%% file: sec4_imple_new.tex
\section{Numerical implementation} 
\label{sec:implementation}
In this section, we outline the main implementation aspects of the
extended virtual element method introduced in 
Sections~\ref{subsec:conforming:VEM:extended}--\ref{subsec:disc}. 
%%% \ref{subsec:partial-enrichment}
For the (nonenriched) virtual element
formulation presented in Section~\ref{subsec:conforming:VEM:regular},
the interested reader can refer to~\cite{Wriggers:2017:CM}.

\smallskip

\subsection{Fully enriched elements with singular fields}
\label{subsec:singular_enr}

To begin with, we assume a \emph{fully enriched} element $\P$, i.e., an
element in which all the $\NP$ nodes are enriched with the two
singular fields \eqref{def:scaledcrackmodes}.
Therefore, on such element we have $4\NP$ degrees of freedom, and we
can represent any virtual displacement field $\vvhX\in\VvhX(\P)$ in
terms of the shape functions of $\VvhX(\P)$
as
$\vvhX=\NvX\DOFS{\vvhX}$ where $\DOFS{\vvhX}\in\REAL^{4\NP}$ is the
vector of the degrees of freedom of $\vvhX$ with respect to the basis
function $\{\phiv_i\}_{i=1}^{4\NP}$ spanning $\VvhX(\P)$ and
$\NvX\in\REAL^{2\times 4\NP}$ is the matrix whose columns contain such
basis functions
\begin{equation}
  \label{eq:NvX:def}
  \NvX \EQUIV
  \left[
    \begin{array}{c} \varphi_1         \\ 0                  \end{array} 
    \begin{array}{c} 0                 \\ \varphi_1          \end{array} 
    \begin{array}{c} \dots             \\ \dots              \end{array} 
    \begin{array}{c} \usxIt\varphi_{1}  \\ \usyIt\varphi_{1}  \end{array} 
    \begin{array}{c} \usxIt\varphi_{2}  \\ \usyIt\varphi_{2}  \end{array} 
    \begin{array}{c} \dots             \\ \dots              \end{array} 
    \begin{array}{c} \usxIIt\varphi_{1} \\ \usxIIt\varphi_{1} \end{array} 
    \begin{array}{c} \usxIIt\varphi_{2} \\ \usyIIt\varphi_{2} \end{array} 
    \begin{array}{c} \dots             \\ \dots              \end{array} 
    \right]
  =
  \left[
    \begin{array} {ccccc}
      \phiv_{1} & \phiv_{2} & \dots & \phiv_{4\NP}   
    \end{array}
    \right].
\end{equation}
Now, we define matrix $\MvX\in\REAL^{2\times8}$, whose columns are the
basis vectors $\mv_\alpha$ of $\PSX{}(\P)$ introduced in~\eqref{eq:PX}
\begin{align}\label{eq:M}
  \MvX \EQUIV
  \left[
    \begin{array} {ccccccccccccc}
      1 & 0 & \eta & \xi & 0    & \eta & \usxIt & \usxIIt \\
      0 & 1 & -\xi &   0 & \eta &  \xi & \usyIt & \usyIIt 
    \end{array}
    \right]
  =
  \left[
    \begin{array} {cccccccc}
      \mv_{1} & \mv_{2} & \mv_{3} & \mv_{4} & \mv_{5} & \mv_{6} & \mv_{7} & \mv_{8}
    \end{array}
    \right].
\end{align}
Hereafter, we conveniently use the notation $\mv_{7}=\uvIt$ and
$\mv_{8}=\uvIIt$.
We represent the action of the projection operator $\PiaX$ on the
virtual basis functions by means of a matrix
$\matPiaX\in\REAL^{8\times4\NP}$.
The $i$-th column of this matrix, denoted by
$\piv^{i}=\big(\pi^{i}_{\alpha}\big)\in\REAL^{8}$, contains the
coefficients of $\PiaX(\phiv_i)$ when the projection is expanded on the
basis $\MvX$ so that
\begin{align}
  \label{eq:PI}
  \PiaX(\phiv_i) = \MvX\piv^{i} = \sum_{\alpha=1}^{8}\mv_{\alpha}\pi^{i}_{\alpha},
\end{align}
which in compact form can be expressed as $\PiaX(\NvX)=\MvX\matPiaX$.

We preliminarily observe that for every $\qvX\in\PSX{}(\P)$ and every
$\vvhX\in\VvhX(\P)$, recalling that $\nabla\cdot\sigmav(\qvX)=0$ and
applying the divergence theorem, we find that
\begin{align}\label{eq:diver}
  \asP(\qvX,\vvhX)
  &
  = \int_{\P} \sigmav(\qvX) : \varepsilonv(\vvhX) d\xv
  = \int_{\P} \sigmav(\qvX) : \nabla\vvhX d\xv \nonumber \\
  &
  = \int_{\P} \nabla \cdot (\sigmav(\qvX)\cdot\vvhX) - \vvhX\cdot(\nabla\cdot\sigmav(\qvX)) d\xv \nonumber \\ 
  &
  = \int_{\partial\P} (\sigmav(\qvX)\cdot\vvhX)\cdot\norP ds 
\end{align}
The boundary integral is always computable, since the integrand is
known on the boundary.

By virtue of \eqref{eq:diver} and recalling the definition of the
elliptic projection operator, we compute the projections of the
virtual shape functions in terms of the basis of $\PSX{}(\P)$.
Indeed, for every $\vvhX\in\VvhX(\P)$ we can write the following
orthogonality condition:
\begin{align}\label{eq:def_compact}
  \asP( \mv_\beta, \PiaX(\vvhX) ) = \asP( \mv_\beta, \vvhX ) \qquad
  \beta = 1,\dots,8.
\end{align}
Then, recall that $\DOFS{\vvhX}=(\vshXi)$ are the $4\NP$ degrees of freedom
of $\vvhX$ with respect to the basis $\{\phiv_i\}_{i=1}^{4\NP}$.
In view of \eqref{eq:PI} and noting that $\vvhX$ and $\DOFS{\vvhX}$
are arbitrary, we find that
\begin{equation}
  \label{eq:computepi}
  \begin{array}{lll}
    {}        & \displaystyle\sum_{i=1}^{4\NP}  \asP\big( \mv_\beta, \PiaX(\phiv_i) \big) \vshXi      \,=\, \sum_{i=1}^{4\NP} \asP( \mv_\beta, \phiv_i ) \vshXi &\quad \beta = 1,\dots,8 \nonumber \\[1.25em]
    \implies  & \displaystyle\sum_{\alpha=1}^{8}\asP\big( \mv_\beta, \mv_\alpha       \big) \pi^i_\alpha \,=\, \asP( \mv_\beta, \phiv_i)                                                   &\quad \beta = 1,\dots,8,\,\,\,\is=1,\ldots,4\NP \nonumber \\[1.5em]
    \implies  & \GvXt\matPiaX = \BvXt,
  \end{array}
\end{equation}
where the companion matrices %%
$\GvXt \in \REAL^{8 \times 8}$ and %%
$\BvXt \in \REAL^{8 \times 4\NP}$ are defined componentwise as
\begin{align*}
  (\GsXt)_{\beta,\alpha} &= \asP( \mv_\beta, \mv_\alpha), \quad \beta,\alpha = 1,\dots,8, \\[1.em]
  (\BsXt)_{\beta,i}     &= \asP( \mv_\beta, \phiv_i  ), \quad \beta = 1,\dots,8, \quad i = 1,\dots,4\NP,
\end{align*}
or in the equivalent compact form by
\begin{align}
  \label{eq:G-B:compact}
  \GvXt = \as( \MvX^T, \MvX )
  \qquad\textrm{and}\qquad
  \BvXt = \as( \MvX^T, \NvX ).
\end{align}
Recalling \eqref{eq:diver}, both matrices $\GvXt$ and $\BvXt$ can be
computed by integrating on the element boundary as follows:
\begin{align}\label{int}
  (\GsXt)_{\beta,\alpha} &= \int_{\P} \sigmav(\mv_\beta) : \varepsilonv( \mv_\alpha )d\Omega = \int_{\partial\P} \big(\sigmav(\mv_\beta\big)\cdot \mv_\alpha\big) \cdot \norP  d\Gamma, \\
  (\BsXt)_{\beta,i}     &= \int_{\P} \sigmav(\mv_\beta) : \varepsilonv( \phiv_i   )d\Omega = \int_{\partial\P} \big(\sigmav(\mv_\beta\big)\cdot \phiv_i  \big) \cdot \norP d\Gamma.
\end{align}
The first three rows of $\GvXt$ and $\BvXt$ are zero, since the small
strain tensor associated to rigid body motions is zero, and therefore
$\GvXt$ is rank deficient.
To overcome this issue we use
conditions~\eqref{eq:cond_extra_xvem:B}-\eqref{eq:cond_extra_xvem:C},
which imposes that the projector preserves the average nodal
translations and rotations.
So, we define the matrices $\GvX=\big((\GsX)_{\beta,\alpha}\big)$ and
$\BvX=\big((\BsX)_{\beta,i}\big)$ as
\begin{align}
  (\GsX)_{\beta,\alpha} &= 
  \left \{
  \begin{array} {ll}
    \displaystyle\frac{1}{\NP}\sum_{j=1}^{\NP} \mv_\alpha(\vv_j)                 & \beta=1,2, \\[1.5em]
    \displaystyle\frac{1}{\NP}\sum_{j=1}^{\NP} \rv(\vv_j)\cdot\mv_\alpha(\vv_j)  & \beta=3,   \\[1.8em]
    (\GsXt)_{\beta,\alpha}                                                      & \beta=4,\dots,8,
  \end{array}
  \right.
  \label{G_matrix}
  \intertext{and}%\\[1em]
  (\BsX)_{\beta,i} &= 
  \left \{
  \begin{array} {ll}
    \displaystyle\frac{1}{\NP}\sum_{j=1}^{\NP} \phiv_i(\vv_j)                  & \beta=1,2, \\[1.5em]
    \displaystyle\frac{1}{\NP}\sum_{j=1}^{\NP} \rv(\vv_j)\cdot\phiv_i(\vv_j)   & \beta=3,   \\[1.8em]
    (\BsXt)_{\beta,i}                                                         & \beta=4,\dots,4\NP.
  \end{array}
  \right .
  \label{B_matrix}
\end{align}
Since matrix $\GvX$ is nonsingular, the projection matrix $\matPiaX$
is the unique solution of the linear system $\GvX\matPiaX=\BvX$.
To derive the representation of the operator $\PiaX$ with respect to
the basis $\{\phiv_i\}_{i=1}^{4\NP}$ spanning $\VvhX(\P)$ we introduce
matrix $\DvX\in\REAL^{4\NP\times8}$, whose $\alpha$-th column
($\alpha=1,\ldots,8$) contains the degrees of freedom of the vector
polynomial $\mv_{\alpha}$.
Therefore, it holds that $\MvX=\NvX\DvX$ and
\begin{align*}
  \PiaX(\NvX) = \MvX\matPiaX = \NvX\,\big(\DvX\matPiaX\big),
  \label{Pi_compute}
\end{align*}
from which we infer that such matrix representation is given by matrix
$\DvX\matPiaX$.
A straightforward calculation yields that $\GvX=\BvX\DvX$ in the \XVEM{}, which is similar in form to the standard relation in the VEM, $\Gv=\Bv\Dv$ ~\cite{BeiraodaVeiga-Brezzi-Marini-Russo:2014}. This
provides a means to verify the correctness of the computation of these
matrices.

The stiffness matrix is given by the sum of a consistency and
stability term,
\begin{equation*}
  \KvPX = \KvPXc + \KvPXs,
\end{equation*}
so that we can evaluate the extended stiffness bilinear form applied
to $\vvhX,\,\wvhX\in\VvhX(\P)$ by using the degrees of freedom
$\DOFS{\vvhX}$ and $\DOFS{\wvhX}$ as follows
\begin{align}
  \ashPX\Big( \vvhX, \wvhX \Big) = \big(\DOFS{\vvhX}\big)^T\,\KvPX\,\DOFS{\wvhX}.
  %\label{eq:KvP:def:2}
\end{align}
For every $\vvhX\in\VvhX(\P)$, we first consider the relations
\begin{equation}\label{Pi_relations}
  \PiaX(\vvhX)
  = \PiaX\big( \NvX\DOFS{\vvhX} \big)
  = \PiaX\big( \NvX \big) \DOFS{\vvhX}
  = \MvX\matPiaX\,\DOFS{\vvhX}
  = \big(\DOFS{\vvhX}\big)^T (\matPiaX)^T (\MvX)^T.
\end{equation}
Recalling \eqref{Pi_relations}, we compute the consistency term as
\begin{align*}
  \asP\bigl( \PiaX( \vvhX ), \PiaX( \wvhX )\bigr)
  &= \big(\DOFS{ \wvhX }\big)^T (\matPiaX)^T \asP\bigl( \MvX^T, \MvX \bigr) \matPiaX \DOFS{\vvhX} \nonumber\\[0.5em]
  &= \big(\DOFS{ \wvhX }\big)^T (\matPiaX)^T \GvXt \matPiaX \DOFS{ \wvhX }.
\end{align*}
By comparison, we see that
\begin{align}
  \KvPXc = (\matPiaX)^T \GvXt \matPiaX.
  \label{consistency}
\end{align}

For the stability term we generalize the so-called  \emph{dofi-dofi}~\cite{BeiraodaVeiga-Brezzi-Cangiani-Manzini-Marini-Russo:2013}
and \emph{D-recipe}~\cite{Beirao-Dassi-Russo:2017} stabilizations by evaluating the second term
in~\eqref{eq:ashPX:def} at the element vertices.
To this end, we first note that
\begin{align*}
  \vvhX(\xv_{\ell})-\PiaX(\vvhX)(\xv_{\ell})
  = \NvX(\xv_{\ell})\DOFS{ \wvhX } - \NvX(\xv_{\ell})\DvX\matPiaX\DOFS{ \wvhX }.
  %%  = \Big(\Jv - \DvXt\matPiaX\Big)_{\ell}\hat{ \wv }_X^{\hh},
\end{align*}
Definition~\eqref{eq:NvX:def} and $\varphi_j(\xv_{\ell}) = \delta_{j\ell} $
implies that
\begin{align*}
  \mbox{$j$-th~column~of~}\NvX(\xv_{\ell})
  &=
  \begin{cases}
    \Big( \varphi_j(\xv_{\ell}),\,                   0                                  \Big)^T & \mbox{for~$1     \leq\js\leq2\NP,\,j$ odd},  \\[0.5em]
    \Big( 0,                                  \varphi_j(\xv_{\ell})                   \Big)^T & \mbox{for~$1     \leq\js\leq2\NP,\,j$ even}, \\[0.5em]
    \Big( \usxIt (\xv_{\ell})\phi_j(\xv_{\ell}),\,  \usyIt(\xv_{\ell})\varphi_j(\xv_{\ell})  \Big)^T & \mbox{for~$1+2\NP\leq\js\leq3\NP$},          \\[0.5em]
    \Big( \usxIIt(\xv_{\ell})\phi_j(\xv_{\ell}),\, \usyIIt(\xv_{\ell})\varphi_j(\xv_{\ell}) \Big)^T & \mbox{for~$1+3\NP\leq\js\leq4\NP$},
  \end{cases}
  \\[1em]
  &=
  \begin{cases}
    \Big( \delta_{j,\ell},                      0                               \Big)^T & \mbox{for~$1     \leq\js\leq2\NP,\,j$ odd},      \\[0.5em]
    \Big( 0,                                  \delta_{j,\ell}                   \Big)^T & \mbox{for~$1     \leq\js\leq2\NP,\,j$ even},     \\[0.5em]
    \Big( \usxIt (\xv_{\ell})\delta_{j,\ell},\,  \usyIt (\xv_{\ell})\delta_{j,\ell} \Big)^T & \mbox{for~$1+2\NP\leq\js\leq3\NP$}, \\[0.5em]
    \Big( \usxIIt(\xv_{\ell})\delta_{j,\ell},\,  \usyIIt(\xv_{\ell})\delta_{j,\ell} \Big)^T & \mbox{for~$1+3\NP\leq\js\leq4\NP$},
  \end{cases}
\end{align*}
for $\ell=1,\ldots,\NP$.
We collect $\NvX(\xv_{\ell})$ in the compact block-diagonal matrix
$\Jv=\textrm{diag}(\Jv_{11},\Jv_{22})\in\REAL^{4\NP\times4\NP}$ such
that $\Jv_{11}=\Iv_{2\NP\times2\NP}$, which is the
$2\NP\times2\NP$-sized identity matrix, and
$\Jv_{22}=\big[\Jv_{22}^{\INTP},\Jv_{22}^{\INTP\INTP}\big]$, with
\begin{align*}
  \Jv_{22}^{\INTP}     &= \textrm{diag}\Big( \uvIt (\xv_{1}), \uvIt (\xv_{2}), \ldots, \uvIt (\xv_{\NP}) \Big),\\[0.5em]
  \Jv_{22}^{\INTP\INTP} &= \textrm{diag}\Big( \uvIIt(\xv_{1}), \uvIIt(\xv_{2}), \ldots, \uvIIt(\xv_{\NP}) \Big).
\end{align*}
Here, $\Jv_{11}\in\REAL^{2\NP\times2\NP}$ collects the values of the
first $2\NP$ columns of $\NvX(\xv_{\ell})$ and is such that the
$\ell$-th row corresponds to the $\ell$-th element vertex.
In turn, the matrix blocks $\Jv_{22}^{\INTP}$ and
$\Jv_{22}^{\INTP\INTP}$ are $2\NP\times\NP$-sized matrices that
collect the $\NP$ columns of $\NvX$ corresponding to $\psivI\Vvha(\P)$
and $\psivII\Vvha(\P)$, and again each row corresponds to a given
element vertex.
Finally, we introduce the matrix $\DvXt=\Jv\DvX$, and we write the
\emph{dofi-dofi} stabilization as
\begin{align}
  \SPX\Big( \vvhX-\PiaX(\vvhX), \wvhX-\PiaX(\wvhX) \Big)
  &= \tau \sum_{\ell=1}^{\NP}
  \Big( \NvX(\xv_{\ell})\DOFS{ \vvhX } - \NvX(\xv_{\ell})\DvX\matPiaX\DOFS{ \vvhX } \Big) \nonumber\\[-0.5em] &\phantom{= \alpha\sum_{\ell=1}^{\NP}}\,\times
  \Big( \NvX(\xv_{\ell})\DOFS{ \wvhX } - \NvX(\xv_{\ell})\DvX\matPiaX\DOFS{ \wvhX } \Big)
  \nonumber\\[0.5em]
  &= \tau
  \Big( \DOFS{ \vvhX }  \Big)^T
  \Big( \Jv - \DvXt\matPiaX \Big)^T \,
  \Big( \Jv - \DvXt\matPiaX \Big)
  \Big( \DOFS{ \wvhX }      \Big) ,
  \label{eq:stab2}
\end{align}
where $\tau$ is a suitable scaling parameter; a possible choice is $\tau = \alpha \, \textrm{trace}(\KvPXc)\slash{4\NP}$, where $\alpha$ is a constant (a sensitivity analysis on the choice of $\alpha$ is presented in the next section). Hence,
\begin{align*}
  \KvPXs=\tau(\Jv-\DvXt\matPiaX)^T\,(\Jv-\DvXt\matPiaX).
\end{align*}
Similarly, the \emph{D-recipe} stabilization, which was originally
proposed in~\cite{Beirao-Dassi-Russo:2017} for the Poisson equation, can be
generalized by taking the stabilization matrix
\begin{align*}
  \KvPXs=(\Jv-\DvXt\matPiaX)^T\,\SvX\,(\Jv-\DvXt\matPiaX),
\end{align*}
where
\begin{align*}
  (\SvX)_{ij} = \delta_{ij} \max\left( \frac{\textrm{trace}(\bm{C})}{3}h_E, (\KvPXc)_{ii} \right)
  \quad i,j=1,\ldots,\NP,
\end{align*}
and with the same definition of $\Jv$ and $\DvXt$ introduced above.

\subsection{Partially enriched elements with singular fields}

Let $\EvX\in\REAL^{4\NP\times(2\NP+2\kP)}$ be the matrix
representation of the extension operator $\extkP$ introduced in
Section~\ref{subsec:partial-enrichment}, so that every vector-valued
field $\vvhXt\in\VvhXt(\P)$ with degrees of freedom
$\DOFS{\vvhXt}\in\REAL^{2\NP+2\kP}$ is remapped into the vector-valued
field $\vvhX\in\VvhX(\P)$ with degrees of freedom
$\DOFS{\vvhX}=\EvX\DOFS{\vvhXt}\in\REAL^{4\NP}$.
Let $\KvPXt\in\REAL^{(2\NP+2\kP)\times(2\NP+2\kP)}$ such that
\begin{align}
  \ashPXt\Big( \vvhXt, \wvhXt \Big) =
  \big(\DOFS{\vvhXt}\big)^T\,\KvPXt\,\DOFS{\wvhXt}.
  \label{eq:KvP:def:2}
\end{align}
Now, starting from definition~\eqref{eq:ashPXt:def} and
using~\eqref{eq:KvP:def:2}, a straightforward calculation yields
\begin{align*}
  \ashPXt\Big( \vvhXt, \wvhXt \Big)
  &= \ashPX \Big( \extkP\vvhX, \extkP\wvhXt \Big)
  = \Big(\DOFS{ \extkP\wvhX }\Big)^T\,\KvPX\,\DOFS{ \extkP\vvhX }
  \\[0.5em]
  &= \Big(\EvX\DOFS{\wvhXt}\Big)^T\,\KvPX\,\EvX\DOFS{\vvhXt}
  = \Big(\DOFS{\wvhXt}\Big)^T\,\big(\EvX\big)^T\KvPX\EvX\,\DOFS{\vvhXt},
\end{align*}
and by comparison with~\eqref{eq:KvP:def:2} it follows that
\begin{align}
  \KvPXt = \big(\EvX\big)^T\KvPX\EvX,
  \label{eq:KvPt:def}
\end{align}
since $\vvhXt$ and $\wvhXt$ are arbitrary.
To conclude this section, we are only left to explain the construction
of the matrix $\EvX$ that embodies definition~\eqref{eq:extkP:def}.
To obtain such a matrix, we take the $4\NP\times4\NP$-size identity
matrix and remove the $2(\NP-\kP)$ columns that corresponds to the
basis functions of the nonenriched vertices in
$\psivI\Vvha(\P)\oplus\psivII\Vvha(\P)$.
Finally we note that when we apply $(\EvX)^T$ to the left of matrix
$\KvPX$ and $\EvX$ to the right of matrix $\KvPX$ we are indeed
selecting the rows and the columns of $\KvPX$ that corresponds to all
the degrees of freedom of $\Vvh(\P)$ and the degrees of freedom of the
enriched nodes of $\psivI\Vvha(\P)\oplus\psivII\Vvha(\P)$.
In the numerical implementation we do not need to construct the
matrix $\EvX$ explicitly and compute $\KvPXt$
using~\eqref{eq:KvPt:def}, since we can simply build the stiffness matrix
$\KvPX$ of the full enrichment case and remove all rows and columns that refer to non-enriched nodal
degrees of freedom.

\subsection{Embedding discontinuities}
\label{sec_discontin}
\begin{figure}[!htb]
  \centering
  \includegraphics[width=0.5\textwidth]{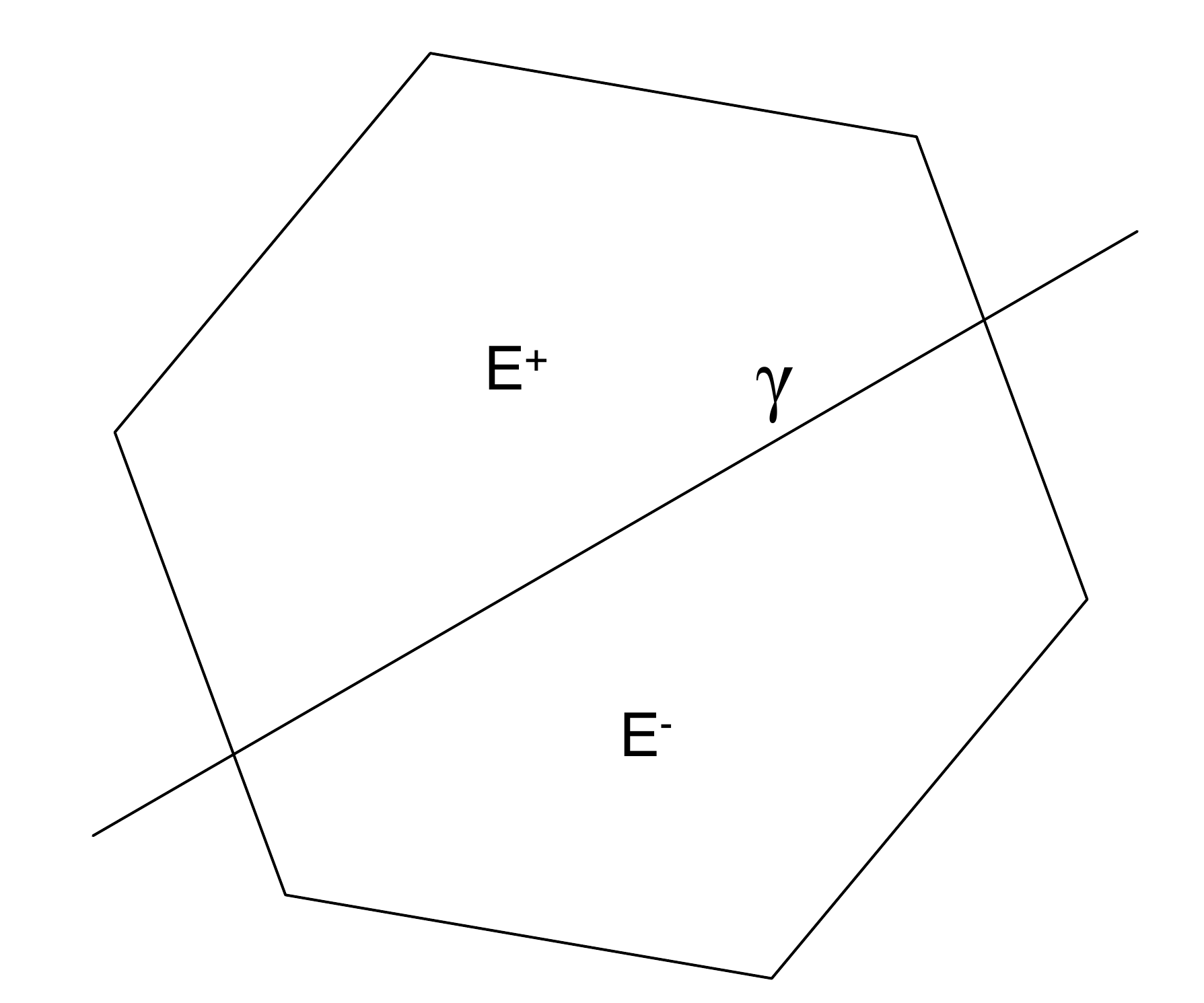}
  \caption{Crack line $\gamma$ cuts element $E$ 
                            into two subelements $E^-$ and $E^+$.}\label{fig:cracked_poly}
\end{figure}

Let us consider the case in which an element $E$, fully enriched
according to the construction outlined in Section
\ref{subsec:singular_enr}, is also cut by a crack $\gamma$ into two
subelements $E^{-}$ and $E^{+}$, see Fig.~\ref{fig:cracked_poly}.
Then, following the approach presented in Section \ref{subsec:disc},
we can compute the projectors ${\matPiaX}^{-}$ and
${\matPiaX}^{+}$ on the two subelements $E^{-}$ and $E^{+}$
generated by the crack line $\gamma$ intersecting element
$E$.
Recalling that, in the present approach, the number of degrees of
freedom is doubled and is equal to $N_{\textrm{dof}}^{\textrm{HH}} =
8N_E$, we denote by $i^{-}$ and $i^{+}$ the degrees of freedom
associated to $\bm{\phi}_i^{-}$ and $\bm{\phi}_i^{+}$ respectively,
with $ i = 1,\dots, 4N_E$.

Then, the consistent part of the stiffness matrix has a diagonal block
structure and is composed by the following submatrix blocks:
\begin{subequations}
  \begin{align*}
    (\KvPXc)_{+,+} &=  ({\matPiaX}^{+})^T {\GvXt}^{+} {\matPiaX}^{+}, \\[0.5em]
    (\KvPXc)_{-,-} &=  ({\matPiaX}^{-})^T {\GvXt}^{-} {\matPiaX}^{-}, %\\[0.5em]
    %% (\bm{K}^{\P}_c)_{-,+} &= \bm{0},
  \end{align*}
\end{subequations}
where matrices ${\GvXt}^{-}$ and ${\GvXt}^{+}$ are the
counterparts of $\GvXt$ computed for $E^{-}$ and $E^{+}$,
respectively.

On the other hand, the general expression for the stabilization part
shares the same diagonal block structures and reads:
\begin{subequations}
  \begin{align*}
    (\KvPXs)_{+,+} &= \tau (\Jv-\DvXt^{+}{\matPiaX}^{+})^T (\Jv - \DvXt^{+}{\matPiaX}^{+}), \\[0.5em]
    (\KvPXs)_{-,-} &= \tau (\Jv-\DvXt^{-}{\matPiaX}^{-})^T  (\Jv - \DvXt^{-}{\matPiaX}^{-}), %%\\[0.5em]
    %% (\bm{K}^{\P}_s)_{-,+} &= \bm{0},
  \end{align*}
\end{subequations}
where matrices $\DvXt^{-}$ and $\DvXt^{+}$ are the counterparts
of $\DvXt$ computed for $E^{-}$ and $E^{+}$, respectively.

\subsection{Computation of stress intensity factors}
\label{SIFS}
In order to determine the susceptibility of a given elastic
two-dimensional body to fracture growth we need to extract appropriate
crack tip parameters such as the $J$-integral and mixed-mode stress
intensity factors.
\begin{figure}[!htb]
  \centering
  \begin{tabular}{c}
    \begin{subfigure}{0.48\textwidth}
      \centering
      \includegraphics[width=\textwidth]{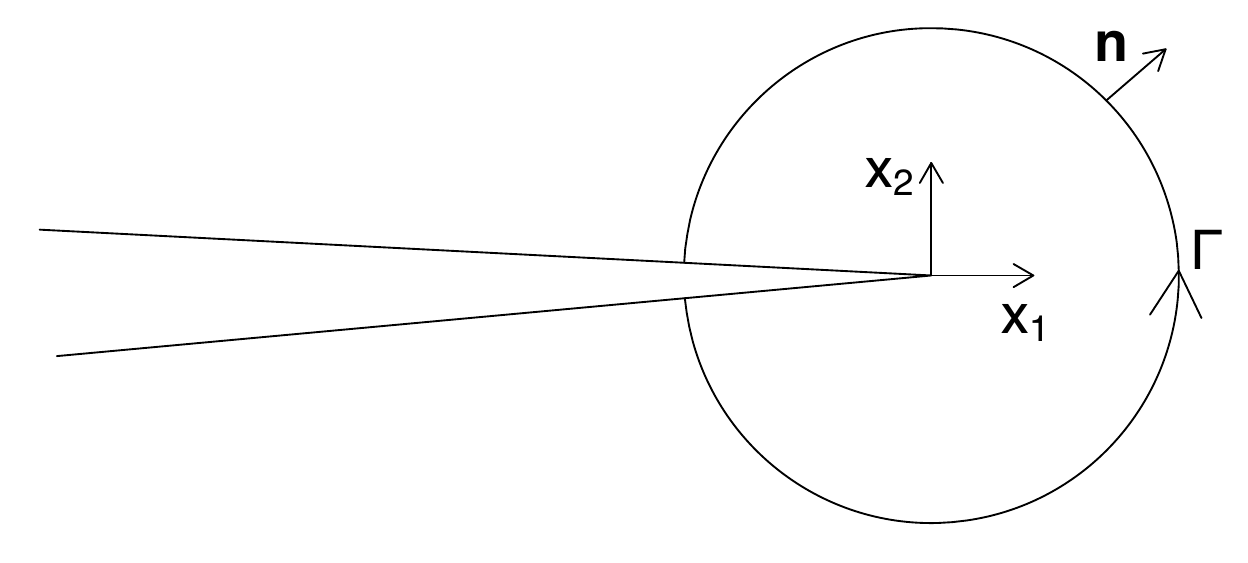}
      \caption{}\label{crack_tip}
      \vspace{5mm}
    \end{subfigure}\hfill \\
    \begin{subfigure}{0.48\textwidth}
      \centering
      \includegraphics[width=\textwidth]{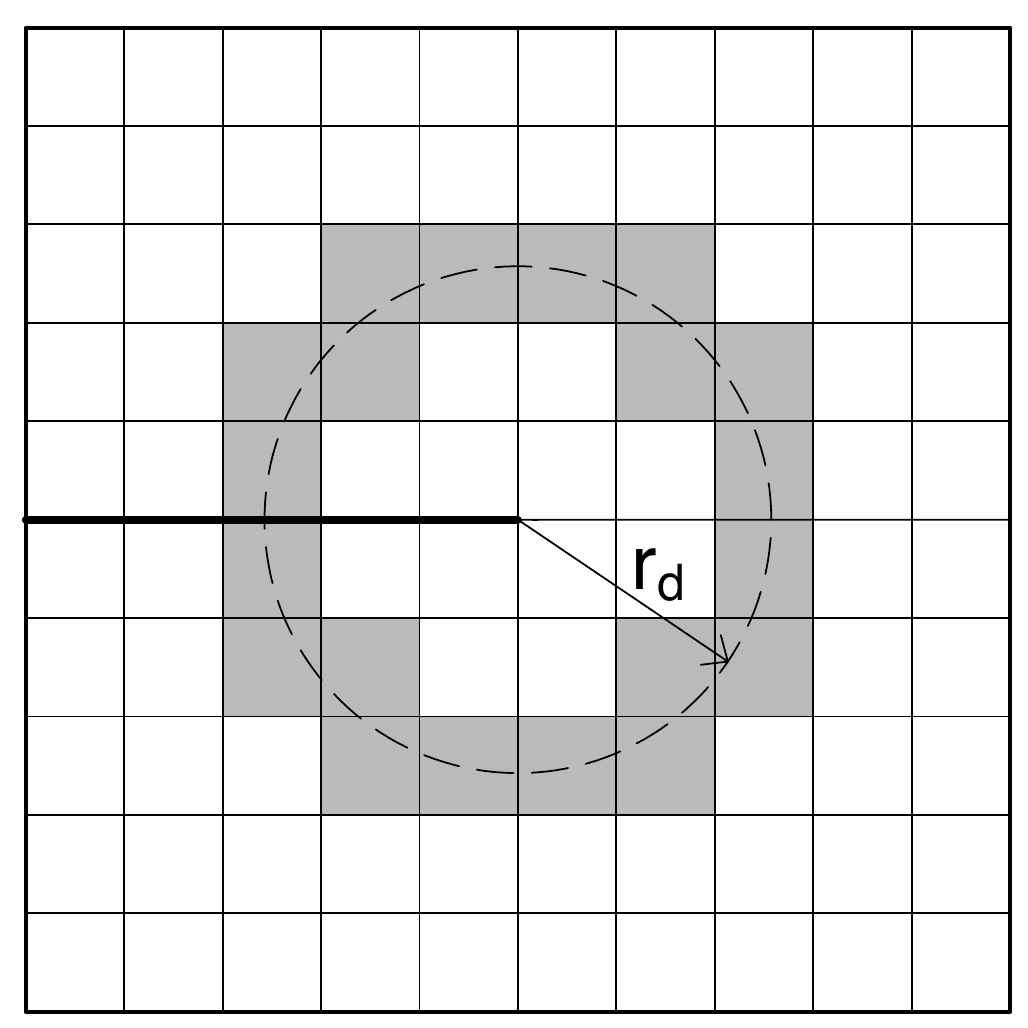}
      \caption{}\label{ring}
    \end{subfigure}\\
  \end{tabular}
  \caption{Local crack tip coordinates (a) and discretized J-domain
    (shaded area) (b).}\label{fig:SIFS}
\end{figure}
We consider a neighborhood of the crack tip, as shown in
Fig.~\ref{crack_tip}.
Given an arbitrary closed path $\Gamma$ around the crack tip, the
$J$-integral is defined as:
\begin{equation}
  \label{J_integral}
  J = \int_\Gamma \left ( W dx_2 - T_i \frac{\partial u_i}{\partial x_1}ds   \right ),
\end{equation}
which is path independent under the assumptions of small deformations,
elastic material behavior and quasi-statically applied 
loads~\cite{Rice:1968}.
In \eqref{J_integral}, $u_i$ is the $i$-th component of the
displacement field, $ds$ is the differential of the arc length of
$\Gamma$, $T_i$ is the $i$-th component of the traction vector along
$\Gamma$ and $W$ is the strain energy density, which is defined as:
\begin{equation}
  \label{W}
  W = \int_0^{\varepsilon_{ij}} \sigma_{ij} d\varepsilon_{ij} =
   \frac{1}{2}  \sigma_{ij} \varepsilon_{ij}.
\end{equation}
However, the $J$-integral in~\eqref{J_integral} is not 
well-suited for
numerical computations, since it is defined on a vanishingly small
closed path.
For this reason, in numerical procedures, Li et al.~\cite{Li:1985} proposed to recast the line integral \eqref{J_integral}
into a domain integral over an annular region $\Omega$, bounded by an
inner closed curve $\partial \Omega_i \to 0$, which contains the crack
tip, and an outer closed curve $\partial \Omega_o$:
\begin{equation}
  \label{J_integral_Li}
  J = \int_\Omega \left ( \sigma_{ij}\frac{\partial u_i}{\partial x_1} - W \delta_{1i} \right )\frac{\partial w}{\partial x_i}d\Omega,
\end{equation}
where $w$ is a suitable weight function that is equal to unity within
the domain bounded by $\partial \Omega_o$ and vanishes on $\partial
\Omega_o$.
Based on this definition, given two equilibrium states denoted by
superscripts (1) and (2), referred to as the present state and an
auxiliary state respectively, the interaction
integral in domain form is given by
\begin{equation}
  \label{interaction_int}
  I^{(1,2)} = \int_\Omega F_j(x_1,x_2)\frac{\partial w}{\partial x_j}d\Omega,
\end{equation}
where
\begin{equation}
  \label{F_int}
  F_j(x_1,x_2) = \sigma_{ij}^{(1)}\frac{\partial u_i^{(2)}}{\partial x_1} + \sigma_{ij}^{(2)}\frac{\partial u_i^{(1)}}{\partial x_1} - W^{(1,2)} \delta_{1j},
\end{equation}
and $W^{(1,2)} = \sigma_{ij}^{(1)} \varepsilon_{ij}^{(2)} =
\sigma_{ij}^{(2)} \varepsilon_{ij}^{(1)} $ is the interaction strain
energy.

The stress intensity factors for mode I and mode II crack opening,
respectively denoted by $K_I$ and $K_{II}$, are computed using the relation
\begin{equation}
  \label{SIFs_int}
  I^{(1,2)} = \frac{2}{E'}\left [K_I^{(1)}K_I^{(2)} + K_{II}^{(1)}K_{II}^{(2)} \right ],
\end{equation}
where $E' = E$ for plane stress conditions and $E' = E/(1-\nu^2)$ for
plain strain conditions.
Indeed, on choosing the auxiliary field corresponding to $K_I = 1$ and
$K_{II} = 0$ allows $K_I$ to be extracted in a straightforward manner
and similarly on selecting the auxiliary field corresponding to $K_I =
0$ and $K_{II} = 1$ allows $K_{II}$ to be computed:
\begin{equation}
  \label{SIFs_comput}
  K_I = \frac{E'}{2}I^{(1,I)}, \quad K_{II} = \frac{E'}{2}I^{(1,II)}.
\end{equation}

However, computing the interaction integral \eqref{interaction_int} is
not straightforward in the X-VEM, since the numerical integration is performed over polygonal elements.
For this reason, after considering a $J$-domain that is
an annular region $\Omega_J$ that consists of a ring of elements that are
intersected by a circle of given radius $r_d$ centered on the crack
tip (i.e., the shaded area in Fig.~\ref{ring}), we apply the
divergence theorem and transform the domain integral
\eqref{interaction_int} into a line integral that is evaluated on the
boundaries of the element~\cite{Nguyen:2018:VEM}:
\begin{equation}
\label{interac_comput}
I^{(1,2)} = \sum_{E \in \Omega_J} \left ( \int_{\partial E} F_j(x_1,x_2)wn_j d\Gamma - \int_{E} \frac{\partial F_j}{\partial x_j}(x_1,x_2)w d\Omega \right ).
\end{equation}
We note that $\nabla \cdot \bm{F} = 0$ % the second integral
in~\eqref{interac_comput}
since the auxiliary fields are equilibrated, and therefore only the boundary integral needs 
to be computed.

Since virtual shape functions are not known in the interior of the
elements, we use the elliptic projection of the solution in terms of
displacements to compute the corresponding deformation field and the
stress components.
Hence, the interaction integral can be finally computed as:
\begin{equation}
  \label{interac_comput2}
  I^{(1,2)} = \sum_{E \in \Omega_J} \int_{\partial E} \left[  \sigma_{ij}(\Pi^a_{E}(u_i^{(1)}))\frac{\partial u_i^{(2)}}{\partial x_1} +\sigma_{ij}^{(2)}\frac{\partial \Pi^a_{E}(u_i^{(1)})} {\partial x_1} -  \widetilde{W}^{(1,2)} \delta_{1j}\right ] w n_j d\Gamma,
\end{equation}
where $\widetilde{W}^{(1,2)} = \sigma_{ij}(\Pi^a_{E}(u_i^{(1)})) \varepsilon_{ij}^{(2)}$. From a computational viewpoint, it is convenient to assume the weight
function $w$ to be equal to unity on all nodes in $\Omega_J$ that lie
within the circle of radius $r_d$, and equal to zero on all nodes
in $\Omega_J$ that lie outside the circle of radius $r_d$.
Along element edges, where integrations are carried out, linear
interpolation of $w$ between its nodal values is adopted.

%% file: sec5_numer.tex
% Hey Emacs, this is -*-latex-*-

\section{Numerical examples}
\label{sec5:numerical_examples}

In order to check the consistency of the X-VEM, we first conduct two
distinct patch tests: an \emph{extended patch test},
addressing the enrichment with singular fields as described in
Section~\ref{subsec:singular_enr}, and a \emph{discontinuous patch test} aimed
at assessing the inclusion of discontinuities in the discrete space by
means of the approach presented in Section~\ref{subsec:disc}.
Then, we test the X-VEM on a benchmark problem to establish the
convergence rate of the method and the accuracy of the stress intensity
factors. Unless stated otherwise,
Young's modulus $E = 10^5$ and Poisson ratio $\nu = 0.3$ are chosen in the numerical computations.
\subsection{Extended patch test}
\label{ext_pt}
The extended patch test ensures that the singular enrichment 
fields in~\eqref{def:scaledcrackmodes} can be exactly reproduced using the
X-VEM.
\begin{figure}[!htb]
  \centering
  \begin{subfigure}{0.45\textwidth}
    \centering
    \includegraphics[width=\textwidth]{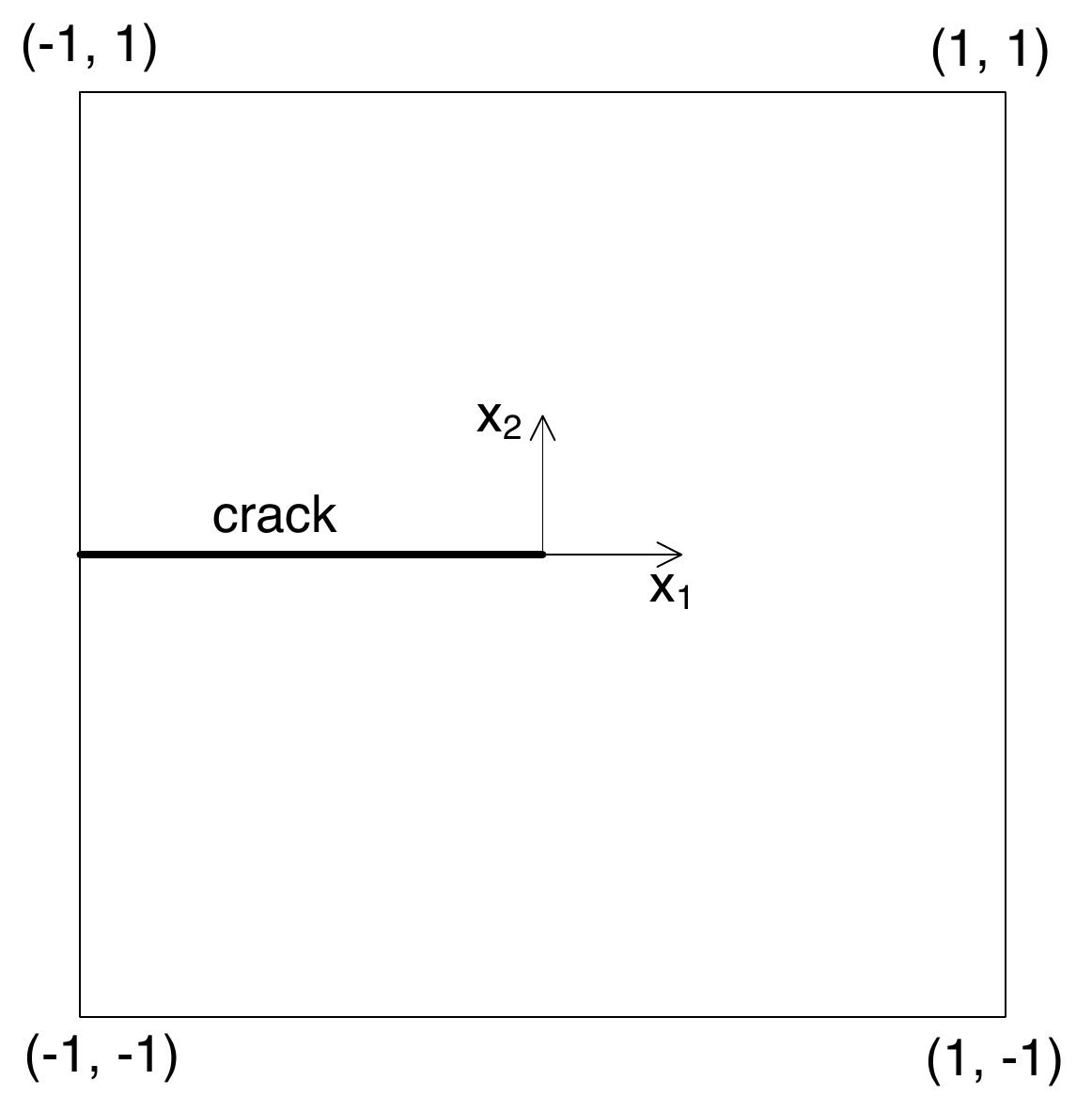}
    \vspace{1mm}
    \caption{}\label{bench_geometry}
  \end{subfigure}\hfill
  \begin{subfigure}{0.52\textwidth}
    \centering
    \includegraphics[width=\textwidth]{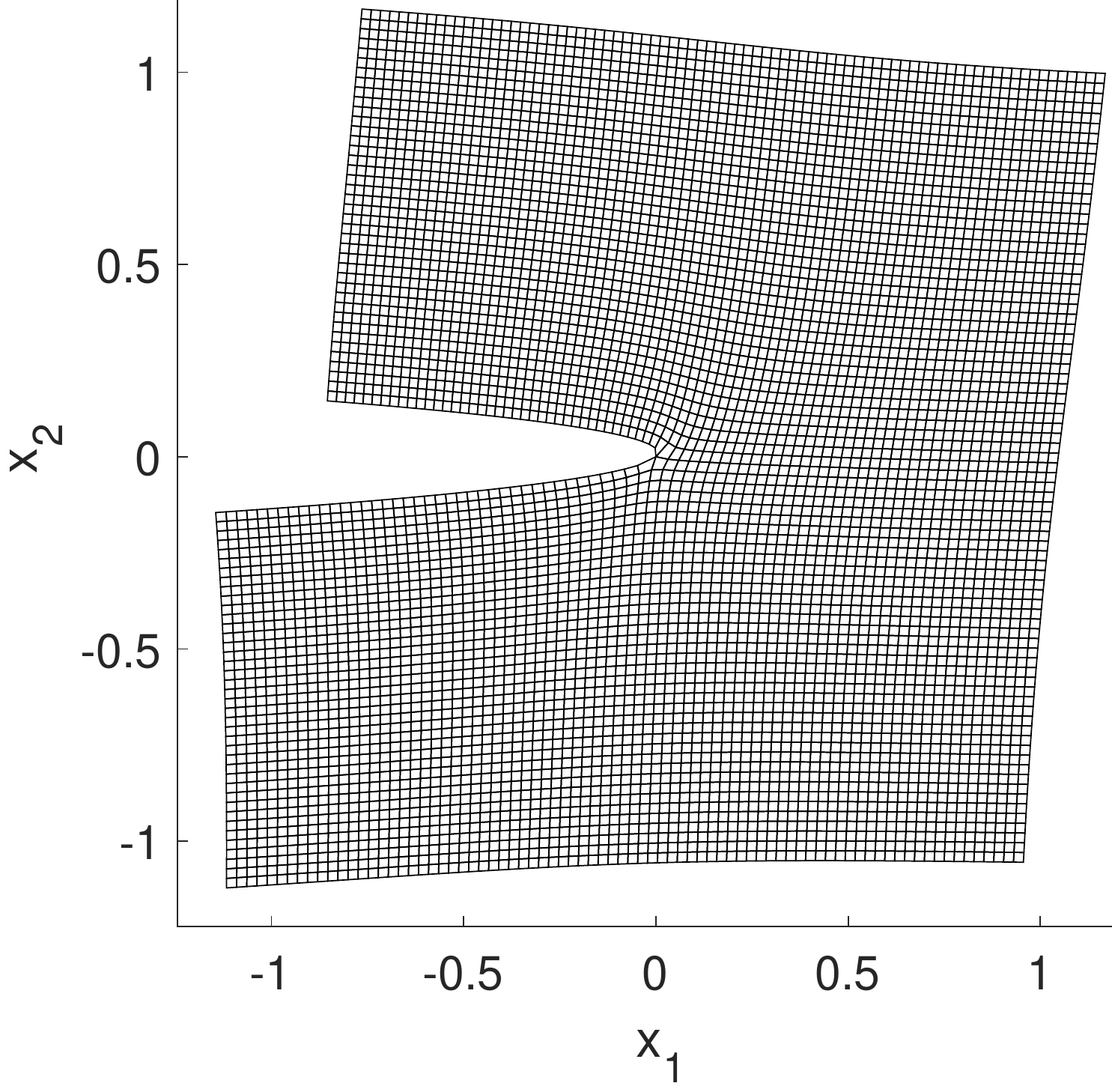}
    \caption{}\label{bench_deformed}
  \end{subfigure}
  \caption{Mixed mode I and mode II crack opening benchmark problem:
    domain geometry (a) and exact deformed shape
    (b).}
  \label{fig:bench_geometry_tot}
\end{figure}
To perform the extended patch test, we consider a square elastic plate that occupies the region $(-1,1)^2$ under plane strain conditions, 
with a horizontal crack of unit length that extends from $(-1,0)$ to
$(0,0)$ (see Fig.~\ref{bench_geometry}).

Both a coarse mesh of $10\times10$ square elements and a coarse mesh
of 64 polygonal elements are considered. For the purposes of the
extended patch test, the crack is modeled explicitly so that we do not
have to embed the discontinuity in the discrete space.
All the nodes in the domain are enriched and the Cartesian components
of the near-tip displacement fields~\eqref{def:scaledcrackmodes} are
imposed on the boundary of the domain by requiring that all the
enriched boundary degrees of freedom are equal to 1 and all the
standard boundary degrees of freedom are equal to 0.
The exact displacement solution field is shown in
Fig.~\ref{bench_deformed}.
As detailed in the previous section, integrals need to be evaluated
over the element boundary only. We adopt a 16-points Gauss quadrature
rule on each element edge.

As a measure for the error of the numerical solution with respect to
the exact solution we adopted the relative error in strain energy,
which is computed as
\begin{equation}
  E(\bm{u}^h) = \frac{|a(\tilde{\bm{u}},\tilde{\bm{u}})-a(\bm{u}^h,\bm{u}^h)|}{a(\bm{u},\bm{u})},
\end{equation}
where 
$\frac{1}{2} a(\bm{u},\bm{u}) = 1.6776885579 \times 10^{-5}$ % 14828 
%% $a(\bm{u},\bm{u})/ 2 = 1.677688557914828e-05$ 
is the strain
energy of the exact solution $\bm{u}$, and $\tilde{\bm{u}}^h$ is the
projection of the discrete solution $\bm{u}^{h}$, which is defined as:
\begin{equation}
  \label{eq:approx_sol}
  \tilde{\bm{u}}^h = \sum_{E \in \mathcal{T}}  \Pi_{E}^a \bm{u}^{h}.
\end{equation}
We also adopt this error measure in the subsequent sections that follow.
In \eqref{eq:approx_sol}, we use the same symbol $\Pi_{E}^a$ to denote
the restriction of the virtual element functions defined on the
element $E$ of the projection operator $\Pi^a$ if $E$ is a nonenriched
element and the projection operator $\Pi^a_X$ if $E$ is an enriched
element.
The choice of using the projection $\tilde{\bm{u}}^h$ of the solution
$\bm{u}^h$ follows from observing that it is not possible to compute
the true energy associated with $\bm{u}^h$, since the virtual
functions are not explicitly
known~\cite{BeiraodaVeiga-Chernov-Mascotto-Russo:2018}.
The relative error in strain energy for the extended patch tests 
is provided in Table~\ref{tab:1}, which clearly shows that the X-VEM delivers sound accuracy in reproducing the enrichment fields, although the error is affected by numerical integration of singular functions.

\begin{table}
  \centering
  \caption{Relative error in strain energy for the extended patch test
           on the $(-1,1)\times(1,1)$ square domain with horizontal crack.}
  \label{tab:1}
  \begin{tabular}{|c|c|}
    \hline
     \multicolumn{1}{c}{\textrm{Mesh}} & \multicolumn{1}{c}{$E(\bm{u}^h)$} \\ \hline
    \multicolumn{1}{c}{$10 \times 10$ square elements}  & \multicolumn{1}{c}{$2 \times 10^{-12}$} \\
    \multicolumn{1}{c}{64 polygonal elements}   
    & \multicolumn{1}{c}{$3 \times 10^{-10}$} \\  \hline
  \end{tabular}
\end{table}

\subsection{Discontinuous patch test}

In order to evaluate the effectiveness and robustness of the X-VEM in
the presence of discontinuities, formulated according to the approach
presented in Section \ref{subsec:disc}, we adopt a suitable patch test
which entails solving a problem whose exact solution is discontinuous
and lies in the discrete space.
We then verify if the extended virtual element approximation matches
such a solution.
To this end, we here adapt the discontinuous patch test first proposed by
Dolbow and Devan~\cite{Dolbow:2004:EEA} in finite strain elasticity to
the present context of plain strain linear elasticity.
The test involves solving the problem of a 2D elastic domain occupying
the unit square domain $\Omega = (0,1)^2$ that is bisected by an horizontal
crack $\gamma$ into two open subdomains $\Omega^- = (0,1)\times (0,1/2)$ and
$\Omega^+ = (0,1)\times(1/2,1)$.
The crack is implicitly included in the model following the
construction proposed in Section~\ref{subsec:disc}.
\begin{figure}
  \centering
  \includegraphics[width=0.6\textwidth]{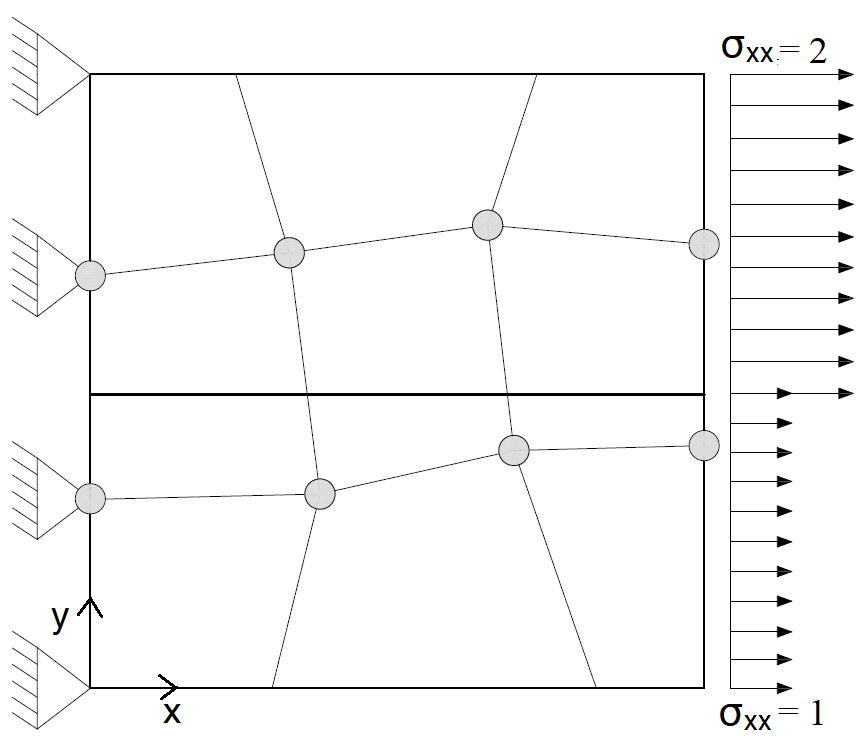}
  \caption{Geometry and loading conditions of the discontinuous patch test.
  }\label{fig:disc_PT}
\end{figure}
For the sake of simplicity, we assume $E = 1$ and $\nu = 0$, so that
the problem is reduced to one dimension.
As boundary conditions, we prescribe zero displacements along the edge
$x=0$, a discontinuous distribution of horizontal tractions along the
edge $x=1$ and zero tractions along the horizontal edges $y=0$ and
$y=1$:
\begin{subequations}\label{eq:disc:patch:DD}
  \begin{align*}
    &\bm{u}(0,y) = \bm{0} , \\
    &\sigma_{xx}(1,y) =
    \begin{cases} 1, & y \leq 1/2\\
      2, & y > 1/2
    \end{cases}, \quad \sigma_{yy}(1,y) = \sigma_{xy}(1,y) = 0, \\
    &\sigma_{yy}(x,0) = \sigma_{xy}(x,0) = 0, \\
    &\sigma_{yy}(x,1) = \sigma_{xy}(x,1) = 0.
  \end{align*}
\end{subequations}
For this problem, whose geometry and boundary conditions are depicted
in Fig.~\ref{fig:disc_PT}, the exact solution is the following
piecewise linear function
\begin{equation}
  \label{exact_sol}
  \bm{u}(x,y) =
  \begin{cases} 
    [x, 0]^T , &(x,y) \in \Omega^-\,,\\[0.5em]
    [2x, 0]^T, &(x,y) \in \Omega^+\,. 
  \end{cases}
\end{equation}
The exact solution \eqref{exact_sol} belongs to the discrete space.
In agreement with the expectations, the extended virtual element
formulation presented in Section~\ref{subsec:disc}, which uses
distinct projector operators on the two subdomains generated by the
horizontal crack, passes the proposed patch test with a relative error in strain energy of $2\times10^{-13}$.

\subsection{Convergence study}
\label{conv_study}
We study the convergence of the X-VEM for the problem of a
two-dimensional square plate under plain strain conditions that contains
a horizontal crack, extending from the boundary to the center of the
specimen.
The boundary conditions are such that mixed-mode conditions
prevail.
%% I and II crack opening takes place.
%%
The geometry of the domain is the same adopted as that for the extended patch
test in Section~\ref{ext_pt} and is shown in
Fig.~\ref{bench_geometry}.
On the boundary of the domain, we apply the exact near-tip
displacement fields \eqref{def:crackmodes}, which are also employed as
enrichment fields for the X-VEM and represent the exact solution for
the problem at hand.
%%
%%Such fields are singular on the crack tip, in the origin of the
%%coordinate system.
%%
%%A Young modulus $E$ equal to $10^5$ and a Poisson coefficient $\nu$
%%equal to $0.3$ are adopted.

\begin{figure}[!htb]
  \centering
  \begin{subfigure}{0.5\textwidth}
    \centering
    \includegraphics[width=\textwidth]{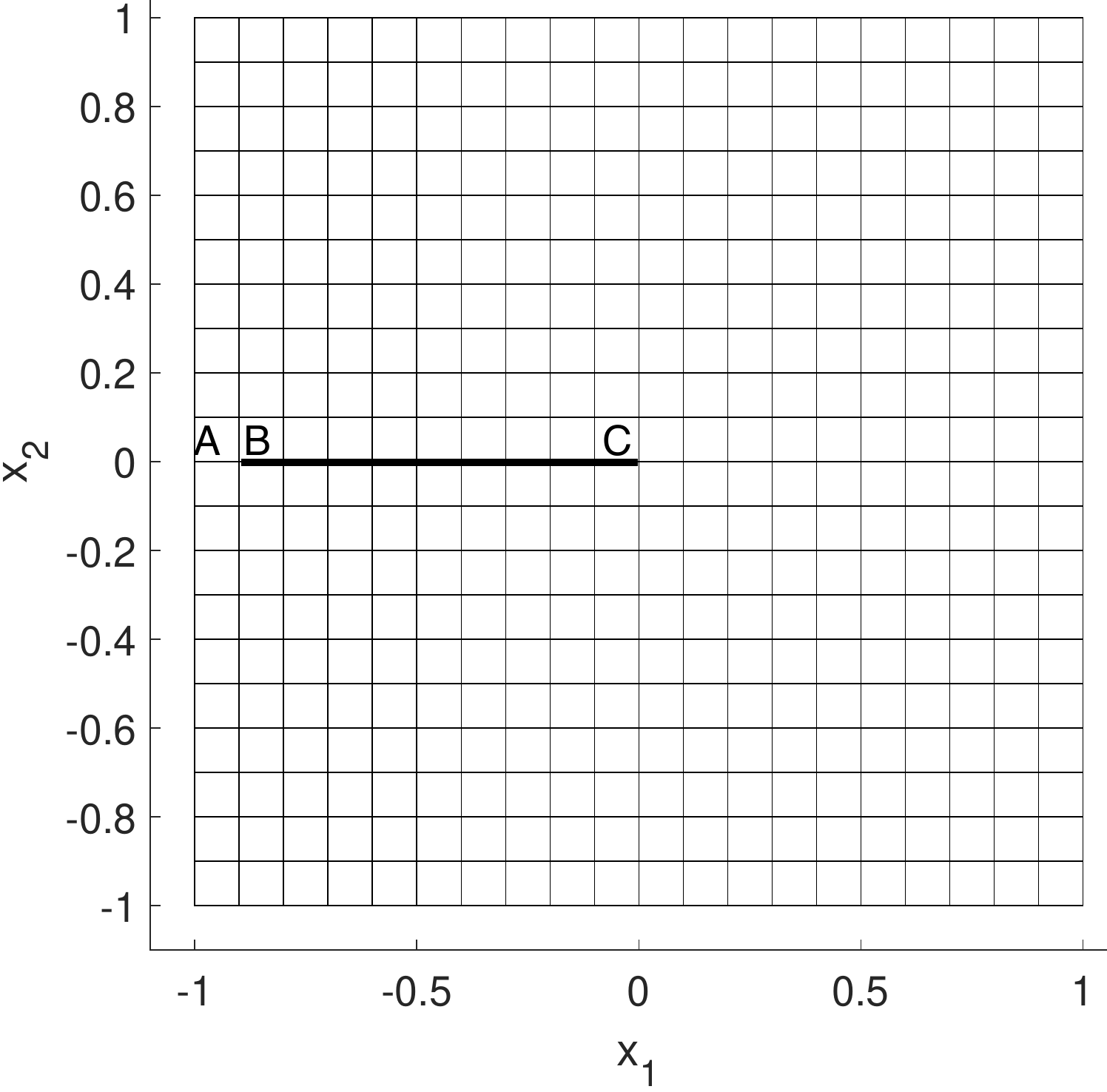}
    \caption{}\label{mesh_quad}
  \end{subfigure}\hfill
  \begin{subfigure}{0.5\textwidth}
    \centering
    \includegraphics[width=\textwidth]{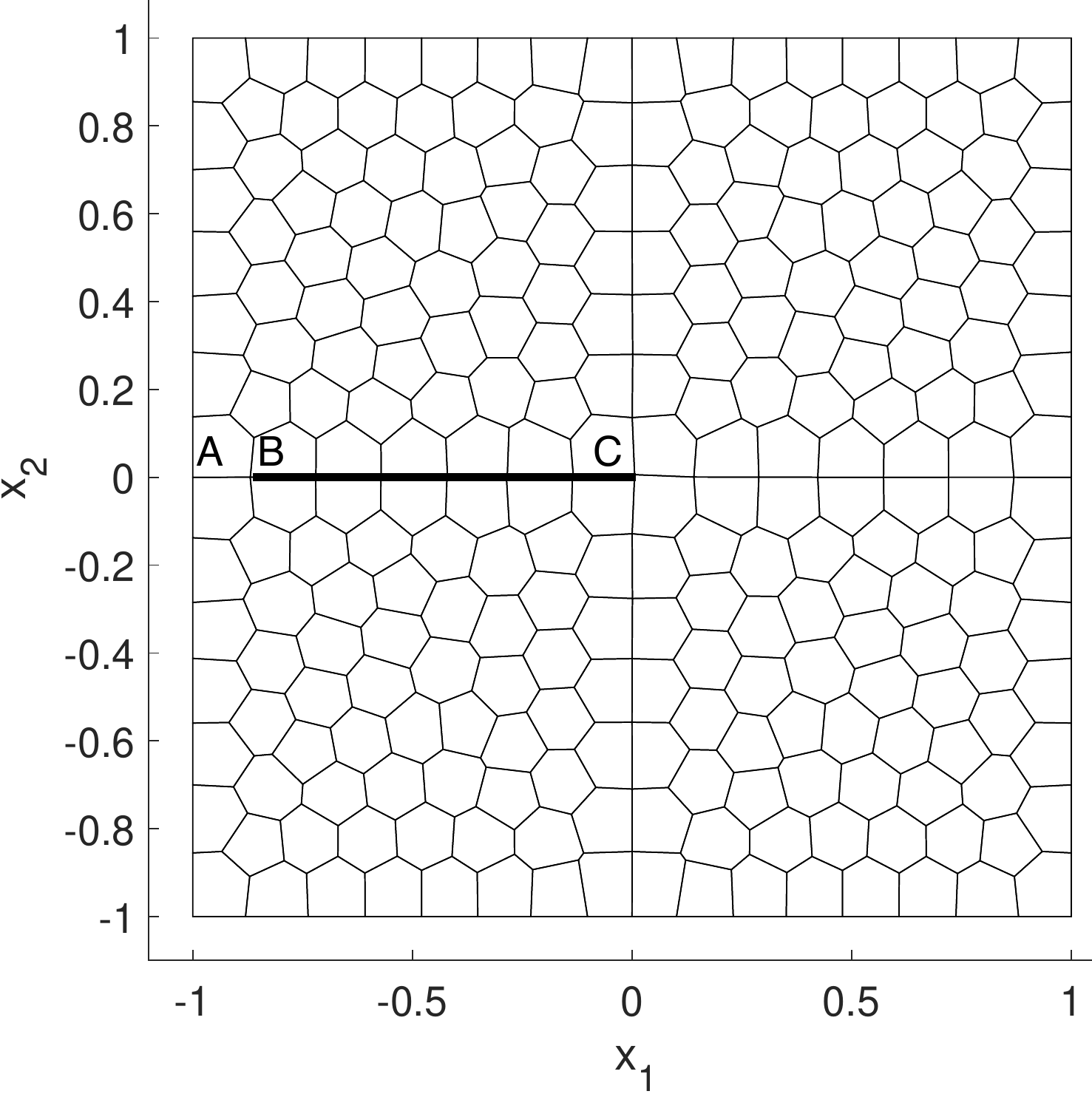}
    \caption{}\label{mesh_poly}
  \end{subfigure}
  \caption{Mixed-mode benchmark problem. (a) Quadrilateral mesh and
  (b) polygonal mesh.}\label{fig:meshes}
\end{figure}

In this study, we consider both quadrilateral and in general
polygonal meshes, see Fig.\ref{fig:meshes}.
Quadrilateral meshes are composed of $10 \times 10$, $20 \times 20$,
$40 \times 40$ and $80 \times 80$ square elements.
For the X-VEM, we use the stabilization in~\eqref{eq:stab2}, where
$\alpha = 1$ is chosen as the scaling parameter.
We generated the polygonal meshes from Voronoi tassellations by using
\texttt{Polymesher}~\cite{Paulino:2007:SMO}.
In order to apply essential boundary conditions, the crack is
explicitly meshed over the first element (AB), while the remaining
part of the crack (BC) is modeled by the X-VEM.

\begin{figure}[!htb]
  \centering
  \includegraphics[width=0.7\textwidth]{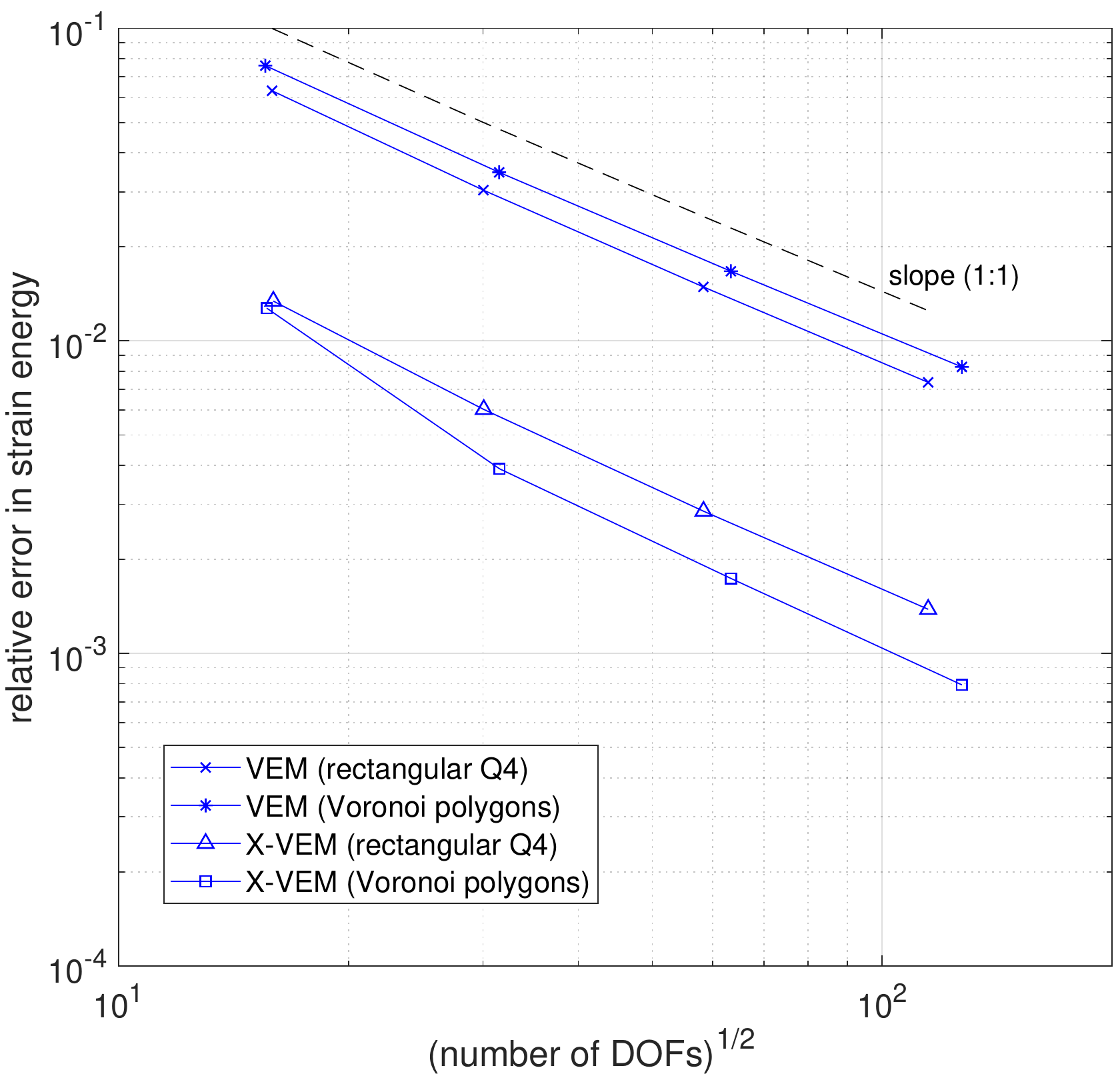}
  \caption{Convergence in strain energy for the mixed-mode benchmark problem. For the X-VEM, only the node at the origin is enriched
  (topological enrichment).
  Comparisons are shown with the standard VEM on quadrilateral and polygonal meshes. All methods converge with a rate close to
  unity.}\label{top_enr_conv}
\end{figure}

To compute the element stiffness matrix $\bm{K}_E$, we implement the
X-VEM of Section \ref{sec:implementation} following two different
strategies: \emph{topological} enrichment 
and \emph{geometric} enrichment.
\begin{figure}[!htb]
  \centering
  \begin{subfigure}{0.45\textwidth}
    \centering
    \includegraphics[width=\textwidth]{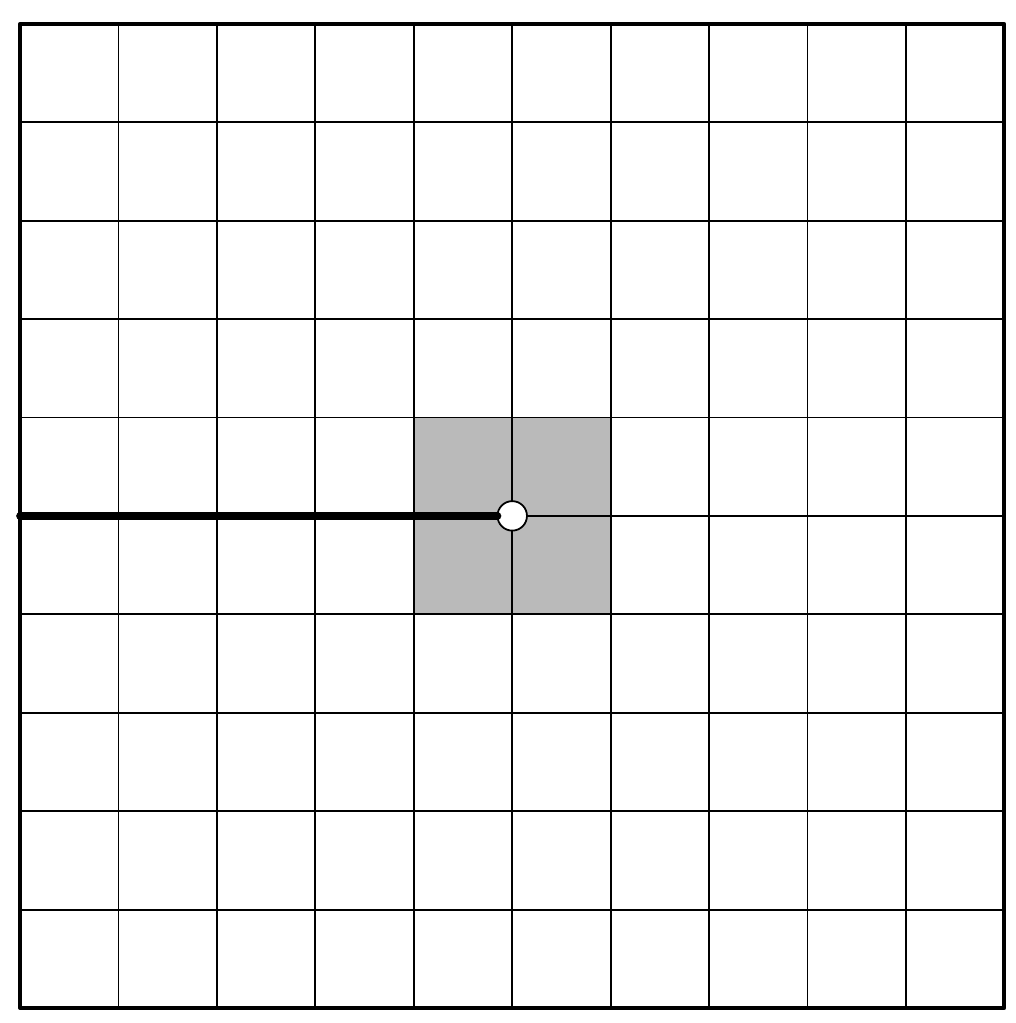}
    \caption{}\label{top_enr}
  \end{subfigure}\hfill
  \begin{subfigure}{0.45\textwidth}
    \centering
    \includegraphics[width=\textwidth]{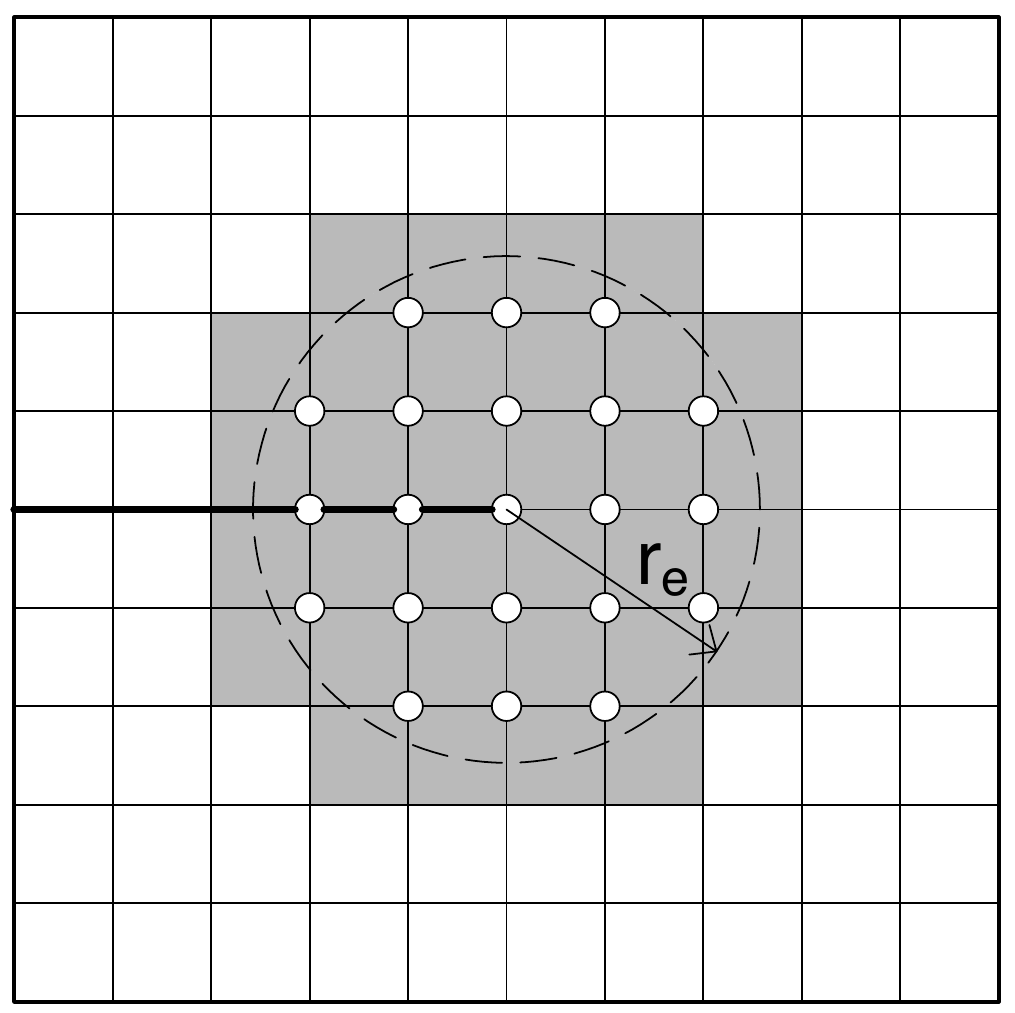}
    \caption{}\label{geom_enr}
  \end{subfigure}
  \caption{Mixed-mode benchmark problem. (a) Topological enrichment and (b) geometric enrichment.}\label{fig:enr_types}
\end{figure}
In the topological enrichment, graphically represented in
Fig.~\ref{top_enr}, we only enrich the node located at the singularity of the solution.
The convergence rate for this problem is given by
$R=\min(2\lambda,2p)$, where $\lambda$ is the order of the singularity
and $p$ the polynomial degree~\cite{Grisvard:1985}.
Since in our case $\lambda = 1/2$ and $p = 1$, we obtain a convergence rate $R=1$ that is non-optimal, as we anticipated in Section ~\ref{subsec:partial-enrichment}.
In fact, this suboptimal convergence rate is also noted in enriched finite
element techniques for fracture problems, cf.~\cite{laborde:2005}.
Figure~\ref{top_enr_conv} shows convergence plots of the relative error
in strain energy. The expected convergence rate $R$ is reported in the
graph. 
Both VEM and X-VEM with topological enrichment converge in strain
energy with a rate close to 1, in agreement with theory.
It turns out that the X-VEM is insensitive to the type of mesh
(quadrilaterals or polygons), and the results from the X-VEM are
consistently more accurate than those from standard VEM.

\begin{figure}[!htb]
  \centering
  \includegraphics[width=0.75\textwidth]{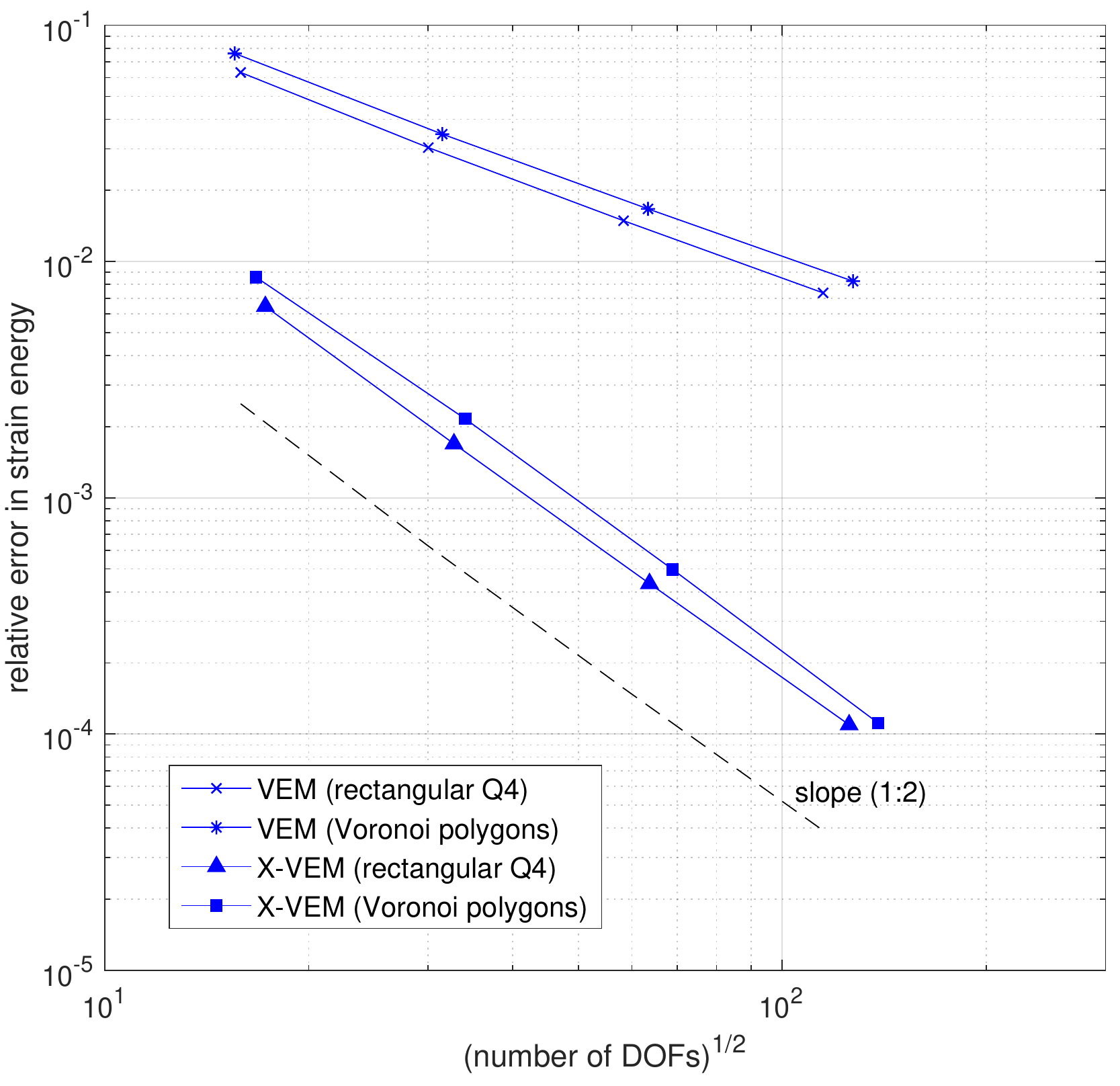}
  \caption{Convergence in strain energy for the mixed-mode benchmark  problem. For the X-VEM, geometric enrichment ($r_e = 0.5$) on quadrilateral
    and polygonal meshes is used.
    Comparisons are made with the standard VEM.  X-VEM converges with a
    rate close to two.}\label{geom_enr_conv}
\end{figure}

\begin{figure}[!htb]
  \centering
  \includegraphics[width=0.75\textwidth]{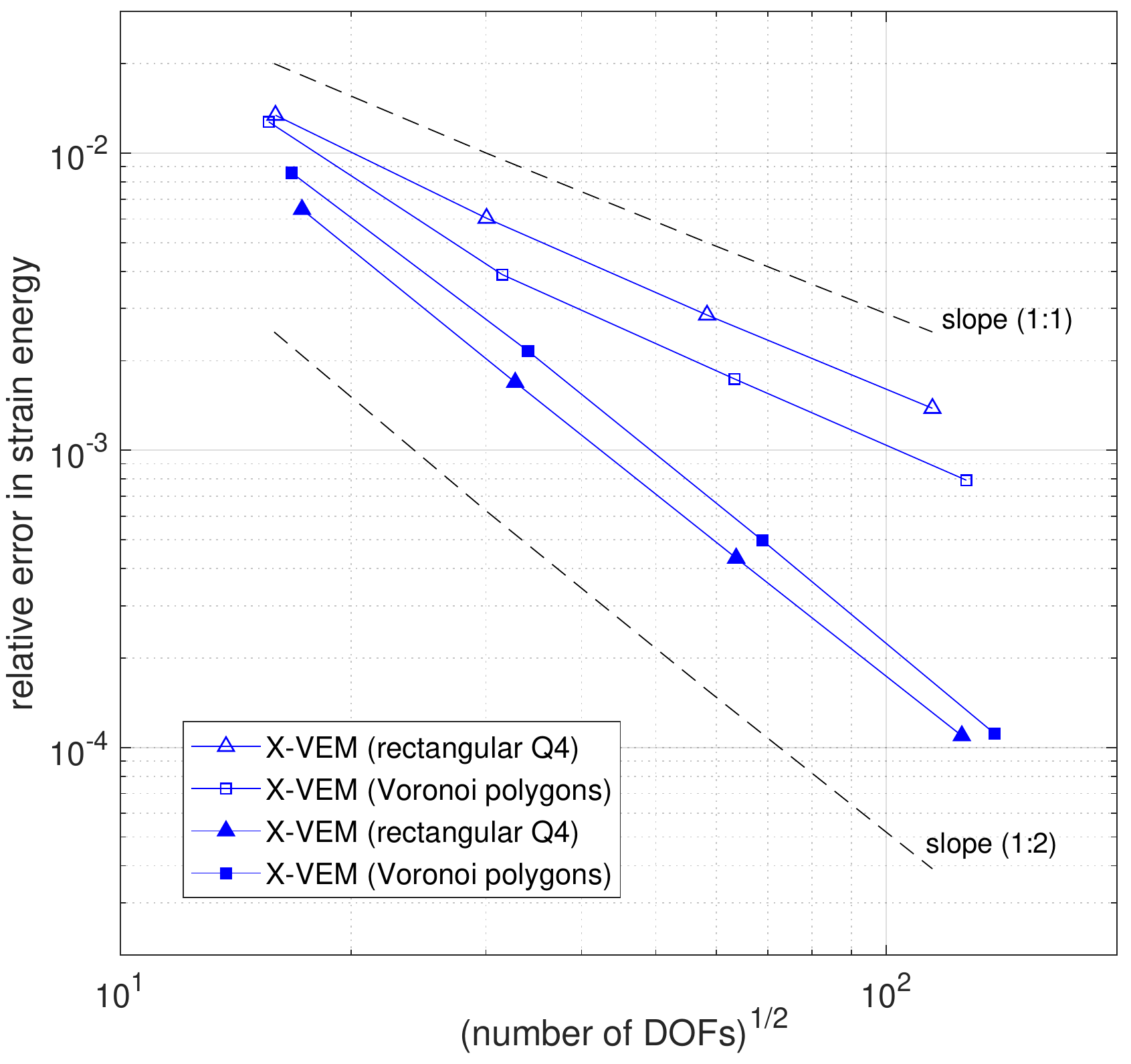}
  \caption{Convergence in strain energy for the mixed-mode benchmark  problem using topological enrichment (hollow markers) and
  geometric enrichment with $r_e = 0.5$ (solid markers) on
  quadrilateral and polygonal meshes.}\label{topo-geom-enr}
\end{figure}

Many prior studies have shown that geometric enrichment,
i.e., enriching all the nodes within a given radius from the
singularity at the crack tip, allows the standard X-FEM for fracture
problems to recover the optimal convergence rate~\cite{laborde:2005,
  Bechet:2005:IIR}.
In order to establish if the proposed X-VEM can deliver the
optimal convergence rate $R = 2$ that is predicted by theory, we enrich
all nodes that are located within a ball of radius $r_e = 0.5$ from the origin
(see Fig.~\ref{geom_enr}).
Figure~\ref{geom_enr_conv} depicts convergence plots for the relative
error in strain energy on quadrilateral and polygonal meshes for the
X-VEM with geometric enrichment.
The convergence rate is close to 2, which is consistent with  theory. 
To provide a clearer picture, Fig.~\ref{topo-geom-enr} shows a
comparison between the convergence plots in strain energy for both
quadrilateral and polygonal meshes.

In order to assess the robustness and the accuracy of the X-VEM in
providing stress intensity factors, we apply the procedure described in Section~\ref{SIFS} to the problem at hand.
For this example, the exact mode I and mode II stress intensity
factors are $K_I = 1$ and $K_{II} = 1$. Both topological and geometric
enrichment are considered.
Stress intensity factors are computed considering a ring of elements
placed at a radius $r_d = 0.4$ from the origin.
\begin{figure}[!p]
  \centering
  \begin{subfigure}[b]{\textwidth}
  \centering
  \includegraphics[width=0.65\textwidth]{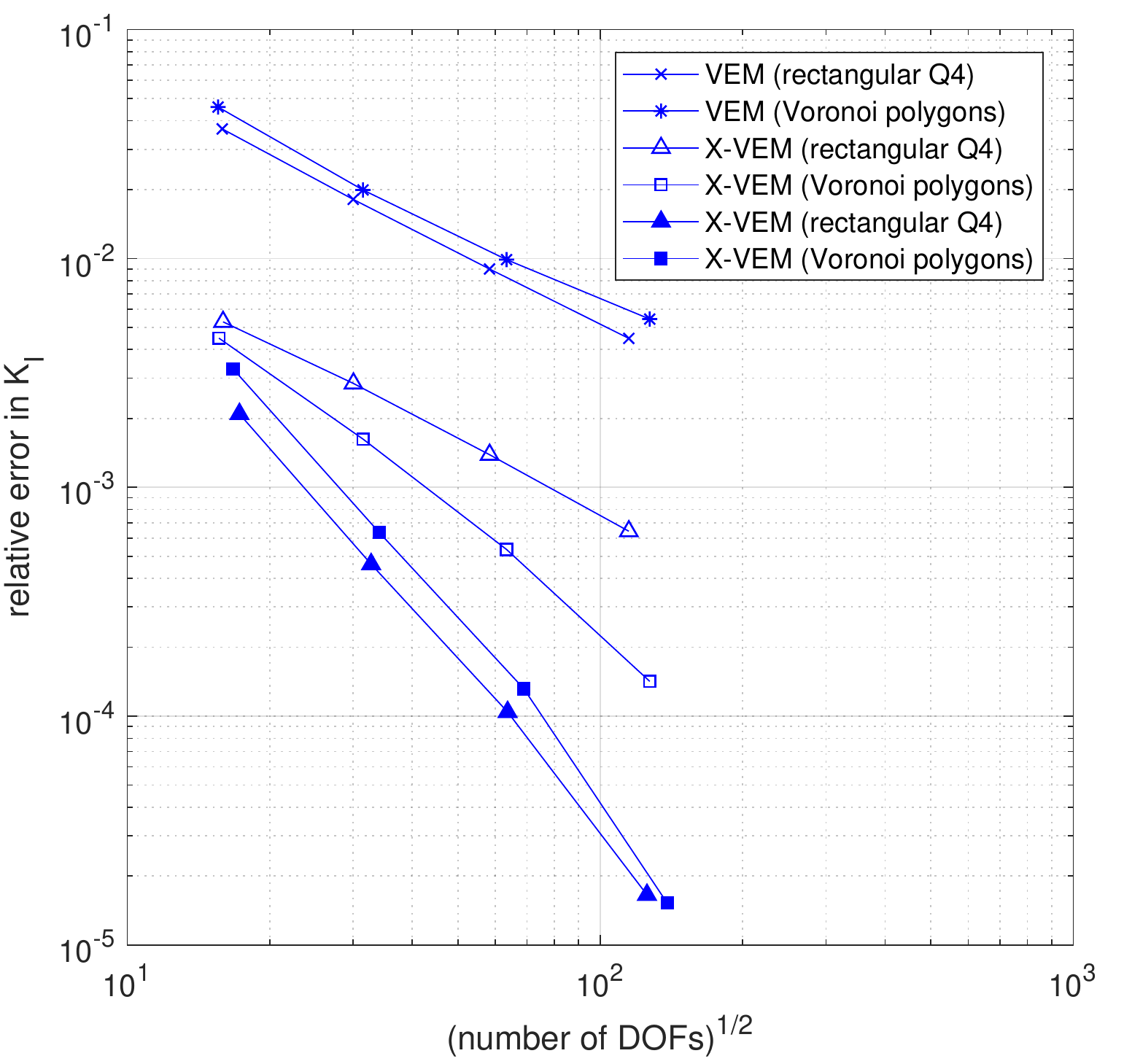}
  \caption{} \label{K1}
  \end{subfigure} 
  \begin{subfigure}[b]{\textwidth}
  \centering
  \includegraphics[width=0.65\textwidth]{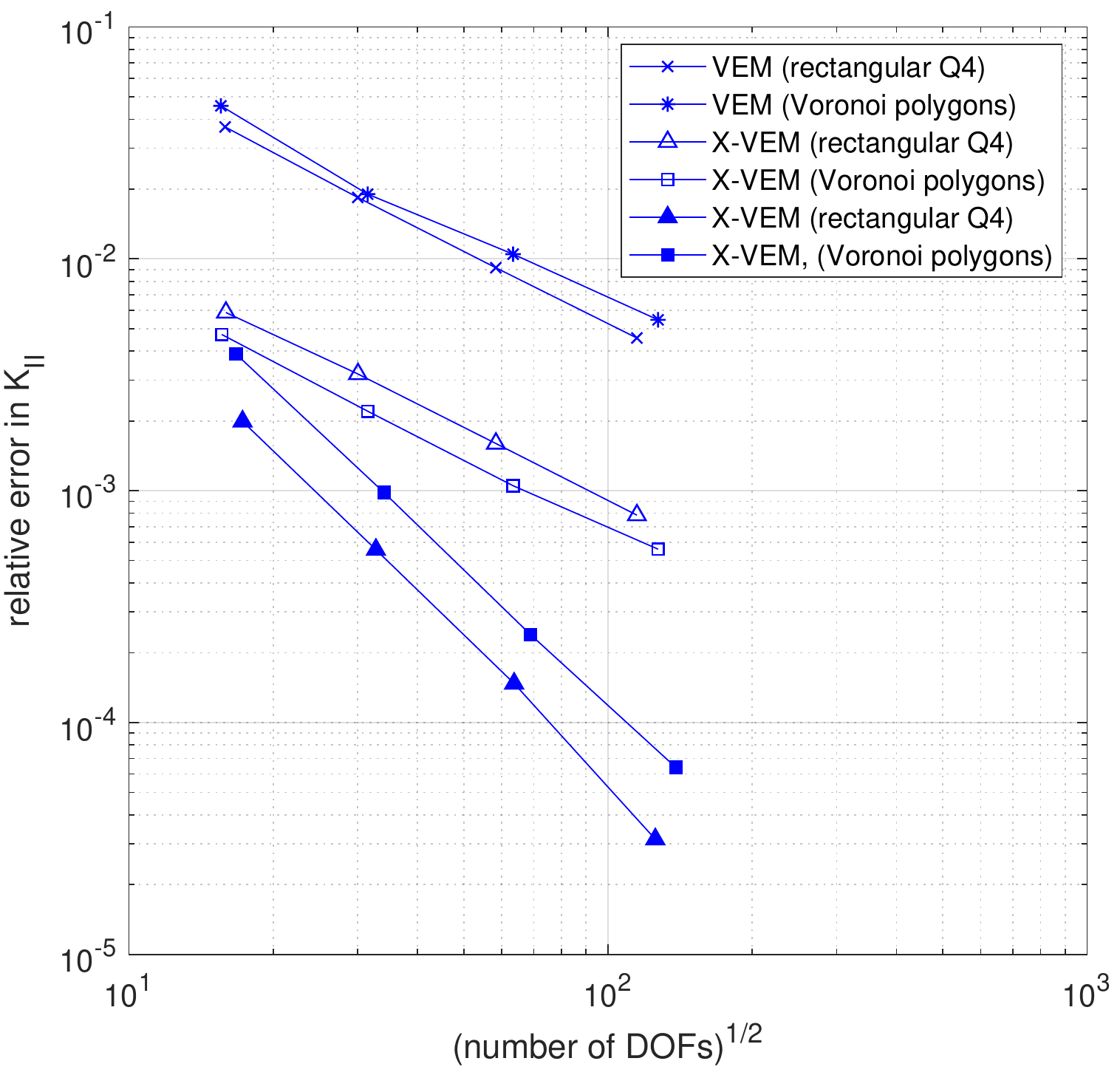}
  \caption{} \label{K2}
  \end{subfigure}
  \caption{Convergence of stress intensity factors for the mixed-mode benchmark problem 
  using
    topological enrichment (hollow markers) and geometric enrichment
    with $r_e = 0.5$ (solid markers) on quadrilateral and
    polygonal meshes. (a) $K_I$ and (b) $K_{II}$.}\label{K1K2}
\end{figure} 

\begin{figure}[!htb]
  \centering
  \includegraphics[width=0.7\textwidth]{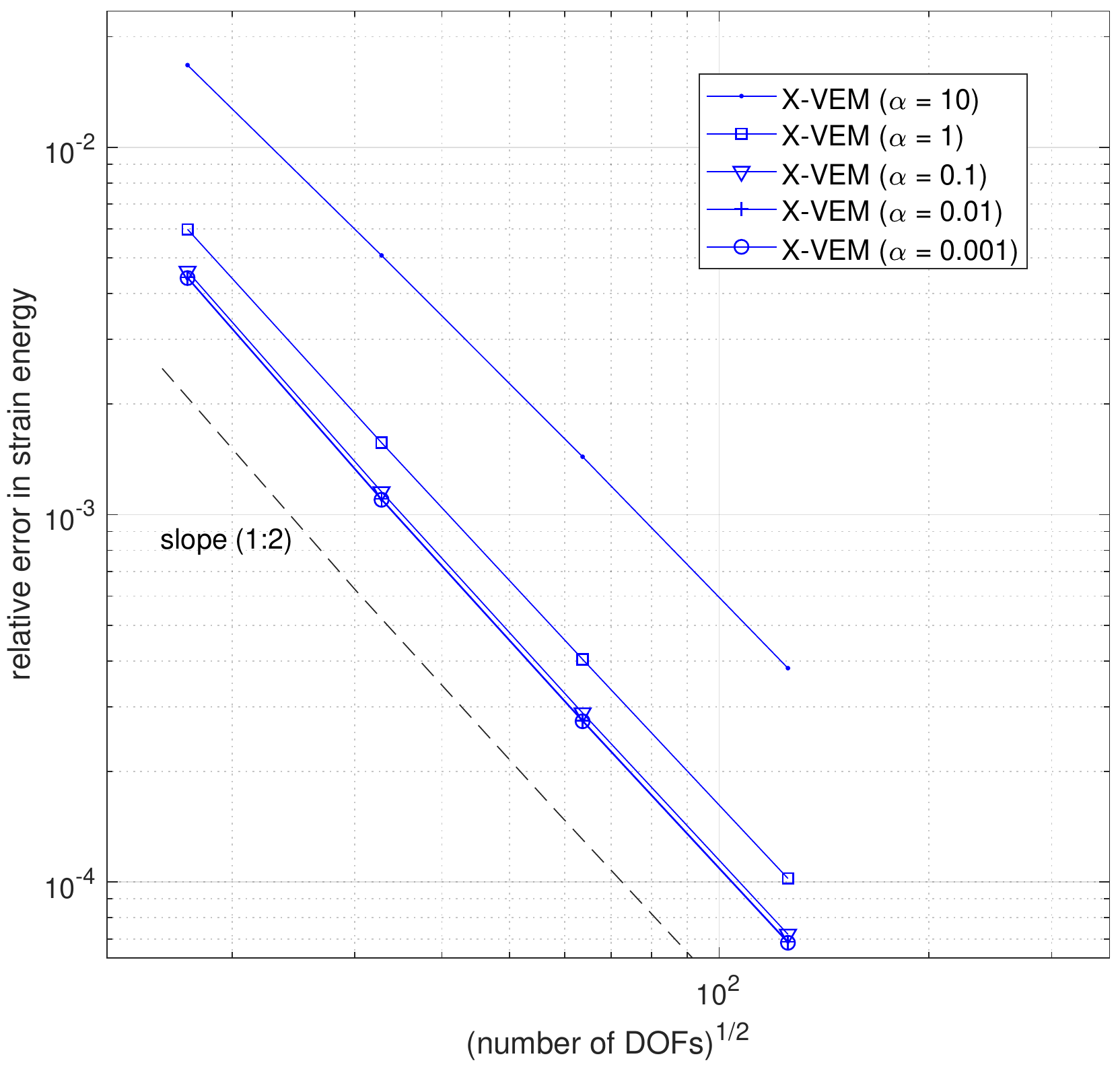}
  \caption{Convergence in strain energy for the mixed-mode benchmark problem. Convergence curves are shown for the X-VEM on a quadrilateral mesh with geometric enrichment and different choices of the stabilization parameter $\alpha$.}\label{alpha}
\end{figure}

\begin{figure}[!p]
  \centering
  \begin{subfigure}[b]{\textwidth}
  \centering
  \includegraphics[width=0.44\textwidth]{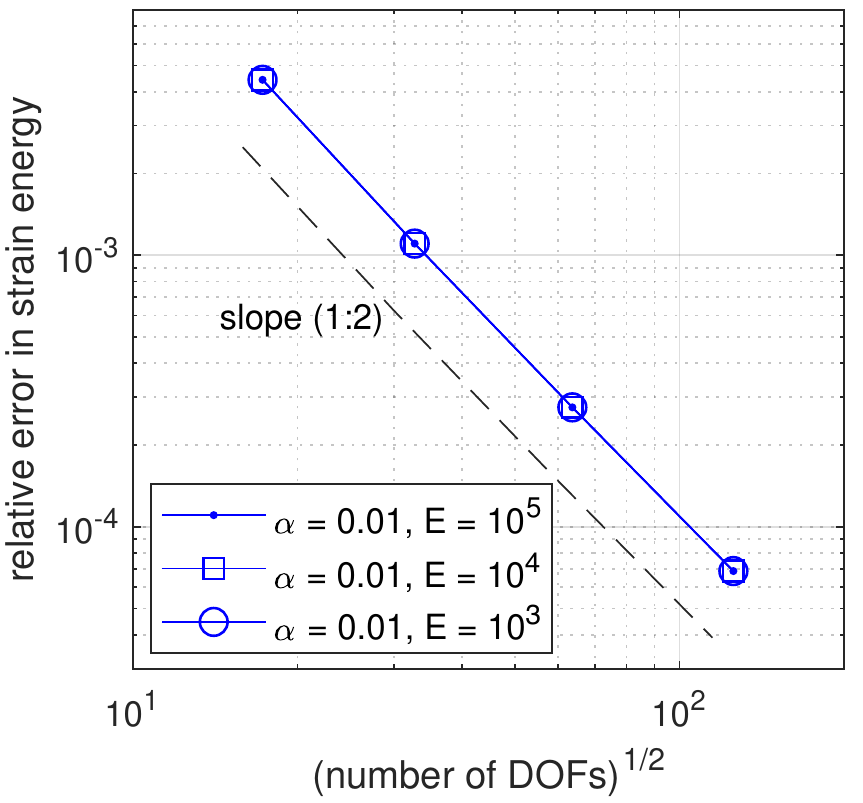}
  \caption{} \label{E-alpha-a}
  \end{subfigure} 
  \begin{subfigure}[b]{\textwidth}
  \centering
  \includegraphics[width=0.44\textwidth]{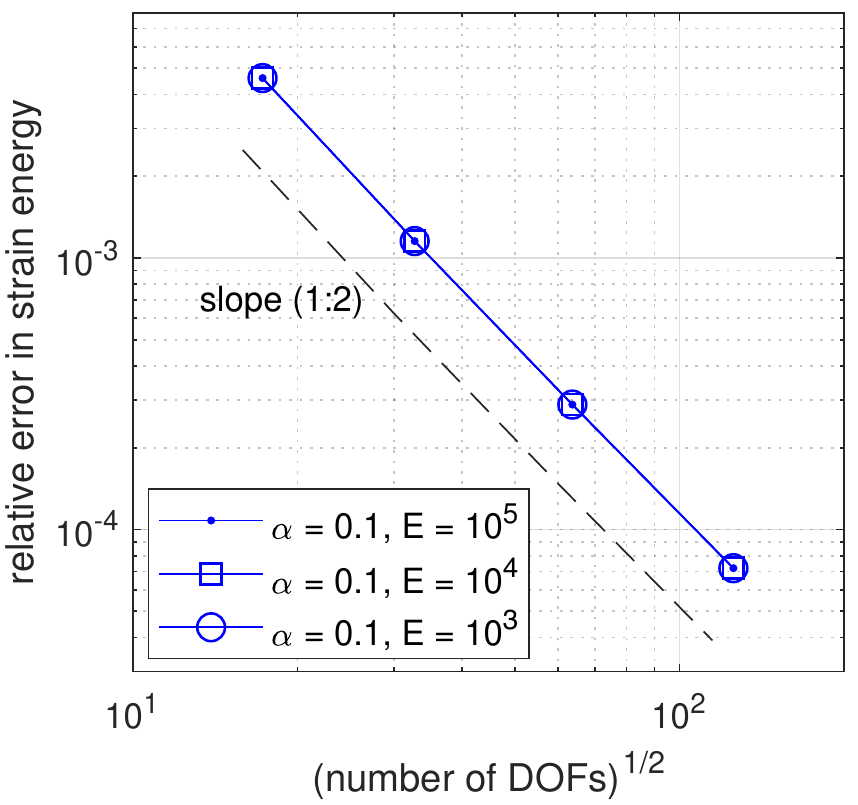}
  \caption{} \label{E-alpha-b}
  \end{subfigure}
   \begin{subfigure}[b]{\textwidth}
  \centering
  \includegraphics[width=0.44\textwidth]{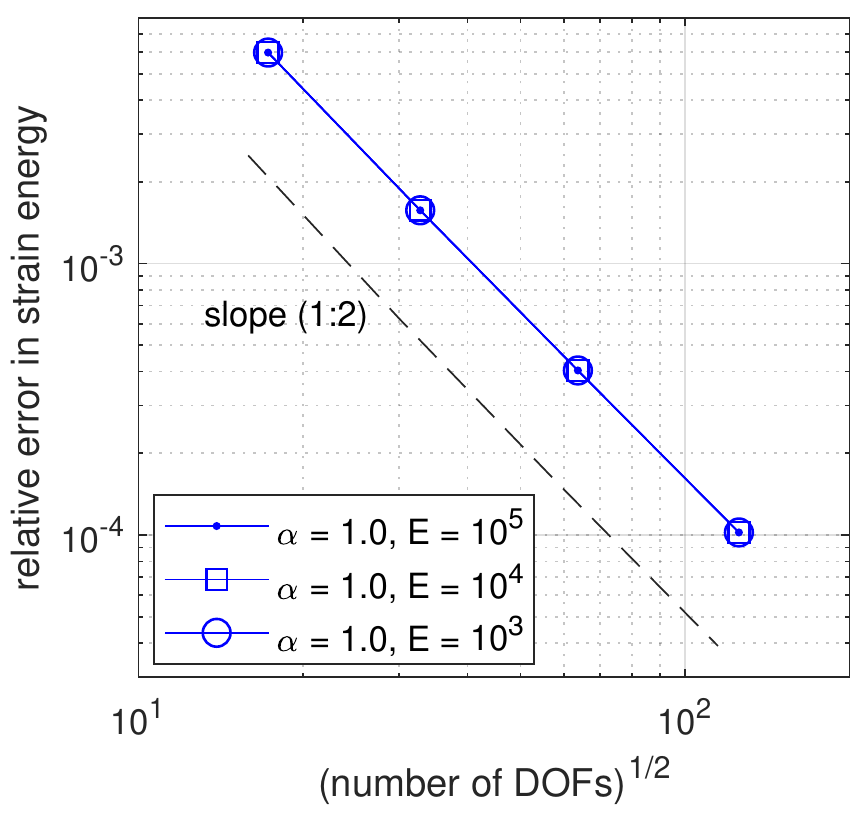}
  \caption{} \label{E-alpha-c}
  \end{subfigure} 
  \caption{Convergence in strain energy for the mixed-mode benchmark problem. Convergence curves are shown for the X-VEM on a quadrilateral mesh with geometric enrichment and different choices of both the stabilization parameter $\alpha$ and the Young modulus $E$: (a) $\alpha = 0.01$, (b) $\alpha = 0.1$, (c) $\alpha = 1.0$}\label{E-alpha}
\end{figure}

Figures~\ref{K1} and~\ref{K2} show the convergence of $K_I$ and
$K_{II}$ on
quadrilateral and polygonal meshes.
Convergence is stable on all the meshes and accuracy is sound. In
particular, geometric enrichment enhances both the convergence rate and
the accuracy.
Finally, in Fig.~\ref{alpha} we investigate the influence of the scaling parameter $\alpha$ in the stabilization. Convergence is optimal for 
$\alpha$ ranging from $10^{-3}$ to $10$, with trends in improved accuracy
towards smaller values of $\alpha$. Moreover, as shown in Fig.~\ref{E-alpha}, $\alpha$ does not need to be adjusted if the Young modulus $E$ is varied: for a given $\alpha$, the accuracy of the method is not significantly influenced by varying $E$.

\subsection{Inclined edge crack in a finite plate under uniform tension}
\label{inclined}
\begin{figure}[!htb]
  \centering
  \begin{subfigure}[b]{0.49\textwidth}
    \centering
    \includegraphics[width=0.7\textwidth]{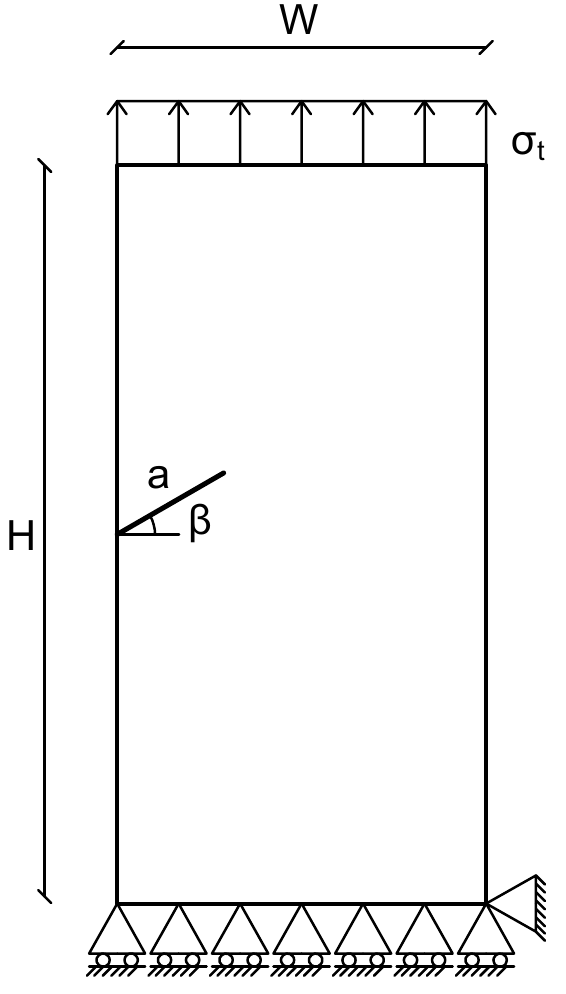}
    \caption{}\label{geom_incl}
  \end{subfigure} \hfill
  \begin{subfigure}[b]{0.49\textwidth}
    \centering
    \includegraphics[width=0.6\textwidth]{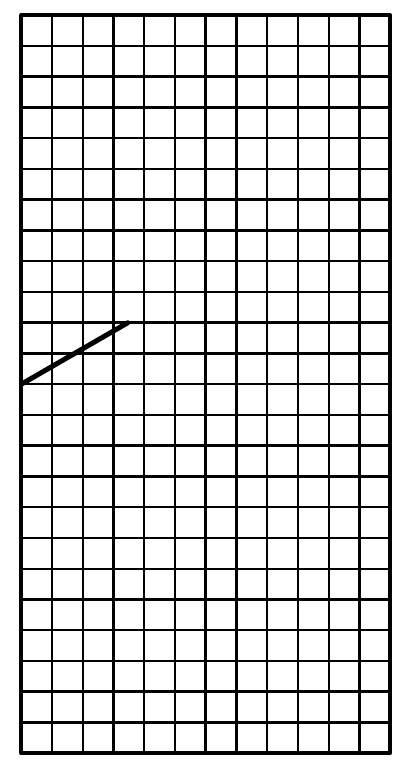}
    \vspace{3mm}
    \caption{}\label{mesh_incl}
  \end{subfigure}
  \caption{Inclined edge crack in a finite plate under uniform tension. (a) Problem geometry and (b) virtual element mesh.}
  \label{fig:inclined}
\end{figure}

We now study the problem of an inclinded edge crack in 
a finite plate under uniform tension. The geometry and boundary conditions are shown in Fig.~\ref{fig:inclined}. The plate width $W = 6$ and  plate height $H = 3$ are chosen. The crack has length $a = 1$ and is inclined at an angle $\beta$. Uniform tractions $\sigma_t = 1$ are applied on the top edge and horizontal rollers are imposed on the bottom edge. % A Young modulus $E = 10^5$ and a Poisson coefficient $\nu = 0.3$ are adopted. 
The exact solution for this problem in the neighborhood of the crack tip is given by a linear combination of the fields~\eqref{def:crackmodes}. However, the exact solution on the whole domain is not known in 
closed-form. 

\begin{figure}[!htb]
  %\centering
  \vspace{-1cm}
  \begin{subfigure}[b]{\textwidth}
  \centering
  \includegraphics[width=0.7\textwidth]{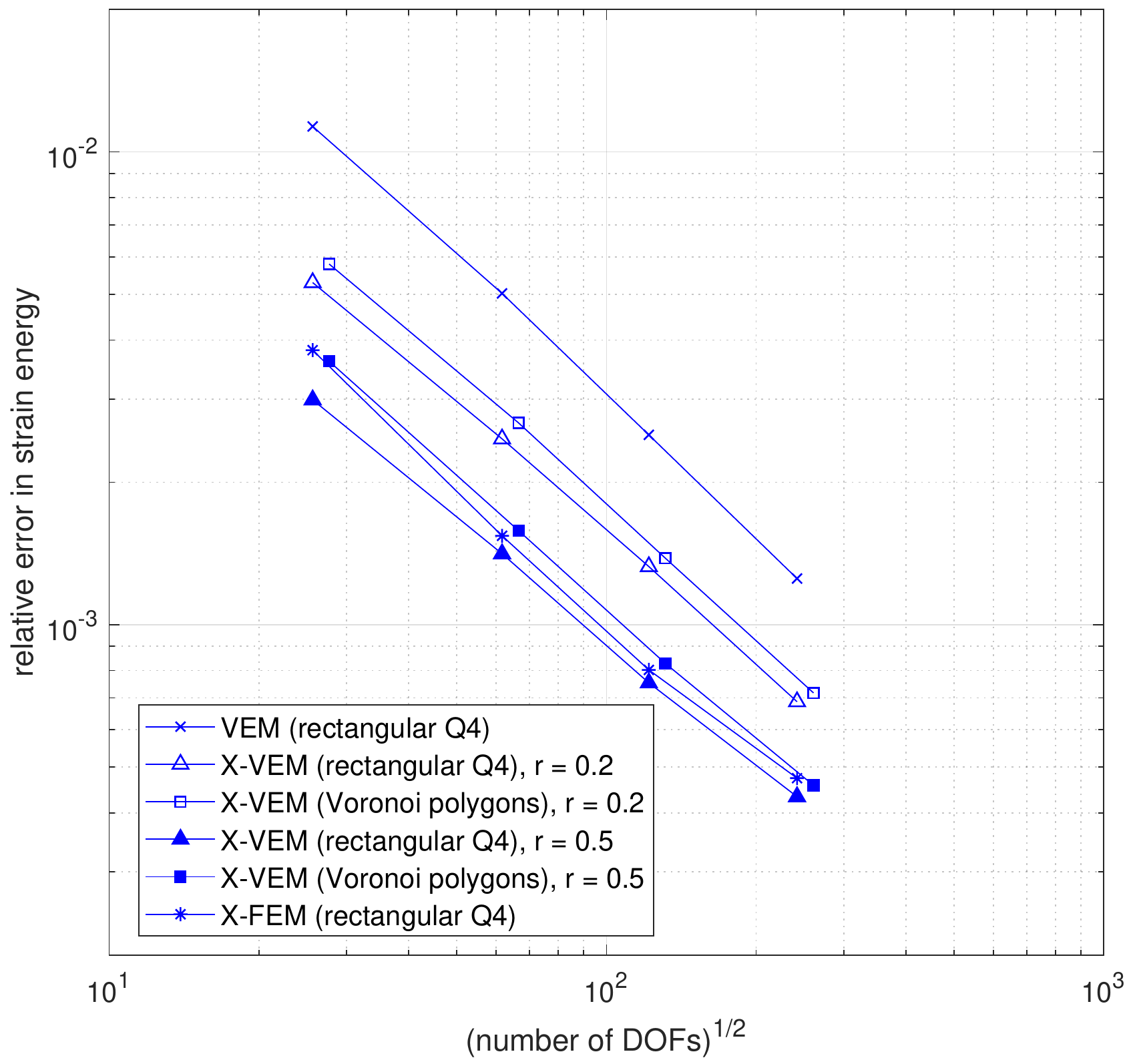}
  \caption{}\label{conv_inclined_beta0}
  \end{subfigure}
  \\
  \begin{subfigure}[b]{\textwidth}
  \centering
  \includegraphics[width=0.7\textwidth]{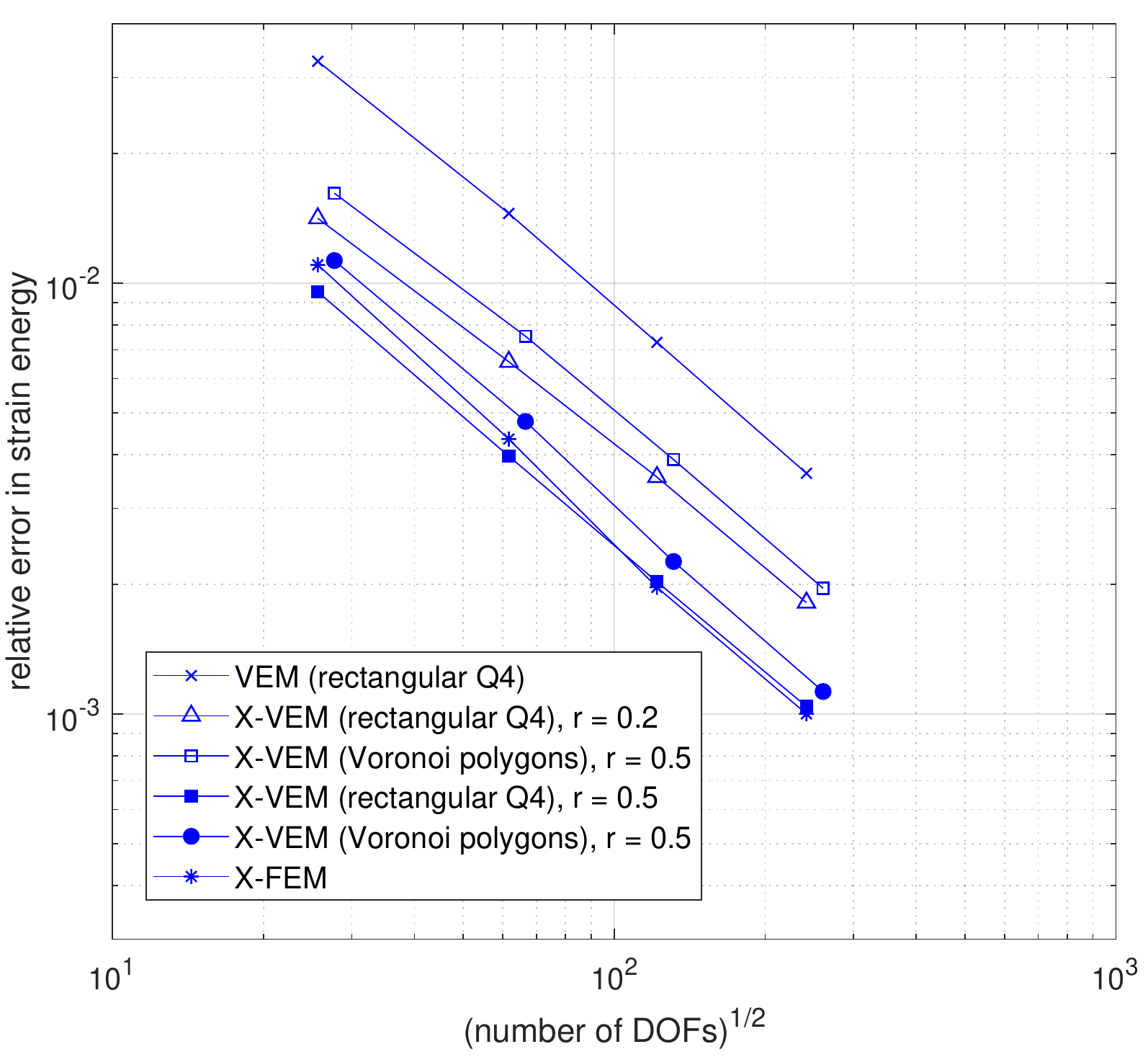}
  \caption{}\label{conv_inclined}
  \end{subfigure}
  \caption{Convergence in strain energy for the problem of an inclined edge crack in a finite plate under uniform tension. Convergence curves are shown for
  quadrilateral and polygonal meshes, with varying enrichment radii.
  (a) $\beta = 0$ and (b) $\beta = \pi / 6 $.
  }\label{conv_inclined_beta0_and_pi6}
\end{figure} 

%\begin{figure}[!htb]
%  \centering
%  \includegraphics[width=0.75\textwidth]{conv_inclined}
%  \caption{Convergence in strain energy for the problem of an inclined % edge crack in a finite plate under uniform tension ($\beta = \pi/6$). % Convergence curves are shown for
%  quadrilateral and polygonal meshes, with varying enrichment radii.}\label{conv_inclined}
% \end{figure}  

\begin{figure}[!htb]
  \centering
  \includegraphics[width=0.7\textwidth]{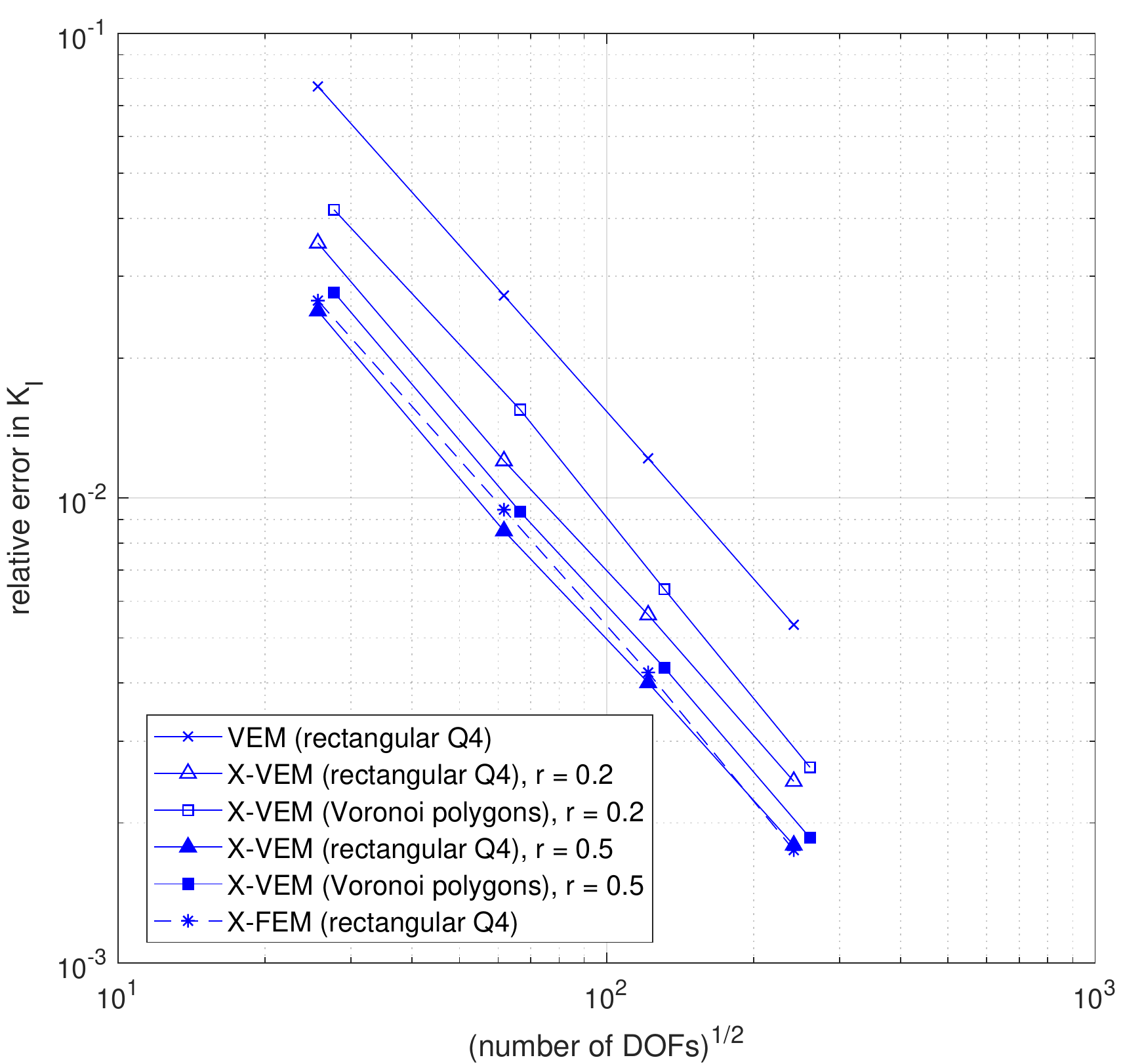}
  \caption{Convergence of $K_I$ for the problem of an inclined edge crack in a finite plate under uniform tension ($\beta = \pi/6$). Convergence curves are shown for 
  quadrilateral and polygonal meshes, with varying enrichment radii.}\label{KI_conv_inclined}
\end{figure}

Figures~\ref{conv_inclined_beta0} and~\ref{conv_inclined} show the convergence plots 
% the problem at hand for 
%% two different values of inclination of the crack, 
(inclination angles, $\beta = 0$ and $\beta = \pi/6$) in terms of the
relative error in strain energy on quadrilateral and polygonal meshes.
The stabilization parameter $\alpha = 0.01$ is chosen. Reference solution for the energy is computed with an overkill mesh of 460.800 elements using the X-FEM. We use meshes of square elements with 
$ h = 1/4$, $1/10$, $1/20$, $1/40$, as well as polygonal (Voronoi) meshes. %  polygonals elements obtained by Voronoi tassellation.  
Convergence of the X-VEM is compared to standard VEM and X-FEM. The X-VEM displays sound accuracy for both mesh types, and is comparable to that obtained with the X-FEM. 
Finally, in Fig.~\ref{KI_conv_inclined}, convergence of $K_I$
for $\beta = \pi/6$ is presented for different types of meshes and enrichment radii. Again, the obtained results are in good agreement with the X-FEM, which can be inferred from Table~\ref{SIF_table}, where the numerical results for $K_I$ and $K_{II}$ are listed for the X-FEM and for the X-VEM with different values of the stabilization parameter $\alpha$. As already noted in Section~\ref{conv_study}, here too smaller values of $\alpha$ improve the accuracy of the X-VEM. 

\begin{table}[!p]
\centering
\caption{Stress intensity factors, $K_I$ and $K_{II}$, for a finite plate with an inclined edge crack under uniform tension. Numerical results are listed for the X-FEM and the X-VEM (different values of the stabilization parameter $\alpha$) on a $60 \times 120$ mesh of square elements.} \label{SIF_table}
  \begin{tabular}{c c c c c| c c c c}
    \hline
    \multicolumn{1}{c}{} &
    \multicolumn{4}{c|}{$K_I$} &
    \multicolumn{4}{c}{$K_{II}$} \\
    \cline{2-9}
    \multirow{1}{*}{$\beta$} &
      \multicolumn{3}{c}{X-VEM} &
      \multicolumn{1}{c|}{X-FEM} &
      \multicolumn{3}{c}{X-VEM} &
      \multicolumn{1}{c}{X-FEM} \\
     % \cmidrule(lr){2-4}
     %  \cmidrule(lr){2-4}
    & $\alpha=0.01$ & $\alpha=0.05$ & $\alpha=0.10$ &  & $\alpha=0.01$ & $\alpha=0.05$ & $\alpha=0.10$ & \\
    \hline
    $\pi/12$ & 2.9351 & 2.9333 & 2.9262 & 2.9349 & 0.4631 & 0.4627 & 0.4615 & 0.4630 \\
    \hline
    $\pi/6$ & 2.3652 & 2.3639 & 2.3582 & 2.3651 & 0.7607 & 0.7603 & 0.7584 & 0.7606 \\
    \hline
    $\pi/4$ & 1.6418 & 1.6408 & 1.6370 & 1.6419 & 0.8333 & 0.8329 & 0.8308 & 0.8334 \\
    \hline
  \end{tabular}
    
\end{table}

%% file: sec6_final.tex
% Hey Emacs, this is -*-latex-*-

\section{Concluding remarks}
\label{conclusions}

We developed a stable and convergent extended virtual element method
for two-dimensional elastic fracture problems, which permits the incorporation of 
crack-tip singularities and discontinuities in the
approximation space.
Inspired by the construction of the X-FEM~\cite{Moes:1999:FEM}, we
augmented the standard virtual element space by means of additional
vectorial basis functions that were constructed using the
asymptotic mode I and mode II crack-tip displacement fields as
enrichment functions.
An extended elliptic projector was proposed that projects the
functions of the extended virtual element space onto the space spanned
by linear polynomials and the enrichment fields.
Crack discontinuities were modeled by decomposing each virtual shape
function as the sum of two discontinuous shape functions, 
following
the approach proposed by Hansbo and Hansbo~\cite{Hansbo:2004:FEM}.
The proposed extended virtual element formulation does not present integration issues,
since all integrals are computed on the elements boundary, where
virtual shape functions are known.
A one-dimensional Gauss quadrature rule proved to be sufficient.
On forming the element projection matrix, the consistency part of
the stiffness matrix was computed and standard VEM procedures were followed
to obtain the element stabilization matrix.
Special attention was required to form the
stabilization matrix on partially enriched elements, and 
an ad-hoc stabilization strategy was devised that delivered accurate results.
Finally, we proposed a procedure for the computation of stress
intensity factors, which entails the discretization of the annular
J-domain with a ring of polygonal elements and the evaluation of the
interaction integral after transforming it to a boundary integral by
means of the divergence theorem.

In order to assess the consistency and the robustness of the proposed
X-VEM, we conducted several numerical tests.
First, we carried out an extended patch test to ensure the consistency
of the method with the mode I and mode II near-tip displacement
fields chosen as enrichments.
We also performed a discontinuous patch test to verify the consistency
of the Hansbo and Hansbo approach~\cite{Hansbo:2004:FEM} that we used
to incorporate the discontinuities.
Then, we performed convergence studies on quadrilateral and polygonal meshes with the X-VEM on the benchmark
problem of an edge crack in a square plate that is subjected to boundary conditions that are consistent with the exact mixed-mode near-tip displacement solutions.
To this end, we considered topological and geometric enrichment strategies.
In particular, we showed that geometric enrichment allows the
method to deliver optimal convergence rates in strain energy as well as in terms of mixed-mode
stress intensity factors. Finally, the X-VEM was used to solve the problem of an inclined crack in a finite plate under uniform tension, and the SIF results from the X-VEM were found to be in good agreement with those computed using standard X-FEM.
In future work, we will investigate the proposed extended virtual element formulation for tracking crack growth in both two-
and three-dimensional elastic media. 

%% file: main_XVEM_elasticity.bbl
\begin{thebibliography}{10}

\bibitem{Melenk:1996:PUF}
J.~M. Melenk and I.~Babu\v{s}ka.
\newblock The partition of unity finite element method: {Basic} theory and
  applications.
\newblock {\em Computer Methods in Applied Mechanics and Engineering},
  139:289--314, 1996.

\bibitem{Babuska:1997:PUM}
I.~Babu\v{s}ka and J.~M. Melenk.
\newblock The partition of unity method.
\newblock {\em International Journal for Numerical Methods in Engineering},
  40:727--758, 1997.

\bibitem{Moes:1999:FEM}
N.~Mo{\"e}s, J.~Dolbow, and T.~Belytschko.
\newblock A finite element method for crack growth without remeshing.
\newblock {\em International Journal for Numerical Methods in Engineering},
  46(1):131--150, 1999.

\bibitem{Tabarraei:2008:EFE}
A.~Tabarraei and N.~Sukumar.
\newblock Extended finite element method on polygonal and quadtree meshes.
\newblock {\em Computer Methods in Applied Mechanics and Engineering},
  197(5):425--438, 2008.

\bibitem{Zamani:2011:EIP}
A.~Zamani and M.~R. Eslami.
\newblock Embedded interfaces by polytope {FEM}.
\newblock {\em International Journal for Numerical Methods in Engineering},
  88:715--748, 2011.

\bibitem{Song:1997:CMAME}
C.~Song and J.~P. Wolf.
\newblock The scaled boundary finite-element method – alias consistent
  infinitesimal finite-element cell method – for elastodynamics.
\newblock {\em Computer Methods in Applied Mechanics and Engineering},
  147:329--355, 1997.

\bibitem{Song:2002:CS}
C.~Song and J.~P. Wolf.
\newblock Semi-analytical representation of stress singularity as occurring in
  cracks in anisotropic multi-materials with the scaled boundary finite-element
  method.
\newblock {\em Computer and Structures}, 80:183--197, 2002.

\bibitem{Song:2018:RSB}
C.~Song, E.~T. Ooi, and S.~Natarajan.
\newblock A review of the scaled boundary finite element method for
  two-dimensional linear elastic fracture mechanics.
\newblock {\em Engineering Fracture Mechanics}, 187:45--73, 2018.

\bibitem{Mousavi:2010:DT}
S.~E. Mousavi and N.~Sukumar.
\newblock Generalized {Duffy} transformation for integrating vertex
  singularities.
\newblock {\em Computational Mechanics}, 45(2--3):127--140, 2010.

\bibitem{Chin:2015:NIH}
E.~B. Chin, J.~B. Lasserre, and N.~Sukumar.
\newblock Numerical integration of homogeneous functions on convex and
  nonconvex polygons and polyhedra.
\newblock {\em Computational Mechanics}, 56(6):967--981, 2015.

\bibitem{chin:2017}
E.~B. Chin, J.~B. Lasserre, and N.~Sukumar.
\newblock Modeling crack discontinuities without element-partitioning in the
  extended finite element method.
\newblock {\em International Journal for Numerical Methods in Engineering},
  86(11):1021--1048, 2017.

\bibitem{BeiraodaVeiga-Brezzi-Cangiani-Manzini-Marini-Russo:2013}
L.~{Beir\~ao~da~Veiga}, F.~Brezzi, A.~Cangiani, G.~Manzini, L.~D. Marini, and
  A.~Russo.
\newblock Basic principles of virtual element methods.
\newblock {\em Mathematical Models \& Methods in Applied Sciences},
  23:119--214, 2013.

\bibitem{Lipnikov:2014}
K.~Lipnikov, G.~Manzini, and M.~Shashkov.
\newblock Mimetic finite difference method.
\newblock {\em Journal of Computational Physics}, 257(Part B.):1163--1227,
  2014.

\bibitem{BeiraodaVeiga-Lipnikov-Manzini:2014}
L.~Beir\~ao~da Veiga, K.~Lipnikov, and G.~Manzini.
\newblock {\em The Mimetic Finite Difference Method}, volume~11 of {\em MS\&A.
  Modeling, Simulations and Applications}.
\newblock Springer, {I} edition, 2014.

\bibitem{Beirao-Brezzi-Marini:2013}
L.~{Beir\~{a}o~da~Veiga}, F.~Brezzi, and D.~Marini.
\newblock Virtual elements for linear elasticity problems.
\newblock {\em SIAM Journal on Numerical Analysis}, 51(2):794--812, 2013.

\bibitem{Paulino:2014:CMAME}
A.~L. Gain, C.~Talischi, and G.~H. Paulino.
\newblock On the virtual element method for three-dimensional linear elasticity
  problems on arbitrary polyhedral meshes.
\newblock {\em Computer Methods in Applied Mechanics and Engineering},
  282:132--160, 2014.

\bibitem{Paulino:2017:CMAME}
H.~Chi, L.~Beir\~{a}o~da Veiga, and G.~H. Paulino.
\newblock Some basic formulations of the virtual element method (vem) for
  finite deformations.
\newblock {\em Computer Methods in Applied Mechanics and Engineering},
  318:142--190, 2017.

\bibitem{Wriggers:2017:CM}
P.~Wriggers, B.~D. Reddy, W.~T. Rust, and B.~Hudobivnik.
\newblock Efficient virtual element formulations for compressible and
  incompressible finite deformations.
\newblock {\em Computational Mechanics}, 60:253--268, 2017.

\bibitem{Paulino:2021:MECC}
K.~Park, H.~Chi, and G.H. Paulino.
\newblock B-bar virtual element method for nearly incompressible and
  compressible materials.
\newblock {\em Meccanica}, 56:1423–--1439, 2012.

\bibitem{Beirao-Lovadina-Mora:2015}
L.~{Beir\~{a}o~da~Veiga}, C.~Lovadina, and D.~Mora.
\newblock A virtual element method for elastic and inelastic problems on
  polytope meshes.
\newblock {\em Computer Methods in Applied Mechanics and Engineering},
  295:327--346, 2015.

\bibitem{Hudobivnik:2019:CM}
B.~Hudobivnik, F.~Aldakheel, and P.~Wriggers.
\newblock A low order 3d virtual element formulation for finite
  elasto–plastic deformations.
\newblock {\em Computational Mechanics}, 63:253–--269, 2019.

\bibitem{Artioli:2017:CMAME}
E.~Artioli, S.~De~Miranda, C.~Lovadina, and L.~Patruno.
\newblock A stress/displacement virtual element method for plane elasticity
  problems.
\newblock {\em Computer Methods in Applied Mechanics and Engineering},
  325:155--174, 2017.

\bibitem{Dassi:2020}
F.~Dassi, C.~Lovadina, and M.~Visinoni.
\newblock A three-dimensional {H}ellinger–{R}eissner virtual element method
  for linear elasticity problems.
\newblock {\em Computer Methods in Applied Mechanics and Engineering},
  364:112910, 2020.

\bibitem{Artioli:2020:CMAME}
E.~Artioli, L.~Beir\~{a}o~da Veiga, and F.~Dassi.
\newblock Curvilinear virtual elements for 2d solid mechanics applications.
\newblock {\em Computer Methods in Applied Mechanics and Engineering},
  359:112667, 2020.

\bibitem{Paulino:2019:CMAME}
K.~Park, H.~Chi, and Paulino G.H.
\newblock On nonconvex meshes for elastodynamics using virtual element methods
  with explicit time integration.
\newblock {\em International Journal for Numerical Methods in Engineering},
  356:669--684, 2019.

\bibitem{Paulino:2020:IJNME}
K.~Park, H.~Chi, and Paulino G.H.
\newblock Numerical recipes for elastodynamic virtual element methods with
  explicit time integration.
\newblock {\em International Journal for Numerical Methods in Engineering},
  121:1--31, 2020.

\bibitem{Antonietti:2021:IJNME}
P.~F. Antonietti, G.~Manzini, I.~Mazzieri, H.~M. Mourad, and M.~Verani.
\newblock The arbitrary-order virtual element method for linear elastodynamics
  models: convergence, stability and dispersion-dissipation analysis.
\newblock {\em International Journal for Numerical Methods in Engineering},
  122:934--971, 2021.

\bibitem{Benedetto-Berrone-Pieraccini-Scialo:2014}
M.~F. Benedetto, S.~Berrone, S.~Pieraccini, and S.~Scial{\`o}.
\newblock The virtual element method for discrete fracture network simulations.
\newblock {\em Computer Methods in Applied Mechanics and Engineering}, 280:135
  -- 156, 2014.

\bibitem{Nguyen:2018:VEM}
V.~M. Nguyen-Thanh, X.~Zhuang, H.~{Ngyyen-Xuan}, T.~Rabczuk, and P.~Wriggers.
\newblock A {Virtual Element Method} for {2D} linear elastic fracture analysis.
\newblock {\em Computer Methods in Applied Mechanics and Engineering},
  340:366--395, 2018.

\bibitem{Hussein:2019}
A.~Hussein, F.~Aldakheel, B.~Hudobivnik, P.~Wrigger, P.A. Guidault, and
  O.~Allix.
\newblock A computational framework for brittle crack-propagation based on
  efficient virtual element method.
\newblock {\em Computer Methods in Applied Mechanics and Engineering},
  159:15--32, 2019.

\bibitem{Artioli:2020}
E.~Artioli, S.~Marfia, and E.~Sacco.
\newblock {VEM}-based tracking algorithm for cohesive/frictional 2d fracture.
\newblock {\em Computer Methods in Applied Mechanics and Engineering},
  365:112956, 2020.

\bibitem{Perugia-Pietra-Russo:2016}
I.~Perugia, P.~Pietra, and A.~Russo.
\newblock A plane wave virtual element method for the {H}elmholtz problem.
\newblock {\em ESAIM: Mathematical Modelling and Numerical Analysis},
  50(3):783--808, 2016.

\bibitem{Benvenuti:2019:CMAME}
E.~Benvenuti, A.~Chiozzi, G.~Manzini, and N.~Sukumar.
\newblock Extended virtual element method for the {Laplace} problem with
  singularities and discontinuities.
\newblock {\em Computer Methods in Applied Mechanics and Engineering},
  356:571--597, 2019.

\bibitem{Chiozzi:2020:MECC}
A.~Chiozzi and E.~Benvenuti.
\newblock Extended virtual element method for the torsion problem of cracked
  prismatic beams.
\newblock {\em Meccanica}, 55:637--648, 2020.

\bibitem{Artioli-Mascotto:2021:CMAME}
E.~Artioli and L.~Mascotto.
\newblock Enrichment of the nonconforming virtual element method with singular
  functions.
\newblock {\em Computer Methods in Applied Mechanics and Engineering},
  385:114024, 2021.

\bibitem{Duarte:2000:CS}
C.~A. Duarte, I.~Babu\v{s}ka, and J.~T. Oden.
\newblock Generalized finite element methods for three-dimensional structural
  mechanics problems.
\newblock {\em Computer and Structures}, 77:215--232, 2000.

\bibitem{Hansbo:2004:FEM}
A.~Hansbo and P.~Hansbo.
\newblock A finite element method for the simulation of strong and weak
  discontinuities in solid mechanics.
\newblock {\em Computer Methods in Applied Mechanics and Engineering},
  193(33-35):3523--3540, 2004.

\bibitem{Grisvard:1992}
P.~Grisvard.
\newblock {\em Singularities in Boundary Value Problems}.
\newblock Masson, Paris, France, 1992.

\bibitem{Areias:2006:ACO}
P.~M.~A. Areias and T.~Belytschko.
\newblock {A comment on the article `A finite element method for simulation of
  strong and weak discontinuities in solid mechanics' by A. Hansbo and P.
  Hansbo [Comput. Methods Appl. Mech. Engrg. 193 (2004) 3523-3540]}.
\newblock {\em Computer Methods in Applied Mechanics and Engineering},
  195:1275--1276, 2006.

\bibitem{Duchon:1977:SMR}
J.~Duchon.
\newblock Splines minimizing rotation-invariant semi-norms in {Sobolev} spaces.
\newblock In {\em Constructive Theory of Functions of Several Variables},
  volume 571 of {\em Lecture Notes in Mathematics}, pages 85--100.
  Springer-Verlag, Berlin, Germany, 1977.

\bibitem{Benvenuti:2021:CM}
E.~Benvenuti and N.~Orlando.
\newblock A mesh-independent framework for crack tracking in elastodamaging
  materials through the regularized extended finite element method.
\newblock {\em Computational Mechanics}, 68:25--49, 2021.

\bibitem{BeiraodaVeiga-Brezzi-Marini-Russo:2014}
L.~{Beir\~{a}o~da~Veiga}, F.~Brezzi, L.~D. Marini, and A.~Russo.
\newblock The hitchhiker's guide to the virtual element method.
\newblock {\em Mathematical Models \& Methods in Applied Sciences},
  24(8):1541--1573, 2014.

\bibitem{Beirao-Dassi-Russo:2017}
L.~{Beir\~{a}o~da~Veiga}, F.~Dassi, and A.~Russo.
\newblock High-order virtual element method on polyhedral meshes.
\newblock {\em Computers and Mathematics with Applications}, 74:1110--1122,
  2017.

\bibitem{Rice:1968}
J.~Rice.
\newblock A path independent integral and the approximate analysis of strain
  concentration by notches and cracks.
\newblock {\em Journal of Applied Mechanics}, 35:379--386, 1968.

\bibitem{Li:1985}
F.~Z. Li, C.~F. Shih, and A.~Needleman.
\newblock A comparison of methods for calculating energy release rates.
\newblock {\em Engineering Fracture Mechanics}, 21(2):405--421, 1985.

\bibitem{BeiraodaVeiga-Chernov-Mascotto-Russo:2018}
L.~{Beir\~{a}o~da~Veiga}, A.~Chernov, L.~Mascotto, and A.~Russo.
\newblock Exponential convergence of the hp virtual element method in presence
  of corner singularities.
\newblock {\em Numerische Mathematik}, 138:581--613, 2018.

\bibitem{Dolbow:2004:EEA}
J.~E. Dolbow and A.~Devan.
\newblock Enrichment of enhanced assumed strain approximations for representing
  strong discontinuities: addressing volumetric incompressibility and the
  discontinuous patch test.
\newblock {\em International Journal for Numerical Methods in Engineering},
  59(1):47--67, 2004.

\bibitem{Paulino:2007:SMO}
C.~Talischi, G.~H. Paulino, A.~Pereira, and F.~M. Menezes.
\newblock {PolyMesher: a general-purpose mesh generator for polygonal elements
  written in {M}atlab}.
\newblock {\em Structural and Multidisciplinary Optimization}, 45(3):309--328,
  2012.

\bibitem{Grisvard:1985}
P.~Grisvard.
\newblock {\em Elliptic Problems in Nonsmooth Domains}.
\newblock Pitman Publishing, Inc, Boston, MA, 1985.

\bibitem{laborde:2005}
P.~Laborde, J.~Pommier, Y.~Renard, and M.~Sala\"{u}n.
\newblock High-order extended finite element method for cracked domains.
\newblock {\em International Journal for Numerical Methods in Engineering},
  64(3):354--381, 2005.

\bibitem{Bechet:2005:IIR}
E.~B\'echet, H.~Minnebo, N.~Mo{\"e}s, and B.~Burgardt.
\newblock Improved implementation and robustness study of the {X-FEM} for
  stress analysis around cracks.
\newblock {\em International Journal for Numerical Methods in Engineering},
  64(8):1033--1056, 2005.

\end{thebibliography}
